\documentclass[review]{elsarticle}
\usepackage{mathtools}
\usepackage{graphicx}
\usepackage{epstopdf}
\usepackage{multirow}
\usepackage{csvsimple}
\usepackage{bm}
\usepackage{float}
\usepackage{amsfonts}
\usepackage[toc,page]{appendix}
\usepackage{amsmath}
\usepackage{caption}
\usepackage{xcolor}
\usepackage[labelformat=simple]{subcaption}

\usepackage{lineno}
\modulolinenumbers[5]

\journal{Journal of \LaTeX\ Templates}

\begin{document}
	
	\begin{frontmatter}
		
		\title{Higher-order accurate diffuse-domain methods for partial differential equations with Dirichlet boundary conditions in complex, evolving geometries}

		
		\author[mymainaddress]{Fei Yu}
		
		\author[mymainaddress]{Zhenlin Guo}
		
		\author[mymainaddress,mymainaddress1]{John Lowengrub}

		\address[mymainaddress]{Department of Mathematics, University of California, Irvine, CA, 92697, United States}
		\address[mymainaddress1]{Department of Mathematics and Department of Biomedical Engineering, University of California, Irvine, CA, 92697, United States}

		\begin{abstract}
			The diffuse-domain, or smoothed boundary, method is an attractive approach for solving partial differential equations in complex geometries because of its simplicity and flexibility. In this method the complex geometry is embedded into a larger, regular domain. The original PDE is reformulated using a smoothed characteristic function of the complex domain and source terms are introduced to approximate the boundary conditions. The reformulated equation, which is independent of the dimension and domain geometry, can be solved by standard numerical methods and the same solver can be used for any domain geometry. A challenge is making the method higher-order accurate. For Dirichlet boundary conditions, which we focus on here, current implementations demonstrate a wide range in their accuracy but generally the methods yield at best first order accuracy in $\epsilon$, the parameter that characterizes the width of the region over which the characteristic function is smoothed. Typically, $\epsilon\propto h$, the grid size. Here, we analyze the diffuse-domain PDEs using matched asymptotic expansions and explain the observed behaviors.  Our analysis also identifies simple modifications to the diffuse-domain PDEs that yield higher-order accuracy in $\epsilon$, e.g., $O(\epsilon^2)$ in the  \(L^2\) norm and $O(\epsilon^p)$ with $1.5\le p\le 2$ in the \(L^{\infty}\) norm. Our analytic results are confirmed numerically in stationary and moving domains where the level set method is used to capture the dynamics of the domain boundary and to construct the smoothed characteristic function. 
			  		\end{abstract}
		
		
	\end{frontmatter}
	
	\section{Introduction}
	Many problems in the physical, biological and engineering sciences involve solving partial differential equations (PDEs) in complex geometries that change their size and shape in time. There are two main approaches for solving such problems. In one approach, interface-fitted meshes are used. Examples include finite element methods (e.g., \cite{000395210500018}), boundary integral (e.g., \cite{Galenko-2018}) and boundary element (e.g., \cite{Feischl-2015}) methods. Because of the challenges associated with generating interface-fitted meshes for complex geometries, especially in three-dimensions, and because in many applications the complex geometry evolves in time, which would require a new mesh to be generated at each time step, embedded domain methods have been developed as an alternative approach. 
	
	In this second approach, the complex domain is embedded into a larger, regular domain and the boundary conditions are approximated by a variety of different techniques. Examples include the adaptive fast multipole accelerated Poisson solver (e.g., \cite{000402481300001}), which combines boundary and volume integral methods in the larger domain, fictitious domain methods (e.g., \cite{A1994MX36800005,000343821300007,000326596600021, 000336620900010}) where Lagrange multipliers are applied in order to enforce the boundary conditions, immersed boundary (e.g., \cite{000226822300011, 000330577100008, 000341311000010,000397072800006}), front-tracking (e.g., \cite{000169338900012,000398874700008,000412788000084}) and arbitrary Lagrangian-Eulerian methods (e.g., \cite{A1997XQ32700018,000381585100032,000360087400046,000414478700020}) utilize separate surface and volume meshes
 where force distributions are interpolated from the surface to the volume meshes, in a neighborhood of the domain boundary, to approximate the boundary conditions. {\color{black} In addition, a number of specialized methods have been designed to achieve better than first order accuracy in the $L^\infty$ norm. These include the immersed interface (e.g.,  \cite{ A1994PB02700004,000432512900032,000431399600012, 000432769500031}{\color{red}\cite{li2013adaptive,Kolahdouz2020}}), ghost fluid (e.g., \cite{000081366600002, 000362379300017,000375799200041,000418229800020}), cut-cell methods (e.g., \cite{000180706400035, 000432481000009,000428483000007,000427393800011}) and Voronoi interface (e.g., \cite{000306367500004,000358796700042}) methods. These methods modify difference stencils near the domain boundary to account for the boundary conditions. Further, the extended finite element method (e.g., \cite{000283202200001,000266373500001,000402217000004}) approximates the boundary conditions by enriching the space of test functions. A related approach is the virtual node method \cite{bedrossian2010second,HELLRUNG20122015}, which can achieve 2nd order accuracy in $L^\infty$ and involves adding (virtual) degrees of freedom along the interface together with Lagrange multipliers to enforce the boundary conditions. Another related approach is the active penalty method \cite{shirokoff2015sharp}, which can achieve higher order accuracy in the $L^\infty$ norm by extending the boundary function off the interface such that the extension matches the physical function and its normal derivatives at the interface. The smoothness of the solution and the number of matched normal derivatives constrains the accuracy of the method. {\color{red}  A disadvantage of most of the methods described above is that either modifications of standard finite element and finite difference software packages are needed, or extra auxiliary equations have to be solved.}}
	
Diffuse-domain methods (DDM), also known as smoothed boundary methods, have emerged as an attractive alternative approach because they are easy to implement and the formulation does not depend on the dimension of the problem or the geometry of the domain. In the DDM, the complex geometry is embedded into a larger, regular domain and a phase-field function is used to provide a smooth approximation to the characteristic function of the complex domain. {\color{red} A parameter $\epsilon$, usually proportional to the grid size, is introduced that characterizes the width of the diffuse interfacial region and typically controls the accuracy of the approximation. The original PDE is reformulated with additional source terms that enforce the boundary conditions. When $\epsilon$ is small, the DDM is most efficient when combined with adaptive mesh refinement to enable the use of small grid cells to resolve the narrow $(\approx O(\epsilon)$) transition layer and large grid cells in the extended, non-physical part of the regular domain. While DDMs are not always as accurate as some of the specialized methods described above and involve the solution of equations with non-constant coefficients and introduce additional length scales, DDMs have the advantage that they can be easily formulated for a wide range of equations and solved using standard uniform or adaptive discretizations combined with iterative matrix solvers, which are normally contained in standard PDE software packages. Further, the same solver can be used for any domain geometry. }



	
	The DDM was first used in \cite{Kockelkoren-2003,Bueno-Orovio-2005,Bueno-Orovio-2006} to solve diffusion equations with Neumann (no-flux) boundary conditions. The DDM was extended to simulate PDEs on surfaces in \cite{Ratz-2006}, to PDEs with Robin and Dirichlet boundary conditions in \cite{li09} and to cases in which bulk and surface equations are coupled \cite{tei09}. Later, in \cite{Yu-2012} and \cite{Poulson-2018} alternate derivations of diffuse-domain methods for such problems were presented. Diffuse-domain methods have been applied to a wide variety of problems that arise, for example, in biology (e.g., \cite{Kockelkoren-2003,Fenton-2005,Aland-2011a,Camley-2013,Chen-2014,Ratz-2014,Ratz-2015,Lowengrub-2016,Camley-2017,Lipkova2019}), in fluid dynamics (e.g., \cite{teigen2011diffuse,ala10,Aland-2011,Aland-2012,Aland-2014}) and in materials science (e.g., \cite{Yu-2016,Ratz-2015,Ratz-2016,Hong-2016,Chadwick-2018}), just to name a few, {\color{red} and have been implemented using a wide spectrum of methods including in-house finite-difference, finite-element and pseudo-spectral algorithms and software packages such as Matlab,  AMDiS \cite{Vey2007}, MRAG \cite{Rossinelli2015}, and BSAM \cite{Feng2018}.}
			
	Using rigorous mathematical theory \cite{fra12,Abels-2015,Schlottbom-2016,Burger-2017}, matched asymptotic expansions and numerical simulations (e.g., \cite{Ratz-2006,li09,tei09,Yu-2012,Poulson-2018}), the DDMs have been shown to converge to the original PDE and boundary conditions as the diffuse interface parameter \(\epsilon\) tends to zero. Further, in \cite{000356209300006} a matched asymptotic analysis for general DDMs with Neumann and Robin boundary conditions showed that for certain choices of the source terms, the DDMs were second-order accurate in $\epsilon$ and in the grid size $h$ {\color{black}in both the \(L^2\) and \(L^\infty\) norms}, taking $\epsilon\propto h$; see the recent paper \cite{Burger-2017} for a rigorous proof.
	
	For Dirichlet boundary conditions it is more challenging to obtain higher-order accurate diffuse domain methods. Current implementations of DDMs demonstrate a wide range in their accuracy but generally the methods yield at best first-order accuracy in $\epsilon$ \cite{fra12,Yu-2012,Schlottbom-2016,Burger-2017,Poulson-2018}. For example, some schemes, including one presented in \cite{li09}, achieve  sub-first order accuracy: \({O}(\epsilon\ln\epsilon)\) {\color{black}in both the \(L^2\) and the \(L^\infty\) norms}. While others, including some presented in \cite{li09,Yu-2012,Poulson-2018}, achieve first-order accuracy: ${O}(\epsilon)$ {\color{black}in both the \(L^2\) and \(L^\infty\) norms}. Further, in \cite{Poulson-2018} it was shown using asymptotic analysis that for several DDM approximations of the diffusion equation with Dirichlet boundary conditions, that while the order of accuracy was ${O}(\epsilon)$ {\color{black} in $L^2$ } the error could be significantly reduced by special choices of the DDM parameters. However, because DDMs for Dirichlet boundary conditions derive from introducing penalty terms in the equation, which are inherently low-order accurate, new modifications must be introduced to achieve higher-order accuracy. 

	In this paper, we analyze, using matched asymptotic expansions, several DDMs for approximating PDEs with Dirichlet boundary conditions. Our analysis explains the wide variation in the observed accuracy of the schemes and identifies simple modifications of the DDMs that can yield higher-order accuracy in $\epsilon$. {\color{black} Similar to the active penalty method \cite{shirokoff2015sharp}, our modifications are designed to enforce continuity of the normal derivative of the physical solution and its extension off the interface. However, we match the normal derivative only at the leading order in the asymptotic expansion in $\epsilon$, which simplifies the method.} Our new schemes are ${O}(\epsilon^2)$ in the  \(L^2\) norm and ${O}(\epsilon^p)$ in the \(L^{\infty}\) norm, with $p$ ranging from $1.5$ to $2$ depending on the range of $\epsilon$ and $h$ considered. For example, in the limit, $\epsilon\to 0$ and $\epsilon/h\to 0$ we obtain $p\to 1.5$. Surprisingly, the schemes can be more accurate when $\epsilon\propto h$. {\color{black} We find this can occur because the truncation/discretization error $u_{h,\epsilon}-u_{\epsilon}$ and the analytic error $u_\epsilon-u$ can have similar magnitudes but opposite signs so that there can be cancellation in the total error $u_{h,\epsilon}-u$. As $h$ decreases, the error eventually becomes dominated by the analytic error.}
	
	 Our formulation is dimension-independent, easy to implement and independent of computational methodology.  Here we use the level-set method \cite{Osher-1988,Gibou-2018} to generate the phase-field, or smoothed characteristic, function, and adaptive finite-difference discretizations with a mass-conserving multigrid method on block-structured adaptive meshes (BSAM) \cite{Feng2018} to obtain highly-efficient schemes. We present examples in $1D$, $2D$ and $3D$ in stationary and moving domains and confirm our analytic results.

	
	The outline of the paper is as follows. In Sec. \ref{sec_ddm}, we analyze DDM approximations to the Poisson equation with Dirichlet boundary conditions and present 1D numerical results that confirm our analysis. In Sec. \ref{heat1d}, we extend our analysis and simulations to time-dependent equations, focusing on the heat equation with external forcing. We consider both stationary and moving domains. We discuss different approaches to model the moving domains (phase-field versus level-set) and present 1D numerical simulations using the  method. {\color{black}2D and 3D numerical results are provided in Secs. \ref{2d} and  \ref{3d}, respectively.} Finally, in Sec. \ref{conclude}, we present conclusions and discuss future work. In Appendixes \ref{app_a} and \ref{app_b}, we provide additional numerical results and theoretical details.
	
	\section{Diffuse domain method to the Poisson equation with Dirichlet boundary conditions}\label{sec_ddm}
	Consider the Poisson equation with Dirichlet boundary condition on a complex (non-regular) domain \(D\):
	\begin{eqnarray}
	\Delta u = f \ \ in \ \ D,\label{p1} \\
	u=g  \ \ on\ \ \partial D.\label{p2}
	\end{eqnarray}
	To approximate Eqs. \eqref{p1} and \eqref{p2}, we consider three DDM approximations defined in \cite{li09} and posed
	on a larger, regular domain \(\Omega\):
	\begin{eqnarray}
	&\textbf{DDM1}: &\nabla\cdot(\phi \nabla u_\epsilon)-\frac{1}{\epsilon^3}(1-\phi)(u_{\epsilon}-g)=\phi f , \label{eq_ddm1}\\
	&\textbf{DDM2}:&\phi \Delta u_\epsilon-\frac{1}{\epsilon^2}(1-\phi)(u_{\epsilon}-g)=\phi f, \label{eq_ddm2}\\
	&\textbf{DDM3}:&\nabla\cdot(\phi \nabla u_\epsilon)-\frac{|\nabla\phi|}{\epsilon^2}(u_{\epsilon}-g)=\phi f, \label{eq_ddm3}
	\end{eqnarray}
	where the source terms \(\frac{1}{\epsilon^3}(1-\phi)(u_{\epsilon}-g)\), \(\frac{1}{\epsilon^2}(1-\phi)(u_{\epsilon}-g)\) and \(\frac{|\nabla\phi|}{\epsilon^2}(u_{\epsilon}-g)\) represent different choices for enforcing the Dirichlet boundary condition \(u=g\). The function \(\phi\) approximates the characteristic function of \(D\), e.g.
	\begin{equation}
	\phi(\bm{x},t)=\frac{1}{2}(1-\tanh{(\frac{3r(\bm{x},t)}{\epsilon})}),\label{phi}
	\end{equation}
	here \(r(\bm{x},t)\) is the signed-distance function to \(\partial D\), which is taken to be negative inside \(D\), and
	 \(\epsilon\) is the interface thickness. 
	{\color{black} The interface thickness $\epsilon$ should be smaller than any physically relevant length scale. In particular, $\epsilon<<1/|\kappa|$ where $\kappa$ is the total curvature of the interface and $\epsilon<||u||_{\partial D,\infty}/||\nabla u||_{\partial D,\infty}$ where the norms denote the $L^\infty$ norm on the boundary $\partial D$. In addition, we require $\epsilon<\delta$, where $\delta=\delta(\kappa)$ is the width of a layer around the boundary where the functions $g$ and $f$ can be smoothly extended outside $D$ constant in the normal direction. } {\color{blue} The distance over which $g$ and $f$ need to be extended, and the distance between $\partial D$ and $\partial\Omega$, depends on the speed at which the approximate characteristic function $\phi$ asymptotes. Generally, this distance should be at least $3\epsilon$ so that $\phi$ is sufficiently close to zero. Accordingly, we take $\delta=3\epsilon$ and find no significant improvement in the errors when larger values of $\delta$ are used. To be safe, we took  the distance between $\partial D$ and $\partial\Omega$ to be $\approx 10\epsilon$ and find this works well for all the numerical simulations in this paper, although we likely could have used a smaller value.}
		
	Note that for DDM3, since the phase field function \(\phi\) and its gradient vanish rapidly outside \(D\), in order to prevent the equation from being ill-posed, we use the following modified gradient instead in the numerical calculation \cite{000356209300006},
	\begin{align}
	|\nabla\hat{\phi}|&=\tau+(1-\tau)|\nabla\phi|\label{mgradphi},
	\end{align}
	where \(\tau=10^{-15}\). 
Note that no such regularization is required for DDMs 1 and 2 because the {\color{black} $1-\phi\approx 1$} terms dominate outside $D$. In all the DDMs we consider, we assume that the boundary function \(g\) is defined in a small \({O}(\epsilon)\) neighborhood of \(\partial D\) using an extension off \(\partial D\) that is constant in the normal direction and that \(f\) is similarly extended smoothly out of the domain \(D\). {\color{black} This can always be done if $\epsilon$ is sufficiently small as described above.}

\subsection{Numerical results for DDM1-3}\label{num_1d_ddm}
	We present representative numerical results in 1D for the DDMs given in Sec. \ref{sec_ddm}. We suppose the original Poisson equation is defined in \(D=[-1.111,1.111]\) with Dirichlet boundary conditions at \(x=\pm1.111\), which are not grid points. We take the forcing function $f=1$ and boundary condition $g=1.111^2/2$ so that the exact solution to Eqs. (\ref{p1})-(\ref{p2}) is $u={x^2}/{2}$. The results from other test cases can be found in Appendix \ref{app_a}.
	
	We solve the DDMs in a larger domain \(\Omega=[-2,2]\) with the same Dirichlet boundary condition as the original Poisson equation. Note that the choice for boundary condition for the outer domain \(\Omega\) should not affect the inner domain \(D\). Therefore the distance between \(\partial D\) and \(\partial \Omega\) should be large compared to \(\epsilon\).
	The equations are discretized on a uniform grid with the standard second-order central difference scheme (adaptive mesh refinement is used later in 2D, see Sec. \ref{2d}). The discrete system is solved using the Thomas method \cite{hoffman2001numerical}. 
	
\begin{figure}[H]
		\begin{subfigure}{.49\textwidth}
			\includegraphics[width=\linewidth]{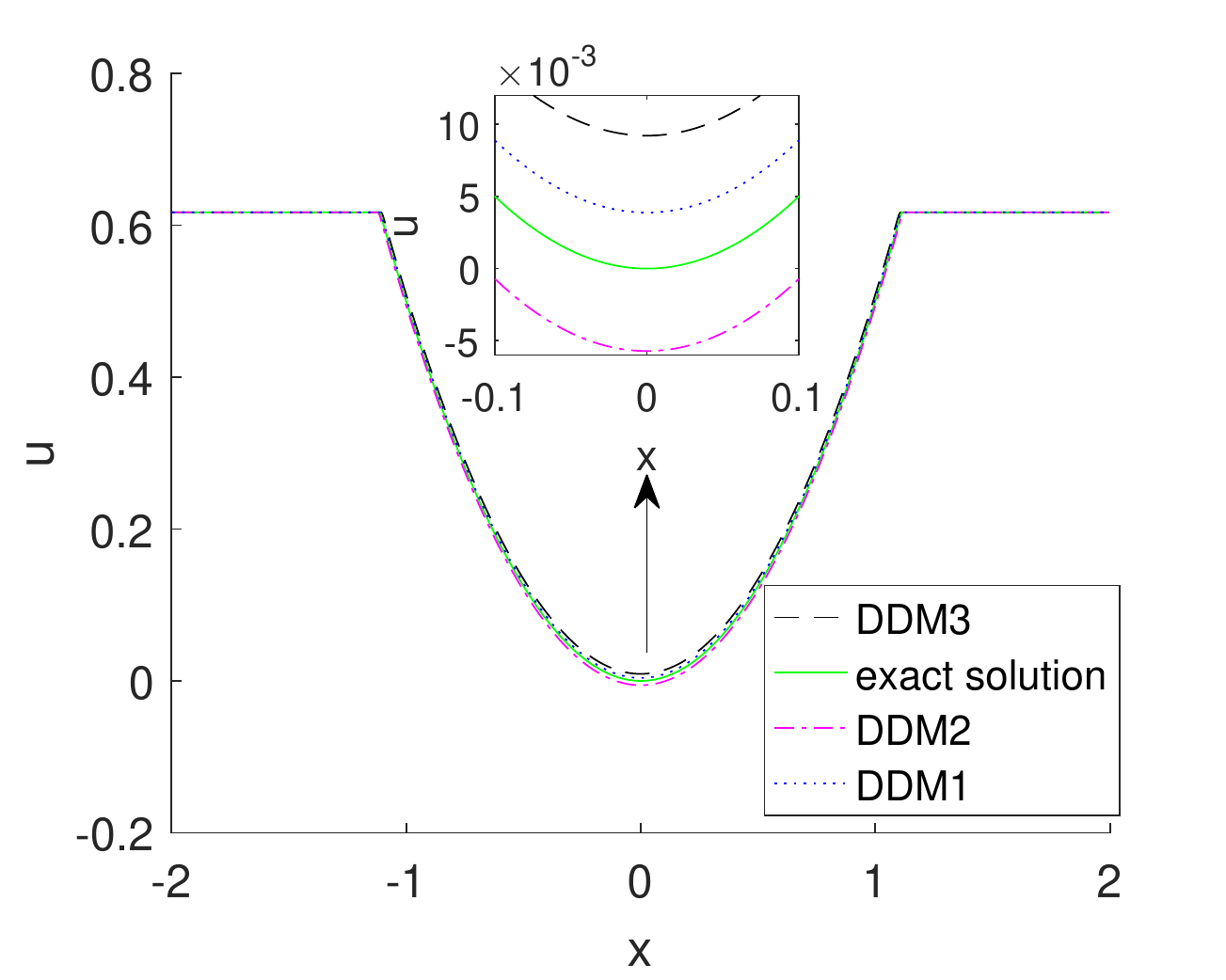}
			\captionsetup{width=0.8\textwidth}
			\caption{}\label{ddm_all}
		\end{subfigure}
		\begin{subfigure}{.49\textwidth}
			\includegraphics[width=\linewidth]{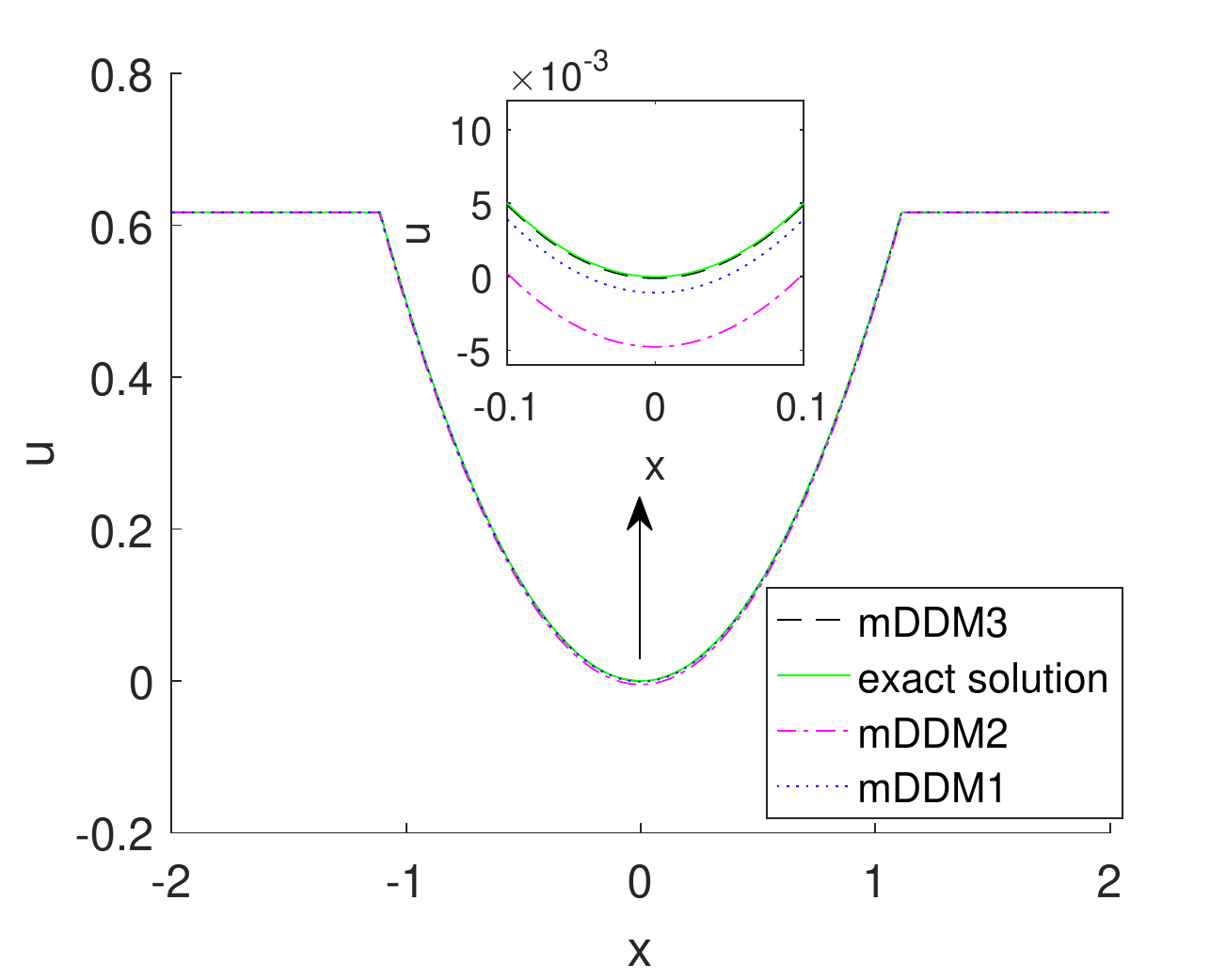}
			\captionsetup{width=0.8\textwidth}
			\caption{}\label{mddm_all}
		\end{subfigure}
		\caption{(a). Numerical solutions of DDM1--3 from Eqs. (\ref{eq_ddm1})-(\ref{eq_ddm3}); (b). The corresponding modified diffuse domain methods mDDM1-3 from Eqs. (\ref{eq_mddm1})-(\ref{eq_mddm3}); with $f=1$, $g=1.111^2/2$, $\epsilon=0.0125$ and $h=\epsilon/4$.}
		\label{ddm_mddm_all}
	\end{figure}

	The numerical solutions of DDM1-3, together with the exact solution, are shown in Fig. \ref{ddm_all}, a close-up of the solutions with \(x\in[-0.1,0.1]\) is shown as an inset. Here, we have taken \(\epsilon=0.0125\) and \(h=\epsilon/4\). 
	%
	%
	We next calculate and compare the error between the simulated
	DDM solutions \(u_{h,\epsilon}\) and the analytic solution \(u\) of the original PDE. We present the errors in both the \(L^2\) and the \(L^{\infty}\) norms, defined as
	\begin{eqnarray}
	E^{(2)}_{\epsilon}=\frac{\lVert \phi(u-u_{h,\epsilon})\rVert_{L^2(\Omega)}}{\lVert\phi u\rVert_{L^2(\Omega)}}, \label{err}
	\end{eqnarray}
	where \(\lVert\cdot \rVert_{L^2}\) is the discrete \(L^2\) norm:
	\begin{equation}
	\lVert u \rVert_{L^2}={\sqrt {1/N\sum_{i=1}^N u_i^2}} \label{l2norm},
	\end{equation}
	where N is the number of grid points in \(\Omega\).
	The error in the \(L^{\infty}\) norm is defined as \cite{A1994PB02700004}
	\begin{equation}
	E^{(\infty)}_{\epsilon}=\frac{\lVert (u-u_{h,\epsilon})\rVert_{L^{\infty}(D)}}{\lVert u\rVert_{L^{\infty}(D)}}, \label{err_max}
	\end{equation}
	and \(\lVert u \rVert_{L^{\infty}(D)}= \max_{1\leq i \leq M} |u_i|\) and M is the number of grid points in \(D\) . Note this is not \(L^{\infty}\) in \(\Omega\).
	The convergence rate in \(\epsilon\) is calculated as 
	\begin{equation}
	\bm{k}=\frac{\log(\frac{E_{\epsilon_i}}{E_{\epsilon_{i-1}}})}{\log(\frac{\epsilon_i}{\epsilon_{i-1}})}\label{conv_order}.
	\end{equation}
	The results are presented in Tabs. \ref{ddm1_case12} and \ref{ddm2_case12}. We observe that DDM2 is 1st order accurate, which is consistent with findings in \cite{Yu-2012,Poulson-2018}, while the convergence rates of DDM1 and DDM3 are less than 1 in both \(L^2\) and \(L^{\infty}\) norms (see also Tabs. \ref{afd1}-\ref{afd6} in Appendix \ref{app_a} for other test cases).
	
\begin{figure}[H]
		\center
		\begin{subfigure}{.6\textwidth}
			\includegraphics[width=\linewidth]{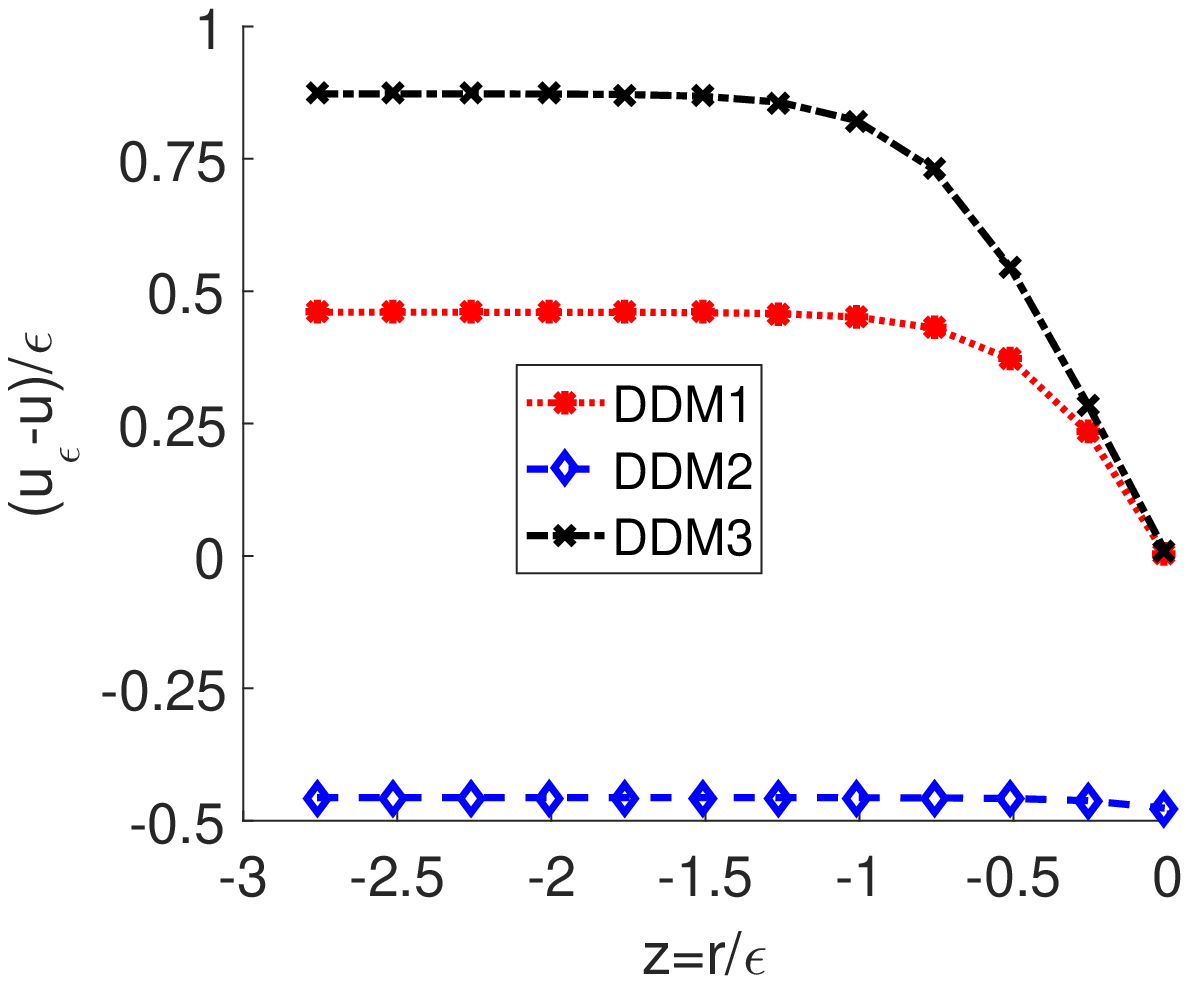}
			\captionsetup{width=0.8\textwidth}
			\caption{}\label{u1_1}
		\end{subfigure}
		\begin{subfigure}{.48\textwidth}
			\includegraphics[width=\linewidth]{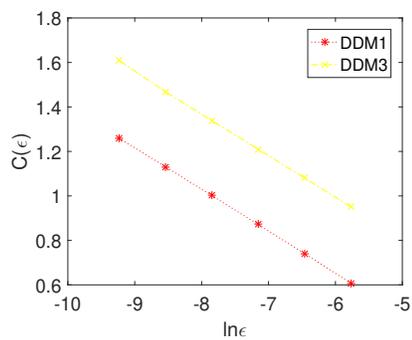}
			\captionsetup{width=0.8\textwidth}
			\caption{}\label{u1_3}
		\end{subfigure}
		\begin{subfigure}{.5\textwidth}
			\includegraphics[width=\linewidth]{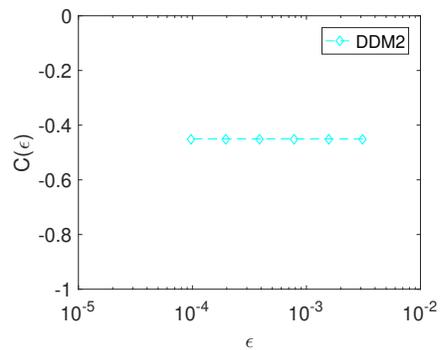}
			\captionsetup{width=0.8\textwidth}
			\caption{}\label{u1_2}
		\end{subfigure}
		
		\caption{Comparisons of the behavior of DDM1, DDM2 and DDM3 near $\partial D$. (a): \(\frac{u_\epsilon-u}{\epsilon}\) with \(\epsilon=0.00625\), (b): \(C(\epsilon)\) for DDM1 and DDM3 (see text), (c): \(C(\epsilon)\) for DDM2.}\label{u1_total}
	\end{figure}

	\begin{table}[H]
		\footnotesize 
		\center
		\begin{tabular}{c||c|c|c|c||c|c|c|c}
			&\multicolumn{4}{c||}{DDM1}&\multicolumn{4}{c}{DDM3}\\ \hline
			\(\epsilon\)&\(E^{(2)}\)&\(\bm{k}\)&\(E^{(\infty)}\)&\(\bm{k}\)&\(E^{(2)}\)&\(\bm{k}\)&\(E^{(\infty)}\)&\(\bm{k}\)
			\begin{filecontents*}{ddm13_case1.csv}
1.10E+01,,,,,,,,
2.00E-01,4.88E-01,0.00 ,2.12E-01,0.00 ,3.82E-02,0.00 ,4.73E-02,0.00 
1.00E-01,1.12E-01,2.12 ,5.26E-02,2.01 ,9.45E-02,-1.31 ,4.42E-02,0.10 
5.00E-02,9.88E-03,3.51 ,7.96E-03,2.72 ,8.03E-02,0.24 ,3.66E-02,0.27 
2.50E-02,1.31E-02,-0.41 ,5.89E-03,0.44 ,5.39E-02,0.58 ,2.44E-02,0.59 
1.25E-02,1.41E-02,-0.11 ,6.31E-03,-0.10 ,3.33E-02,0.69 ,1.50E-02,0.70 
6.25E-03,1.04E-02,0.44 ,4.66E-03,0.44 ,1.97E-02,0.76 ,8.84E-03,0.76 
3.13E-03,6.80E-03,0.62 ,3.04E-03,0.62 ,1.14E-02,0.79 ,5.09E-03,0.80 
1.56E-03,4.16E-03,0.71 ,1.86E-03,0.71 ,6.42E-03,0.82 ,2.87E-03,0.82 
\end{filecontents*}
\csvreader[]{ddm13_case1.csv}{}
			{\\\hline \csvcoli & \csvcolii & \csvcoliii & \csvcoliv & \csvcolv & \csvcolvi & \csvcolvii & \csvcolviii & \csvcolix}
		\end{tabular}
		\caption{The \(L^2\) and \(L^{\infty}\) errors in DDM1 and DDM3 with \(h=\epsilon/4\).}\label{ddm1_case12}
	\end{table}
	\begin{table}[H]
		\footnotesize 
		\center
		\begin{tabular}{c||c|c|c|c}
			&\multicolumn{4}{c}{DDM2}\\ \hline
			\(\epsilon\)&\(E^{(2)}\)&\(\bm{k}\)&\(E^{(\infty)}\)&\(\bm{k}\)
			\begin{filecontents*}{ddm2_case1.csv}
1.10E+01,,,,
2.00E-01,3.42E-01,0.00 ,1.65E-01,0.00 
1.00E-01,1.70E-01,1.01 ,7.95E-02,1.05 
5.00E-02,8.45E-02,1.01 ,3.88E-02,1.04 
2.50E-02,4.17E-02,1.02 ,1.91E-02,1.02 
1.25E-02,2.07E-02,1.01 ,9.51E-03,1.00 
6.25E-03,1.03E-02,1.00 ,4.82E-03,0.98 
3.13E-03,5.16E-03,1.00 ,2.40E-03,1.00 
1.56E-03,2.58E-03,1.00 ,1.20E-03,1.01 
\end{filecontents*}
\csvreader[]{ddm2_case1.csv}{}
			{\\\hline \csvcoli & \csvcolii & \csvcoliii & \csvcoliv & \csvcolv}
		\end{tabular}
		\caption{The \(L^2\) and \(L^{\infty}\) errors in DDM2 with \(h=\epsilon/4\).}\label{ddm2_case12}
	\end{table}
	To understand this behavior, the scaled errors near the boundary \(x=1.111\) are shown in Fig. \ref{u1_total}.
	In Fig. \ref{u1_1}, we plot \(\frac{u_\epsilon-u}{\epsilon}\) versus the stretched inner variable \(z=\frac{r}{\epsilon}\) (\(r = x-1.111\)) using \(\epsilon=0.00625\) and \(N=2560\), which is obtained from the grid size \(h=\frac{\epsilon}{4}\). When \(z<<0\) (inside the domain), \(\frac{u_\epsilon-u}{\epsilon}\) tends to a constant, \(C(\epsilon)\). 
	In Fig. \ref{u1_3}, the constants \(C(\epsilon)\) from DDM1 and DDM3 are plotted as a function of \(\ln(\epsilon)\). Linear fits give that \(C(\epsilon)\approx -0.186\ln(\epsilon)-0.475\) for DDM1 and \(C(\epsilon)\approx -0.188\ln(\epsilon)-0.143\) for DDM3, which suggests that the orders of accuracy of DDM1 and DDM3 are \(\epsilon\ln(\epsilon)\) as \(\epsilon\rightarrow 0\). 
	In Fig. \ref{u1_2}, we plot the limit value \(C(\epsilon)\) from DDM2 as a function of \(\epsilon\). The plots suggest that \(C(\epsilon)\rightarrow -0.45\neq 0\), which implies DDM2 is 1st order accurate.
	\subsection{Matched asymptotic analysis of DDMs}\label{sec_ana}
	In this section, we perform a matched asymptotic analysis in 1D for DDM1-3 in order to explain the numerical results in Sec. \ref{num_1d_ddm}. This analysis can be easily extended to higher dimensions. Without loss of generality, we also assume the boundary is located at \(x=0\) and we take the signed distance (near the boundary) to be \(r=x\), again $x<0$ denotes the interior of $D$. We consider the expansion of diffuse-domain variables in powers of the interface thickness \(\epsilon\) in regions close to and far from the interface, which are known as inner (\(\hat u\)) and outer (\(\bar u\)) expansions, respectively.  On each side of the interface, there exists an outer expansion, here labeled \(\bar u_1(x;\epsilon)\) where \(x<0\) and \(\phi=1\), and \(\bar u_2(x;\epsilon)\) where \(x>0\) and \(\phi=0\) (see Fig. \ref{k1} for an illustration).
	Clearly for all three DDMs, we have
	\begin{equation}
	\bar u_2=g.\label {bar_u_2}
	\end{equation} 
	Note that there could be multiple layers near the boundary and hence multiple inner expansions may be required to match the two outer expansions.
	\subsubsection{Analysis of DDM2} \label{sec_ddm2}
	We first present the matched asymptotic analysis of DDM2. In 1D, DDM2 is,
	\begin{equation}
	\phi\frac{d^2}{dx^2}u_\epsilon-\frac{1}{\epsilon^2}(1-\phi)(u_\epsilon-g)=\phi f \label{ddm20},
	\end{equation}
	We assume the outer expansion \(\bar u_1\) satisfies, 
	\begin{equation}
	\bar u_1(x;\epsilon)=\bar u_1^{(0)}(x)+\epsilon \bar u_1^{(1)}(x)+...\label{outer_ddm2}
	\end{equation}
	Plugging into Eq. \eqref{ddm20} and assuming that neither \(f\) nor \(g\) depends on \(\epsilon\), we have,
	\begin{align}
	&\frac{d^2}{dx^2} \bar u_1^{(0)} = f, \label{out1}\\
	&\frac{d^2}{dx^2} \bar u_1^{(k)} = 0, k=1,2,... \label{out2}
	\end{align}
	Now, if \(\bar u_1^{(0)}=g\) on \(\partial D\) so that \(\bar u_1^{(0)}\) is the unique solution to Eqs. \eqref{p1} and \eqref{p2}, then DDM2 reduces to the Poisson equation with Dirichlet boundary condition at leading order. The \(L^2\) convergence rate is then determined by the next leading order term, \(\bar u_1^{(1)}\). 
	In order to determine the boundary conditions for \(\bar u_1^{(0)}\) and \(\bar u_1^{(1)}\), we need to analyze the behavior of DDM2 (Eq. \eqref{ddm20}) near \(\partial D\) using inner expansions and then match the inner and outer expansions in a region of overlap.
	
\begin{figure}
		\center
		\includegraphics[width=0.6\textwidth]{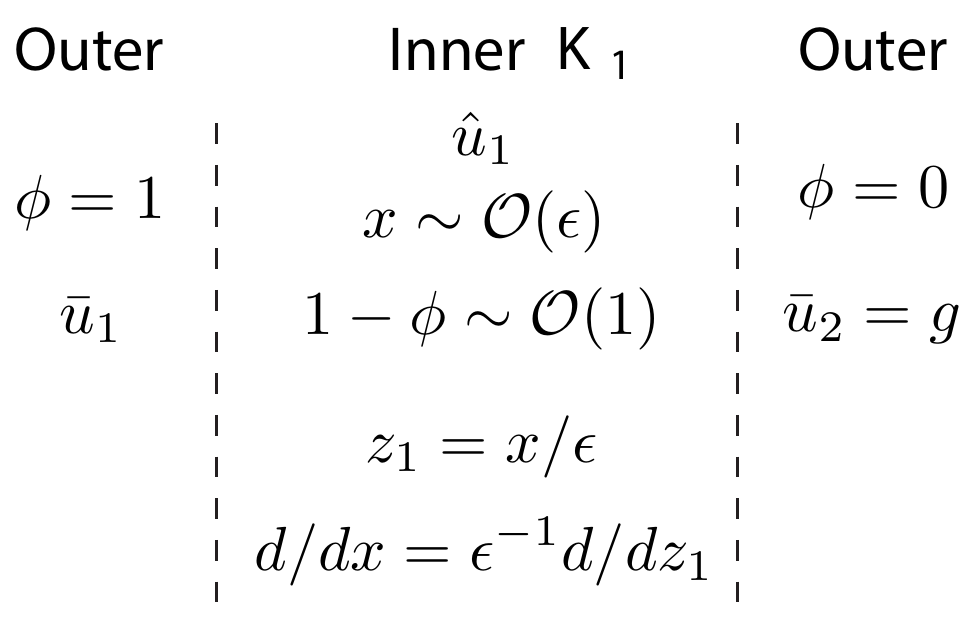}
		\caption{Inner layer \(K_1\)}
		\label{k1}
	\end{figure}

	To get the inner equations, we first rewrite DDM2 by substituting \(\phi\) with Eq. \eqref{phi} and dividing by \(\phi\):
	\begin{equation}
	\frac{d^2}{dx^2}u_{\epsilon}-\frac{e^{6x/\epsilon}}{\epsilon^2}(u_{\epsilon}-g) = f. \label{ddm2}
	\end{equation}
	Consider Eq. \eqref{ddm2} in a region \(K_1\) (see Fig. \ref{k1}), where \(x\sim{O}(\epsilon)\) and \(e^{6x/\epsilon}\sim {O}(1)\) and we introduce a stretched variable,
	\begin{equation}
	z_1=\frac{x}{\epsilon}.
	\end{equation}
	In a local coordinate system near $\partial D$, the derivatives become
	\begin{eqnarray}
	\frac{d}{dx} = \frac{1}{\epsilon}\frac{d}{dz_1},\label{loc_d1}\\
	\frac{d^2}{dx^2} = \frac{1}{\epsilon^2}\frac{d^2}{dz_1^2},\label{loc_d2}
	\end{eqnarray} 
	The inner expansion associated with the stretched variable \(z_1\) is 
	\begin{equation}
	\hat u_1(z_1;\epsilon)=\hat u^{(0)}_1(z_1)+\epsilon \hat u_1^{(1)}(z_1)+...\label{inner_ddm2}
	\end{equation}
	To obtain the matching conditions for each outer solution \(\bar u_i\) (\(i=1,2\)), we assume that there is a region of overlap
	where both the inner and the outer expansions are valid. In this region, if we evaluate the outer expansion in the inner coordinates, this must
	match the limits of the inner solutions away from the interface, that is,
	\begin{align}
	&\hat{u}_1(z_1;\epsilon) \simeq  \bar u_1(\epsilon z_1;\epsilon), \text{ as } z_1\rightarrow-\infty~~{\rm and} ~~\epsilon z_1\to 0^-,\label{match1}\\
	& \hat{u}_1(z_1;\epsilon)  \simeq \bar u_2(\epsilon z_1;\epsilon)=g, \text{ as } z_1\rightarrow+\infty~~{\rm and} ~~\epsilon z_1\to 0^+. \label{match2}
	\end{align}
	Here a single inner expansion \(\hat{u}_1\) is able to match both outer expansions up to \({O}(\epsilon)\) for DDM2 as shown below. As we see later in Sec. \ref{ddm13}, this is not the case for DDM1 and DDM3 and an additional layer needs to be introduced. 
	Combining Eqs. (\ref{bar_u_2}), (\ref{outer_ddm2}), (\ref{inner_ddm2})-(\ref{match2}), we have the following asymptotic matching conditions up to \({O}(\epsilon)\),
	\begin{align}
	&\hat u_1^{(0)}(z_1)=\bar u_1^{(0)}(0), \text{ as } z_1\rightarrow-\infty,\label{mc1}\\
	&\hat u_1^{(1)}(z_1)=\bar u_1^{(1)}(0)+z_1\frac{d}{dx} \bar u_1^{(0)}(0), \text{ as } z_1\rightarrow-\infty,\label{mc2}
	\end{align}
	and
	\begin{align}
	&\hat u_1^{(0)}(z_1)=g, \text{ as } z_1\rightarrow+\infty, \label{mc3}\\
	&\hat u_1^{(1)}(z_1)=0, \text{ as } z_1\rightarrow+\infty. \label{mc4}
	\end{align}
		
	Plugging Eq. \eqref{inner_ddm2} into Eq. \eqref{ddm2} and using derivatives in the local coordinate system (Eqs. \eqref{loc_d1} and \eqref{loc_d2}), we obtain the following equations for the inner expansion \(\hat{u}_1\):
	\begin{align}
	&\text{At }{O}(\epsilon^{-2})\text{, }&\frac{d^2}{dz_1^2}\hat u_1^{(0)}-e^{6z_1}(\hat u_1^{(0)}-g)=0,\label{ddm2_in_eq1}\\
	&\text{At }{O}(\epsilon^{-1})\text{, }&\frac{d^2}{dz_1^2}\hat u_1^{(1)}-e^{6z_1}\hat u_1^{(1)}=0.\label{ddm2_in_eq2}
	\end{align}
	Note that the solution to the following homogeneous ordinary differential equation,
	\begin{equation}
	y''-e^{6x}y=0,\label{homo_bassel}
	\end{equation}
	is given by
	\begin{equation}
	y = C_1 I_0(\frac{e^{3x}}{3}) + C_2 K_0(\frac{e^{3x}}{3}), \label{mod_bessel0}
	\end{equation}
	where \(I_0\) and \(K_0\) are the modified Bessel functions of the first and second kind, respectively, and \(C_1\) and \(C_2\) are constants. Clearly, \(\hat u_1^{(0)}=g\) is a solution to Eq. \eqref{ddm2_in_eq1}. Therefore, 
	\begin{align}
	&\hat u_1^{(0)}(z_1) =g+ C_1 I_0(\frac{e^{3z_1}}{3}) + C_2 K_0(\frac{e^{3z_1}}{3}), \label{ddm2_in_sol1}\\
	&\hat u_1^{(1)}(z_1) =C_3 I_0(\frac{e^{3z_1}}{3}) + C_4 K_0(\frac{e^{3z_1}}{3}),\label{ddm2_in_sol2}
	\end{align}
	where \(C_i\)'s are constants. The modified Bessel functions satisfy:
	\begin{align}
	&\lim_{z_1\rightarrow -\infty}I_0(\frac{e^{3z_1}}{3}) = 1,\label{mb1}\\
	&\lim_{z_1\rightarrow +\infty}I_0(\frac{e^{3z_1}}{3}) = +\infty,\label{mb2}\\
	&\lim_{z_1\rightarrow -\infty}K_0(\frac{e^{3z_1}}{3}) \sim -\ln(\frac{1}{2}\frac{e^{3z_1}}{3})-\gamma\sim -3z_1+\ln 6-\gamma,\label{mb3}\\
	&\lim_{z_1\rightarrow +\infty}K_0(\frac{e^{3z_1}}{3}) = 0,\label{mb4}
	\end{align}
	where \(\gamma\approx 0.5772\) is the Euler\(-\)Mascheroni constant. From Eq. \eqref{mb3} and the matching condition (\ref{mc1}) we conclude that \(C_2=0\). Using Eq. \eqref{mb2} and the matching conditions as \(z_1\rightarrow+\infty\) (Eqs. \eqref{mc3} and \eqref{mc4}), we have \(C_1=C_3=0\). Putting everything together and using \eqref{mc2}, we find
	\begin{align}
	&\hat u_1^{(0)}=g,\label{ddm2_in_sol11}\\
	&\hat u_1^{(1)}=-\frac{A}{3}K_0(\frac{e^{3z_1}}{3})\sim Az_1+\frac{A}{3}(-\ln 6+\gamma), \text{ as } z_1\rightarrow-\infty, \label{ddm2_in_sol22}
	\end{align}
	where \(A=\frac{d}{dx}\bar u^{(0)}(0)\) is the derivative of the exact solution at the boundary. By the matching conditions (Eqs. \eqref{mc1} and \eqref{mc2}), the two leading terms of the outer solution \(\bar u_1\) satisfies
	\begin{align}
	&\bar u_1^{(0)}(0)=g,\label{ddm2_out_sol1}\\
	&\bar u_1^{(1)}(0)=\frac{A}{3}(-\ln 6 + \gamma). \label{ddm2_out_sol2}
	\end{align}
	Therefore, \(\bar u_1^{(0)}=u\), which is the exact solution of the Poisson equation with Dirichlet boundary condition (Eqs. \eqref{p1} and \eqref{p2}), and DDM2 is first order accurate in \(\epsilon\) in the \(L^2\) norm since \(\bar u_1-u\sim {O}(\epsilon)\) if \(A\neq 0\). As for the \(L^{\infty}\) norm, we also need to consider the error at the boundary , e.g.,  $\hat u_1(0)-g\approx\epsilon \hat  u_1^{(1)}(0)\sim {O}(\epsilon)$ for DDM2. Thus, DDM2 is first order accurate in \(\epsilon\) in the \(L^{\infty}\) norm as well if \(A\neq 0\). If \(A=0\), then the asymptotic analysis suggests better than first-order accuracy. In fact when $A=0$, DDM2 is ${O}\left(\epsilon^2\right)$ accurate in both \(L^2\) and \(L^\infty\) as shown in Appendix \ref{app_a}.
	
	To confirm our analysis, we compute \(\bar u_1^{(1)}(0)\approx C(\epsilon)\), where $C(\epsilon)\approx(u_\epsilon-u)/\epsilon$ as $z_1<<0$ is shown in Fig. \ref{u1_2}. Using \(A=1.111\), both the asymptotic analysis and the numerical results agree and give  \(\bar u_1^{(1)} \approx -0.450\) (see Tab. \ref{asy_th2} in Appendix \ref{app_c} for other cases).

	\subsubsection{Analysis of DDM1 and DDM3}\label{ddm13}
	We now extend the matched asymptotic analysis in Sec.\ref{sec_ddm2} to DDM1 and DDM3. In 1D, DDM1 reads,
	\begin{equation}
	\frac{d}{dx}(\phi \frac{d}{dx} u_\epsilon)-\frac{1}{\epsilon^3}(1-\phi)(u_\epsilon-g)=\phi f \label{ddm10},
	\end{equation}
	Substituting \(\phi\) from Eq. \eqref{phi}, we have,
	\begin{equation}
	\frac{d^2}{dx^2}u_{\epsilon} - \frac{6}{\epsilon}\frac{1}{1+e^{-6x/\epsilon}}\frac{d}{dx}u_{\epsilon}-\frac{e^{6x/\epsilon}}{\epsilon^3}(u_{\epsilon}-g) = f. \label{ddm1}
	\end{equation}
	We assume that
	\begin{equation}
	\bar u_1(x;\epsilon)=\bar u_1^{(0)}(x)+\epsilon \bar u_1^{(1)}(x;\ln \epsilon)+o(\epsilon). \label{outer1}
	\end{equation}
	Note that \(\ln \epsilon\) is included in the outer expansion, which we will explain later. Plugging Eq. \eqref{outer1} into the DDM1 (Eq. \eqref{ddm10}), again we have,
	\begin{align}
	&\frac{d^2}{dx^2} \bar u_1^{(0)} = f, \label{out1_1}\\
	&\frac{d^2}{dx^2} \bar u_1^{(k)} = 0, k=1,2,... \label{out2_1}
	\end{align}
	
\begin{figure}
		\center
		\includegraphics[width=0.8\textwidth]{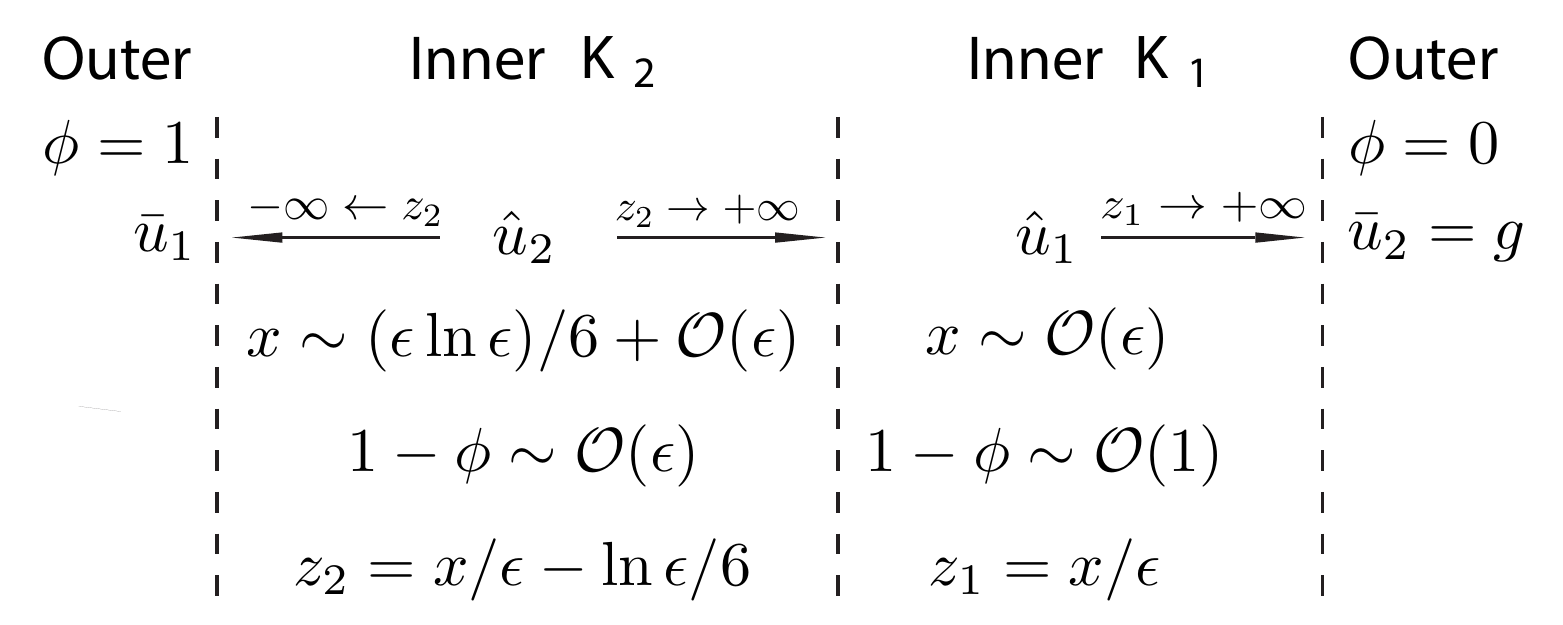}
		\caption{Inner layers \(K_1\) and \(K_2\)}\label{k12}
	\end{figure}

	As for inner expansion, we first consider DDM1 (Eq. \eqref{ddm1}) in \(K_1\), where \(x\sim {O}(\epsilon)\) and \(e^{6x/\epsilon}\sim {O}(1)\). We introduce the same stretched variable,
	\begin{equation}
	z_1=\frac{x}{\epsilon},\label{z1}
	\end{equation}
	and the inner expansion,
	\begin{equation}
	\hat u_1(z_1;\epsilon)=\hat u^{(0)}_1(z_1)+\epsilon \hat u_1^{(1)}(z_1)+...,\label{inner1}
	\end{equation}
	Plugging into Eq. \eqref{ddm1},  we have the following equations, 
	\begin{align}
	&\text{At }{O}(\epsilon^{-3})\text{, }&e^{6z_1}(\hat u_1^{(0)}-g)=0,\label{bc11}\\
	&\text{At }{O}(\epsilon^{-2})\text{, }&\frac{d^2}{dz_1^2}\hat u_1^{(0)}-\frac{6}{1+e^{-6z_1}}\frac{d}{dz_1}\hat u_1^{(0)}-e^{6z_1}\hat u_1^{(1)}=0,\label{bc12}\\
	&\text{At }{O}(\epsilon^{-1})\text{, }&\frac{d^2}{dz_1^2}\hat u_1^{(1)}-\frac{6}{1+e^{-6z_1}}\frac{d}{dz_1}\hat u_1^{(1)}-e^{6z_1}\hat u_1^{(2)}=0,\label{bc13}\\
	&\text{At }={O}(1)\text{, }&\frac{d^2}{dz_1^2}\hat u_1^{(2)}-\frac{6}{1+e^{-6z_1}}\frac{d}{dz_1}\hat u_1^{(2)}-e^{6z_1}\hat u_1^{(3)}=f,\label{bc14}\\
	\end{align}
	It follows that
	\begin{align}
	\hat u_1^{(0)}=g,\label{ddm1_u0}\\
	\hat u_1^{(1)}=\hat u_1^{(2)} = 0,\label{ddm1_u1}\\
	\hat u_1^{(3)}=-fe^{-6z_1}.\label{ddm1_u3}
	\end{align}
	Thus, 
	\begin{equation}
	\hat u_1(z_1;\epsilon)=g-\epsilon^3fe^{-6z_1}+o(\epsilon^3).\label{ddm1_out_u1}
	\end{equation}
	Clearly, \(\hat u_1\) satisfies two of the matching conditions as \(z_1\rightarrow +\infty\) (Eqs. \eqref{mc3} and \eqref{mc4}). However, \(\displaystyle \lim_{z_1 \to -\infty} \hat u_1(z_1;\epsilon)=-\infty\), which implies \(\hat u_1\) does not satisfy the other two matching conditions (Eqs. \eqref{mc1} and \eqref{mc2}) and hence there exists another layer, namely \(K_2\), and another inner solution \(\hat u_2\) (see Fig. \ref{k12}).
	Accordingly, we introduce the layer $K_2$ where \(x\sim \alpha\epsilon\ln(\epsilon)+{O}(\epsilon)\). Here we choose \(\alpha=1/6\) and consider the two leading terms in the inner expansions, other choices of \(\alpha\) will result different layers, but their inner expansions are the same up to \({O}(\epsilon)\). We consider a new stretched variable, 
	\begin{equation}
	z_2=\frac{x-\epsilon(\ln\epsilon)/6}{\epsilon}=\frac{x}{\epsilon}-\frac{\ln\epsilon}{6}\label{k2}
	\end{equation}
	Similar to Eq. \eqref{match1}, we derive the corresponding matching conditions in \(K_2\), that is
	\begin{equation}
	\bar u_1(\epsilon z_2+\frac{\epsilon\ln\epsilon}{6}) \simeq \hat u_2(z_2), \text{ as } z_2\rightarrow -\infty~~{\rm and}~~\epsilon z_2\to 0^-. \label{mc_k2}
	\end{equation}
	It follows that
	\begin{align}
	&\lim_{z_2\rightarrow -\infty}\hat u_2^{(0)}(z_2)=\bar u_1^{(0)}(0),\label{mc_k2_1}\\
	&\lim_{z_2 \rightarrow -\infty}\hat u_2^{(1)}(z_2)=\bar u_1^{(1)}(0) + ( z_2+\frac{\ln\epsilon}{6})\frac{d}{dx}\bar u_1^{(0)}(0). \label{mc_k2_2}
	\end{align}
	As $z_2\to +\infty$, we need to match the inner solutions $\hat u_2(z_2)$ with $\hat u_1(z_1)$ from $K_1$. This is described below.
	
	Plugging in the inner expansion into DDM1 (Eq. \eqref{ddm1}) and using that  \(e^{6x/\epsilon}=\epsilon e^{6z_2}\), we obtain:
	\begin{align}
	&\text{At }{O}(\epsilon^{-2})\text{, }&\frac{d^2}{dz_2^2}\hat u_2^{(0)}-e^{6z_2}(\hat u_2^{(0)}-g)=0,\label{bc21}\\
	&\text{At }{O}(\epsilon^{-1})\text{, }&\frac{d^2}{dz_2^2}\hat u_2^{(1)}-6e^{6z_2}\frac{d}{dz_2}\hat u_2^{(0)}-e^{6z_2}\hat u_2^{(1)}=0.\label{bc22}
	\end{align}
	The general solution to Eq. \eqref{bc21} is
	\begin{align}
	&\hat u_2^{(0)}(z_2) =g+ C_1 I_0(\frac{e^{3z_2}}{3}) + C_2 K_0(\frac{e^{3z_2}}{3}), \label{ddm1_in_u20}
	\end{align}
	where \(I_0\) and \(K_0\) are the modified Bessel functions of the first and second kind, respectively, and \(C_1\), \(C_2\) are constants.
	
	Then we need to match \(\hat u_2^{(0)}(z_2)\) with the leading order of the inner solution in \(K_1\) (Eq. \eqref{ddm1_u0}) as \(z_2\rightarrow +\infty\), that is
	\begin{equation}
	\lim_{z_2 \rightarrow +\infty} \hat u_2^{(0)}(z_2)=  g. \label{mc_k2_ddm1_1}
	\end{equation}
	It follows that \(C_1=0\). From Eq. \eqref{mc_k2_1} and \eqref{mb3}, we obtain \(C_2=0\).
	Thus,
	\begin{equation}
	\hat u_2^{(0)}(z_2)=g,\label{k2_u0}
	\end{equation}
	which implies
	\begin{equation}
	\bar u_1^{(0)}(0)=\lim_{z_2 \to -\infty}\hat u_2^{(0)}(z_2)=g.\label{k2_u10}
	\end{equation}
	Hence, DDM1 recovers the Poisson equation with Dirichlet boundary condition at leading order and \(\bar u_1^{(0)}\) is the exact solution to Eqs. \eqref{p1} and \eqref{p2}.
	
	Plugging Eq. \eqref{k2_u0} into Eq. \eqref{bc22}, we have
	\begin{equation}
	\hat u_2^{(1)}(z_2) = C_3 I_0(\frac{e^{3z_2}}{3}) + C_4 K_0(\frac{e^{3z_2}}{3}), \label{mod_bessel1}
	\end{equation}
	where \(C_3\) and \(C_4\) are constants. \(\hat u_2^{(1)}\) needs to match \(\hat u_1^{(1)} = 0\) in \(K_1\) as \(z_2\rightarrow+\infty\), that is 
	\begin{align}
	&\lim_{z_2 \rightarrow +\infty}\hat u_2^{(1)}(z_2)=0, \label{mc_k2_ddm1_2}
	\end{align}
	Hence, \(C_3=0\).
	On the other side where \(z_2\rightarrow-\infty\), by the matching condition Eq. \eqref{mc_k2_2}, we obtain
	\begin{align}
	&\lim_{z_2 \rightarrow -\infty}\hat u_2^{(1)}(z_2)=\bar u_1^{(1)}(0) + A( z_2+\frac{\ln\epsilon}{6}),\label{mc_k2_ddm1_3}
	\end{align}
	where \(A=\frac{d}{dx}\bar u_1^{(0)}(0)\) is the derivative of exact solution at the boundary \(x=0\).
	Therefore, taking the limit of Eq. (\ref{mod_bessel1}) as $z_2\rightarrow -\infty$ and equating the result to Eq. (\ref{mc_k2_ddm1_3}),
	it follows that
	\begin{equation}
	\bar u_1^{(1)}(0)=-\frac{A}{6}\ln\epsilon+\frac{A}{3}(-\ln 6+ \gamma)\sim {O}(\ln\epsilon) \label{ddm1_out_u2}
	\end{equation}
	Hence, the convergence rate of DDM1 in the \(L^2\) norm is \(\epsilon\ln(\epsilon)\) if \(A\neq0\). Although the error at the boundary is \(\hat u_1(0)-g\sim {O}(\epsilon^3)\) as seen from Eq. \eqref{ddm1_out_u1}),  the \(L^\infty\) error is dominated by \(\epsilon\ln(\epsilon)\).
	
	As for DDM3, Replacing \(\frac{1-\phi}{\epsilon}\) with \(\frac{|\nabla \phi|}{\epsilon^2}\) and conducting an analogous analysis, we obtain 
	\begin{align}
	\bar u_1^{(0)}(0)=g, \label{ddm3_out_u1}\\
	\bar u_1^{(1)}(0)=-\frac{A}{6}\ln\epsilon+\frac{A}{3}(-\ln(\frac{1}{2}\sqrt{\frac{2}{3}})+\gamma)\sim {O}(\ln \epsilon). \label{ddm3_out_u2}
	\end{align}
	Thus, DDM3 is also of \({O}(\epsilon\ln\epsilon)\) in both \(L^2\) and \(L^\infty\) norms if \(A\neq0\). If \(A=0\), our analysis suggests that DDM1 and DDM3 are better than first-order accurate. In fact, we find the schemes are ${O}\left(\epsilon^2\left(\ln \epsilon\right)^2\right)$ accurate (see Appendix \ref{app_a}). 
	
	To validate our analysis, we plot \(C(\epsilon)\approx(u_\epsilon-u)/\epsilon\) versus \(\ln(\epsilon)\) as in Fig. \ref{u1_3} and compute the slope through a linear fit using the numerical results in Sec. \ref{num_1d_ddm}. For both DDM1 and DDM3, our asymptotic analysis suggests that the slope is \(-A/6\approx -0.185\) using $A=1.111$. Numerically, we obtain the slopes $-0.186$ and $-0.188$ from DDM1 and DDM3, respectively, which agrees well with the theory (see Tab. \ref{ubar1_mddm13} in Appendix \ref{app_c} for other cases).

	\subsection{Development and analysis of higher-order, modified DDMs} \label{ana_mddm}
	In order to achieve higher order accuracy, we need to guarantee that the first order term in the outer expansion (\(\bar u_1\)) vanishes, that is, \(\bar u_1^{(1)}=0\). From the matching condition, Eq. \eqref{mc2}, we observe that if \(\hat u_1^{(1)}(z_1)\) behaves as \(Az_1\) as \(z_1\to -\infty\), where \(A\) is the derivative of the exact solution at the boundary, then \(\bar u_1^{(1)}(0)\) must be 0 and hence \(\bar u_1^{(1)}(x)=0\) since $d^2\bar u_1^{(1)}/dx^2 =0$ from Eqs. (\ref{out2}) and (\ref{out2_1}). Therefore, if we can modify the original DDMs in a way such the inner solution \(\hat u_1(z_1)\) satisfies,
	\begin{align}
	\lim_{z_1 \to \pm\infty}\hat u_1^{(0)}(z_1) &= g,\\
	\lim_{z_1 \to +\infty}\hat u_1^{(1)}(z_1) &= 0,\\
	\lim_{z_1 \to -\infty}\hat u_1^{(1)}(z_1) &\sim Az_1,
	\end{align}
	we should be able to achieve higher-order accuracy based on the asymptotic analysis.  Accordingly, this suggests that to achieve higher-order accuracy we may modify the DDMs (referred to as the mDDMs) in the following way:
	\begin{eqnarray}
	&\textbf{mDDM1}: &\nabla\cdot(\phi \nabla u_\epsilon)-\frac{1}{\epsilon^3}(1-\phi)(u_{\epsilon}-g-r\bm{n}\cdot \nabla u_\epsilon)=\phi f  , \label{eq_mddm1}\\
	&\textbf{mDDM2}:&\phi \Delta u_\epsilon-\frac{1}{\epsilon^2}(1-\phi)(u_{\epsilon}-g-r\bm{n}\cdot \nabla u_\epsilon)=\phi f, \label{eq_mddm2}\\
	&\textbf{mDDM3}:&\nabla\cdot(\phi \nabla u_\epsilon)-\frac{|\nabla\phi|}{\epsilon^2}(u_{\epsilon}-g-r\bm{n}\cdot \nabla u_\epsilon)=\phi f, \label{eq_mddm3}
	\end{eqnarray}
	where \(\bm{n} = -\frac{\nabla \phi}{|\nabla \phi|}\) is the normal vector and the new term $r\mathbf{n}\cdot\nabla u_\epsilon$ is designed to cancel the contribution of the derivative of the exact solution (e.g., the $A$ terms) in the outer solutions $\bar u^{(1)}_1$. {\color{black} This is similar to, but simpler than, the active penalty method in \cite{shirokoff2015sharp}, which introduces penalty terms to match the solution and normal derivatives across the physical boundary. }{\color{blue} The modified boundary condition $g+r\mathbf{n}\cdot\nabla u_\epsilon$ is simply a linear approximation of the solution near the boundary in the normal direction to the boundary.}
	Next, we analyze these methods in the following two sections to determine their actual order of accuracy.
	\subsubsection{Asymptotic analysis of mDDM1 and mDDM3}\label{aa_mddm1_3}
\begin{figure}
		\center
		\includegraphics[width=.6\textwidth]{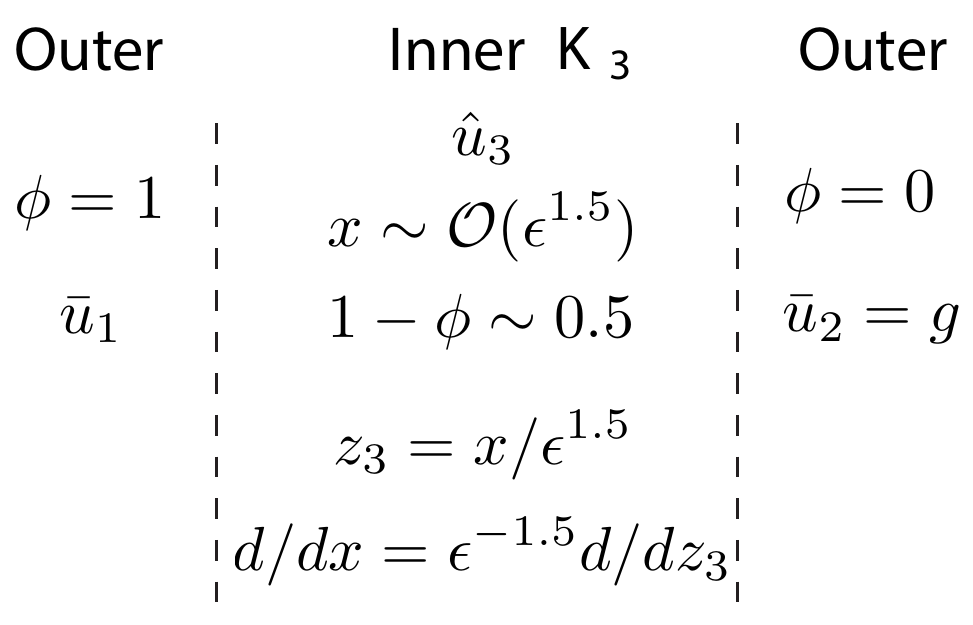}
		\caption{Inner layer \(K_3\)}\label{k3}
	\end{figure}

	Similar to Sec. \ref{sec_ana}, we assume \(x=0\) is the boundary and \(r=x\). Consider mDDM1 in 1D:
	\begin{equation}
	\frac{d}{dx}(\phi \frac{d}{dx} u_\epsilon)-\frac{1}{\epsilon^3}(1-\phi)(u_\epsilon-g-x\frac{d}{dx}u_\epsilon)=\phi f \label{mddm10},
	\end{equation}
	Substituting \(\phi\) with Eq. \eqref{phi}, we obtain
	\begin{equation}
	\frac{d^2}{dx^2}u_{\epsilon} - \frac{6}{\epsilon}\frac{1}{1+e^{-6x/\epsilon}}\frac{d}{dx}u_{\epsilon}-\frac{e^{6x/\epsilon}}{\epsilon^3}(u_{\epsilon}-g-x\frac{d}{dx}u_\epsilon) = f. \label{mddm11}
	\end{equation}
	After examining possible ways of performing matched asymptotic expansions for Eq. \eqref{mddm11}, the simplest approach is to choose a new inner variable \(z_3=\frac{x}{\epsilon^{1.5}}\) and the corresponding inner solution \(\hat u_3(z_3)\) in the region \(K_3\) (see Fig. \ref{k3}), where \(x\sim \epsilon^{1.5}\) and \(e^{6x/\epsilon}\approx 1\). In the local coordinate system, the derivatives become
	\begin{eqnarray}
	\frac{d}{dx} = \frac{1}{\epsilon^{1.5}}\frac{d}{dz_3},\label{mddm1_loc_d1}\\
	\frac{d^2}{dx^2} = \frac{1}{\epsilon^3}\frac{d^2}{dz_3^2},\label{mddm1_loc_d2}
	\end{eqnarray} 
	Similar to Eqs. \eqref{match1} and \eqref{match2}, we develop matching conditions for \(K_3\), that is,
	\begin{align}
	&\hat{u}_3(z_3;\epsilon) \simeq  \bar u_1(\epsilon^{1.5} z_3;\epsilon) , \text{ as } z_3\rightarrow-\infty, ~~~{\rm and}~~~\epsilon^{1.5}z_3\rightarrow 0^- \label{match_mddm1_1}\\
	& \hat{u}_3(z_3;\epsilon)\simeq\bar u_2(\epsilon^{1.5} z_3;\epsilon)=g , \text{ as } z_3\rightarrow+\infty, ~~~{\rm and}~~~\epsilon^{1.5}z_3\rightarrow 0^+. \label{match_mddm1_2}
	\end{align}
	Since we use the scale \({\epsilon^{1.5}}\), it is natural to include half powers in the outer and inner expansions:
	\begin{align}
	\bar u_1(x;\epsilon)&=\bar u_1^{(0)}(x)+\epsilon^{0.5} \bar u_1^{(0.5)}(x)+\epsilon \bar u_1^{(1)}(x)+\epsilon^{1.5} \bar u_1^{(1.5)}(x)...,\label{outer_mddm1}\\ 
	\hat u_3(z_3;\epsilon)&=\hat u_3^{(0)}(z_3)+\epsilon^{0.5} \hat u_3^{(0.5)}(z_3)+\epsilon \hat u_3^{(1)}(z_3)+\epsilon^{1.5} \hat u_3^{(1.5)}(z_3)....\label{inner_mddm1}
	\end{align}
	Plugging into Eqs. \eqref{match_mddm1_1} and \eqref{match_mddm1_2}, we obtain the matching conditions,
	\begin{align}
	&\hat u_3^{(0)}(z_3) = g, \text{ as }z_3\rightarrow +\infty,\label{mddm1_mc1}\\
	&\hat u_3^{(k)}(z_3) = 0, k=0.5, 1, {\color{black}1.5},...,\text{ as }z_3\rightarrow +\infty,\label{mddm1_mc2}
	\end{align}
	and
	\begin{align}
	&\hat u_3^{(k)}(z_3) = \bar u_1^{(k)}(0), k=0, 0.5, 1, \text{ as } z_3\rightarrow -\infty,\label{mddm1_mc3}\\
	&\hat u_3^{(1.5)}(z_3) = \bar u_1^{(1.5)}(0)+z_3\frac{d}{dx}\bar u_1^{(0)}(0),\text{ as } z_3\rightarrow -\infty.\label{mddm1_mc4}
	\end{align}
	Plugging Eq. \eqref{inner_mddm1} into Eq. \eqref{mddm11} we derive the following inner equation at the leading order \({O}(\epsilon^{-3})\):
	\begin{align}
	\frac{d^2}{dz_3^2}\hat u_3^{(0)}-(\hat u_3^{(0)}-g-z_3\frac{d}{dz_3}\hat u_3^{(0)})=0,\label{1_i_1}
	\end{align}
	Clearly, \(\hat u_3^{(0)}=g\) is a solution to Eq. \eqref{1_i_1}. Note that the general solution to the following homogeneous ordinary differential equation,
	\begin{equation}
	y''+xy'-y=0, \label{hode1}
	\end{equation}
	involves a linear combination of parabolic cylinder functions ($D_{-2}(x)$) \cite{bender2013advanced}:
	\begin{equation}
	y = e^{-x^2/4}(C_1D_{-2}(x)+C_2D_{-2}(-x)),
	\end{equation}
	where \(C_1\) and \(C_2\) are constants. Thus,
	\begin{equation}
	\hat u_3^{(0)} = g + e^{-z_3^2/4}(C_1D_{-2}(z_3)+C_2D_{-2}(-z_3)),
	\end{equation}
	where \(D_{-2}(z_3)\sim z_3^2e^{-z_3^2/4}\) as \(z_3\rightarrow +\infty\) and \(D_{-2}(z_3)=-\sqrt{2\pi}z_3e^{z_3^2/4}\) as \(z_3\rightarrow -\infty\). By the matching conditions (Eqs. \eqref{mddm1_mc1} and \eqref{mddm1_mc3}), we have $C_1=C_2=0$ and
	\begin{align}
	\hat u_3^{(0)}=g, \label{1_uhat_11}\\
	\bar u_1^{(0)}(0)=g. \label{1_ubar_11}
	\end{align}
	Hence mDDM1 recovers the Poisson equation at the leading order, e.g., \(\bar u_1^{(0)}=u\), the exact solution. 
	
	At the next order \({O}(\epsilon^{-2.5})\) in Eq. \eqref{mddm11}, we have,
	\begin{align}
	\frac{d^2}{dz_3^2}\hat u_3^{(0.5)}-(\hat u_3^{(0.5)}-z_3\frac{d}{dz_3}\hat u_3^{(0.5)})=0.\label{1_i_2}
	\end{align}
	Thus,
	\begin{equation}
	\hat u_3^{(0.5)} = e^{-z_3^2/4}(C_1D_{-2}(z_3)+C_2D_{-2}(-z_3)).
	\end{equation}
	Again, using the matching conditions (Eqs. \eqref{mddm1_mc2} and \eqref{mddm1_mc3}), we derive \(\hat u_3^{(0.5)}=0\). 
	A similar argument at \({O}(\epsilon^{-2})\) in Eq. \eqref{mddm11} gives  \(\hat u_3^{(1)}=0\).
	
	At \({O}(\epsilon^{-1.5})\) in Eq. \eqref{mddm11}, we obtain,
	\begin{align}
	\frac{d^2}{dz_3^2}\hat u_3^{(1.5)}-(\hat u_3^{(1.5)}-z_3\frac{d}{dz_3}\hat u_3^{(1.5)})=0.\label{1_i_4}
	\end{align}
	Hence,
	\begin{equation}
	\hat u_3^{(1.5)}=e^{-z_3^2/4}(C_1D_{-2}(z_3)+C_2D_{-2}(-z_3)).
	\end{equation}
	Applying the matching conditions (Eqs. \eqref{mddm1_mc2} and \eqref{mddm1_mc4}), we get 
	\begin{align}
	\hat u_3^{(1.5)} = -\frac{Ae^{-z_3^2/4}}{\sqrt{2\pi}}D_{-2}(z_3)\sim Az_3 \text{ as }z_3\rightarrow-\infty\label{1_uhat_3},
	\end{align}
	where \(A=\frac{d}{dx}\bar u_1^{(0)}(0)\) is the derivative of the exact solution at the boundary. Plugging into Eq. \eqref{mddm1_mc4}, we conclude
	\begin{equation}
	\bar u_1^{(1.5)}(0)=0. \label{1_ubar_1.5}
	\end{equation}
	Thus, the asymptotic analysis suggests that \(\bar u_1(x) = \bar u_1^{(0)}(x)+{O}(\epsilon^2)\) and mDDM1 is 2nd order accurate in \(\epsilon\) in the \(L^2\) norm as desired.

	To analyze the error in the \(L^\infty\) norm, we need to consider the error at the boundary. At the boundary, \(\hat u_3(0)-g \sim \epsilon^{1.5}\hat u_3^{(1.5)}(0) = -\frac{A\epsilon^{1.5}}{\sqrt{2\pi}}\). This suggests that mDDM1 is 1.5 order accurate in \(\epsilon\) in the \(L^\infty\) norm.
	
	This analysis can be easily extended to higher dimensions and to mDDM3 to obtain the same conclusions. The only difference with respect to mDDM3 is that \(\frac{1-\phi}{\epsilon^3}\) is replaced with \(\frac{|\nabla \phi|}{\epsilon^2}\), which gives \(\hat u_3(0) - g \sim \epsilon^{1.5}\hat u_3^{(1.5)}(0)=
	-\frac{A\epsilon^{1.5}}{\sqrt{6\pi}}\). The errors in the \(L^2\) norm and \(L^\infty\) norm are still suggested to be $2$nd and $1.5$ orders, respectively. 
	
	To confirm our analysis, we compute \(\hat u_3^{(1.5)}(0)\) numerically in the following way. We plot the values of \((u_\epsilon-u)/\epsilon^{1.5}\) at the boundary (\(z_3=0\)) versus \(\sqrt{\epsilon}\), where \(u_\epsilon\) is the numerical solution, and we find the y-intercept of a quadratic fit. Since \(A=1.111\), our asymptotic analysis gives that the values of \(\hat u_3^{(1.5)}(0)\) are -0.443 and -0.256 for mDDM1 and mDDM3, respectively. We obtain -0.442 and -0.257 from the numerical results of mDDM1 and mDDM3, respectively, which is consistent with our theory (see Tab. \ref{uhat1_mddm13} in Appendix \ref{app_c} for other cases).
	
	\subsubsection{Asymptotic analysis of mDDM2} \label{mddm2_ana}
	Even though mDDM2 is obtained using the same modification as in mDDM1 and mDDM3, surprisingly as our analysis below suggests and numerical results confirm (see Sec. \ref{mddm_nu}), mDDM2 is only first-order accurate in both $L^2$ and $L^\infty$. In 1D, the mDDM2 is:
	\begin{equation}
	\phi \frac{d^2}{dx^2}u_\epsilon-\frac{1}{\epsilon^2}(1-\phi)(u_\epsilon-g-x\frac{d}{dx}u_\epsilon)=\phi f. \label{mddm20}
	\end{equation}
	Substituting \(\phi\) from Eq. \eqref{phi}, we obtain
	\begin{equation}
	\frac{d^2}{dx^2}u_{\epsilon} -\frac{e^{6x/\epsilon}}{\epsilon^2}(u_{\epsilon}-g-x\frac{d}{dx}u_\epsilon) = f. \label{mddm21}
	\end{equation} 
	Considering Eq. \eqref{mddm21} in \(K_1\), we use the inner variable \(z_1=\frac{x}{\epsilon}\) to derive the following inner equations for \(\hat u_1\):
	\begin{align}
	&\text{At }{O}(\epsilon^{-2})\text{, }&\frac{d^2}{dz_1^2}\hat u_1^{(0)}-e^{6z_1}(\hat u_1^{(0)}-g-z_1\frac{d}{dz_1}\hat u_1^{(0)})=0,\label{mddm2_in_eq1}\\
	&\text{At }{O}(\epsilon^{-1})\text{, }&\frac{d^2}{dz_1^2}\hat u_1^{(1)}-e^{6z_1}(\hat u_1^{(1)}-z_1\frac{d}{dz_1}\hat u_1^{(1)})=0.\label{mddm2_in_eq2}
	\end{align} 
	Clearly, \(\hat u_1^{(0)}=g\) is a solution to Eq. \eqref{mddm2_in_eq1}, and the general solution to the homogeneous ordinary differential equation below,
	\begin{equation}
	y''-e^{6x}(y-xy')=0,\label{mddm2_ode}
	\end{equation}
	is given by,
	\begin{align}
	y=C_1x+C_2(-e^{e^{6x}(1-6x)/36}-x\int_0^xh(t)dt),
	\end{align}
	where \(h(x)=e^{e^{6x}(1-6x)/36+6x}\) (see Appendix \ref{app_d} for details).
	It follows that
	\begin{align}
	\hat u^{0}_1(z_1) &= g + C_1z_1+C_2(-e^{e^{6z_1}(1-6z_1)/36}-z_1\int_0^{z_1}h(t)dt),\\
	\hat u^{1}_1(z_1) &=  C_3z_1+C_4(-e^{e^{6z_1}(1-6z_1)/36}-z_1\int_0^{z_1}h(t)dt),
	\end{align}
	where the \(C_i\) are constants and 
	\begin{align}
	\lim_{z_1\rightarrow -\infty}e^{e^{6z_1}(1-6z_1)/36} = 1,\\
	\lim_{z_1\rightarrow +\infty}e^{e^{6z_1}(1-6z_1)/36} = 0,\\
	\lim_{z_1\rightarrow -\infty}
	\int_0^{z_1}h(t)dt \approx 0.17,\\
	\lim_{z_1\rightarrow +\infty}
	\int_0^{z_1}h(t)dt \approx 2.75.
	\end{align}
	Combining with the matching conditions as \(z_1\rightarrow\pm\infty\) (Eqs. \eqref{mc1}-\eqref{mc4}) and solving for the constants \(C_i\), we obtain that
	\begin{align}
	\hat u^{(0)}_1(z_1) &= g,\\
	\hat u^{(1)}_1(z_1) &\approx \frac{A}{2.92}(-e^{e^{6z_1}(1-6z_1)/36}+z_1\int_{z_1}^\infty h(t)dt),\label{mddm2_inner1}
	\end{align}
	where as before $A=\frac{d}{dx}\bar u_1^{(0)}(0)$. Therefore,
	\begin{align}
	\lim_{z_1 \to -\infty}\hat u^{(1)}(z_1)\sim Az_1-\frac{A}{2.92},
	\end{align}
	which implies 
	\begin{align}
	\bar u_1^{(1)}(0) \approx -\frac{A}{2.92}\neq 0, \text{ if }A\neq 0.\label{mddm2_bar_u1}
	\end{align}
	This suggests mDDM2 is 1st order accurate in \(\epsilon\) in the \(L^2\) norm when $A\ne 0$.
	Plugging \(z_1=0\) into Eq. \eqref{mddm2_inner1} and we have
	\begin{align}
	\hat u_1^{(1)}(0) \approx -\frac{Ae^{1/36}}{2.92},\label{mddm2_hat_u1}
	\end{align}
	which implies that mDDM2 is also 1st order accurate in the \(L^\infty\) norm.
	
	To confirm our analysis,  we compute \(\bar u_1^{(1)}\) numerically, which is given by  \((u_\epsilon-u)/\epsilon\) as \(z<<0\), using the numerical solution \(u_\epsilon\) in Sec. \ref{mddm_nu}. In addition, we interpolate and plot the boundary values  of \((u_\epsilon-u)/\epsilon\) at $z_1=0$ versus \(\epsilon\), and then calculate \(\hat u^{(1)}(0)\) from the y-intercept of a linear fit. Since \(A=1.111\), our analysis indicates that \(\bar u_1^{(1)}\approx-\frac{A}{2.92}\approx-0.381\) and \(\hat u^{(1)}(0)\approx -\frac{Ae^{1/36}}{2.92}\approx-0.391\). The numerical results yield \(\bar u_1^{(1)}\approx-0.380\) and \(\hat u^{(1)}(0)\approx-0.391\), which agrees very well with the theory (see Tab. \ref{ubar1__mddm2} in Appendix \ref{app_c} for other cases).
	
	\subsection{Numerical results for mDDM1-3} \label{mddm_nu}
	We adopt the same problem setup and an analogous discretization as in Sec. \ref{num_1d_ddm}. The results using $h=\epsilon/4$ are presented in Fig. \ref{mddm_all}. Clearly mDDM3 and mDDM1 are more accurate than DDM1-3 and mDDM2. Further, mDDM2 and DDM2 (Fig. \ref{ddm_all}) display similar levels of accuracy as suggested by our analysis. 
	
	To test the order of accuracy in $\epsilon$ predicted by theory, we next take $h=\epsilon^{1.5}/4$ in order to resolve the inner layer $K_3$, as guided by our asymptotic analysis in Sec. \ref{aa_mddm1_3}. Later, in Secs. \ref{heat1d} and \ref{2d}, we will test convergence with $\epsilon\propto h$. In Tabs. \ref{mddm1_case12} and \ref{mddm2_case12}, we present convergence results for mDDM1-3 for case 1 (see also Tabs. \ref{af1}-\ref{af6} in Appendix \ref{app_b} for other cases). Consistent with our analysis, we observe that mDDM1 and mDDM3 are approximately 2nd order accurate in the \(L^2\) norm but approximately 1.5 order accurate in the \(L^\infty\) norm. Also as predicted, mDDM2 is only 1st order accurate in both norms. Therefore, in practice we recommend using mDDM1 or mDDM3. 
	
	However, taking \(h\sim\epsilon^{1.5}\) is often too constraining for solving problems in 2D and 3D. Instead, one may choose \(h=\epsilon/c\). Doing this, one can obtain results with orders of accuracy ranging from 1.5 to 2 (1.5 for large c and approximately 2 for \(c\approx 1\); see Tabs. \ref{ddmts1} and \ref{ddmts2} in Sec. \ref{heat1d} and Tabs. \ref{mDDMt3_mov_l2}, \ref{mDDMt3_mov_inf} and \ref{err_2dt_mov} in Appendices).   
	
	\begin{table}[H]
		\footnotesize 
		\center
		\begin{tabular}{c||c|c|c|c||c|c|c|c}
			&\multicolumn{4}{c||}{mDDM1}&\multicolumn{4}{c}{mDDM3}\\ \hline
			\(\epsilon\)&\(E^{(2)}\)&\(\bm{k}\)&\(E^{(\infty)}\)&\(\bm{k}\)&\(E^{(2)}\)&\(\bm{k}\)&\(E^{(\infty)}\)&\(\bm{k}\)
			\begin{filecontents*}{mddm13_case1.csv}
1.10E+01,,,,,,,,
2.00E-01,4.88E-01,0.00 ,2.10E-01,0.00 ,1.55E-01,0.00 ,8.23E-02,0.00 
1.00E-01,1.50E-01,1.70 ,6.59E-02,1.67 ,3.77E-02,2.04 ,2.23E-02,1.88 
5.00E-02,4.46E-02,1.75 ,1.98E-02,1.74 ,8.54E-03,2.14 ,5.96E-03,1.91 
2.50E-02,1.26E-02,1.83 ,5.58E-03,1.83 ,1.89E-03,2.17 ,2.03E-03,1.55 
1.25E-02,3.34E-03,1.91 ,1.55E-03,1.85 ,4.16E-04,2.19 ,6.15E-04,1.72 
6.25E-03,8.61E-04,1.95 ,4.91E-04,1.65 ,8.74E-05,2.25 ,2.29E-04,1.42 
\end{filecontents*}
\csvreader[]{mddm13_case1.csv}{}
			{\\\hline \csvcoli & \csvcolii & \csvcoliii & \csvcoliv & \csvcolv & \csvcolvi & \csvcolvii & \csvcolviii & \csvcolix}
		\end{tabular}
		\caption{The \(L^2\) and \(L^{\infty}\) errors for mDDM1 and mDDM3 with \(h=\epsilon^{1.5}/4\).}\label{mddm1_case12}
	\end{table}
	\begin{table}[H]
		\footnotesize 
		\center
		\begin{tabular}{c||c|c|c|c}
			&\multicolumn{4}{c}{mDDM2}\\ \hline
			\(\epsilon\)&\(E^{(2)}\)&\(\bm{k}\)&\(E^{(\infty)}\)&\(\bm{k}\)
			\begin{filecontents*}{mddm2_case1.csv}
1.10E+01,,,,
2.00E-01,4.58E-01,0.00 ,2.19E-01,0.00 
1.00E-01,2.24E-01,1.03 ,1.06E-01,1.05 
5.00E-02,1.11E-01,1.02 ,5.19E-02,1.03 
2.50E-02,5.47E-02,1.02 ,2.59E-02,1.00 
1.25E-02,2.72E-02,1.01 ,1.29E-02,1.01 
6.25E-03,1.36E-02,1.00 ,6.45E-03,1.00 
\end{filecontents*}
\csvreader[]{mddm2_case1.csv}{}
			{\\\hline \csvcoli & \csvcolii & \csvcoliii & \csvcoliv & \csvcolv }
		\end{tabular}
		\caption{The \(L^2\) and \(L^{\infty}\) errors for mDDM2 with \(h=\epsilon^{1.5}/4\).}\label{mddm2_case12}
	\end{table}
	
	
	\section{Time-dependent problems}\label{heat1d}
	We next extend the mDDMs to simulate time-dependent PDEs in a moving domain \(D(t)\) with Dirichlet boundary conditions. As we discuss later, our approach and analysis holds for much more general time-dependent equations, but we consider the diffusion equation here for simplicity of presentation:
	\begin{align}
	\partial_t u&=\Delta u + f  \text{ in }  D(t)\label{heat1},\\
		u&=g \text{ on }  \partial D(t), \label{heat2}\\
	u(\mathbf{x},0)&=u_0(\mathbf{x}) \text{ in }  D(t). \label{heat3}
		\end{align}
	
	\subsection{Derivation and analysis of high-order, modified DDMs for time-dependent PDEs}\label{mddmt3}
	To approximate Eqs. \eqref{heat1} and \eqref{heat2}, we formulate the diffuse domain model using an approximation (mDDMt3), which is analogous to mDDM3 in Sec. \ref{aa_mddm1_3}:
	\begin{equation}
	\textbf{mDDMt3:} ~~~\partial_t(\phi u_\epsilon)=\nabla(\phi \nabla u_\epsilon)-\frac{1}{\epsilon^2}|\nabla \phi|(u_\epsilon-g-r\bm{n}\cdot \nabla u_\epsilon)+\phi f \label{tddm1}
	\end{equation}
	where as before \(r\) is the signed distance function, with $r<0$ denoting the interior of $D$. Since the domain may be time-dependent, \(\phi\) may depend on time and thus $\phi$ needs to be in the time derivative \cite{li09}. Note that in time-dependent problems, the initial condition needs to be extended outside $D$ as well although this extension need not be constant in the normal direction.
	
	The asymptotic analysis for the time-dependent problem is very similar to that for the Poisson equation presented in Sec. \ref{aa_mddm1_3}. The only difference is that now the time derivative must be analyzed as well. Assuming the inner and outer expansions and matching conditions hold as in Sec. \ref{aa_mddm1_3}, from the leading order outer equation, we find that \(\bar u_1^{(0)}\) satisfies Eq. (\ref{heat1}).
	
	In the inner expansion, taking \(z_3=x/\epsilon^{1.5}\), the time derivative can be written as: 
	\begin{equation}
	\partial_t=-\frac{v}{\epsilon^{1.5}}\partial_{z_3}+o(1) \label{time_d}
	\end{equation}
	where \(v\) is the normal velocity of the domain boundary. 
	Since \(\partial_t\) is of order \(\epsilon^{-1.5}\), the time derivative will not affect the higher-order terms in the inner expansions of mDDM3.
	Therefore, we can still derive \(\hat u_3^{(0)} = g\) and \(\hat u_3^{(0.5)} = \hat u_3^{(1)} = 0\). Plugging into the inner equation at the next leading order (\({O}(\epsilon^{-1.5})\)), we obtain the same solution for \(\hat u_3^{(1.5)}\) as that for the time-independent equation. Thus, mDDMt3 should be ${O}(\epsilon^2)$ accurate in $L^2$ and ${O}(\epsilon^{1.5})$ accurate in $L^\infty$
 for the time-dependent (diffusion) equation. Analogously, mDDM1 has the same order of accuracy as mDDM3 and mDDM2 is ${O}(\epsilon)$ accurate in both $L^2$ and $L^\infty$ (results not shown).
 	
	\subsection{Numerical results using mDDMt3 for time-dependent problems}\label{1D time dependent}
	For time-dependent problems, 
	we first assume the domain does not change in time. Then, taking \(f=\cos(x)\cos(t)+\cos(x)\sin(t)\) and \(g=\cos(1.111)\sin(t)\), we obtain the exact solution \(u=\cos(x)\sin(t)\) on the domain \(D\)=[-1.111,1.111]. Therefore we take the initial $u_\epsilon(x,0)=0$ and the constant normal extension of $f$ is obtained analytically.
	We use central difference discretizations in space (as in Sec. \ref{num_1d_ddm}) and the Crank-Nicholson method discretization in time to solve Eq. \eqref{tddm1} on the larger domain \(\Omega=[-2,2]\) with \(\phi\) from Eq. \eqref{phi}. We still use the Thomas algorithm to solve the tridiagonal matrix system. We calculate the errors at \(t=1\) in the \(L^2\) and \(L^\infty\) norms by setting $dt=h$ where  \(dt\) and \(h\) are the time step and grid size, respectively. Here, instead of taking $h\propto \epsilon^{1.5}$ as in the previous section, we instead take $h={\epsilon}/{c}$ where $c$ is a constant.
	
\begin{figure}[H]
		\includegraphics[width=.75\linewidth]{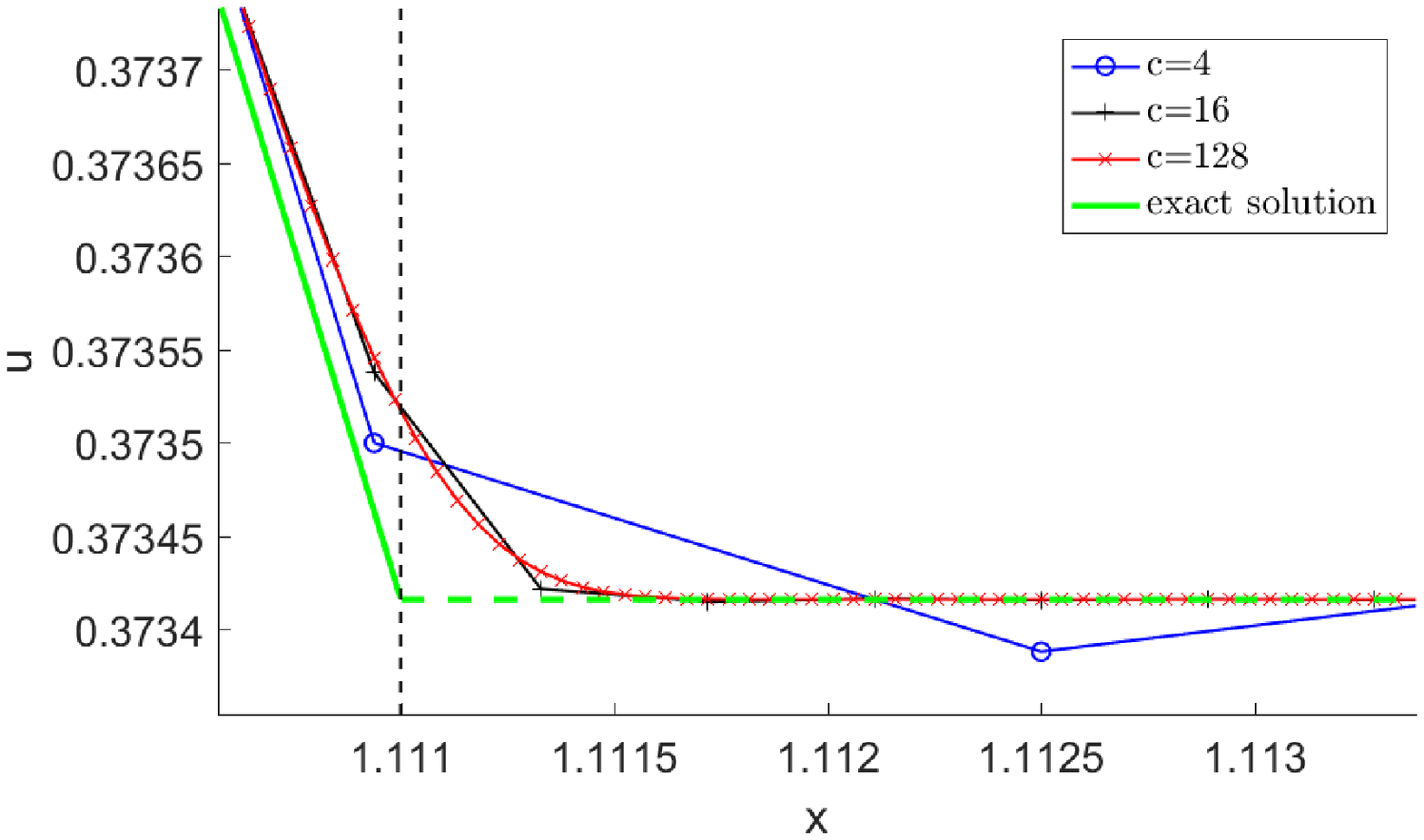}
		\caption{Numerical solutions of mDDMt3 near the boundary \(x=1.111\) on the 1D stationary domain with different \(c\). The green dashed line denotes the extension (e.g., constant in the normal direction) of the exact solution out of the domain. The black dashed line marks the right boundary $x=1.111$.}\label{comp_c_mddmt3_00625_sta}
	\end{figure}

	Various choices of \(c\) have been tested and here we present results using \(c=4,16,128\) in Tabs. \ref{ddmts1}, \ref{ddmts2}. 
	We observe that mDDMt3 is 2nd order accurate in the \(L^2\) norm for all choices of \(c\). However, in $L^\infty$, we find that when $c=4$, the schemes are roughly 2nd order accurate in $\epsilon$ while if $c=128$ the order of accuracy decreases to $1.5$, consistent with our theory. In Fig.  \ref{comp_c_mddmt3_00625_sta}, we show the numerical solutions of mDDMt3 near the boundary $x=1.111$ for difference choices of \(c\). When $c=4$, the boundary layer $K_3$ is not resolved and hence the numerical solution is less smooth near the boundary than is the continuous solution of mDDMt3. {\color{black} Nevertheless, the solution with $c=4$ is actually closer to the exact solution close to the boundary of $D$. This occurs because the discretization error, $u_{h,\epsilon}-u_{\epsilon}$, and analytic error, $u_\epsilon-u$, have opposite signs but similar magnitudes so that they can at least partially cancel one another in the total error $u_{h,\epsilon}-u=(u_{h,\epsilon}-u_{\epsilon})+(u_\epsilon-u)$, yielding higher order accuracy than might be expected. As $h$ decreases, the total error eventually becomes dominated by the analytic error, which results in errors that scale like ${O}(\epsilon^{1.5})$ in \(L^\infty\) as predicted by theory. See App. \ref{res_trunc_analytic} for details.}

	Outside of $D$, the mDDMt3 solutions tend to the extension of the exact solution, that is constant in the normal direction, because the boundary condition at $\partial\Omega$ is equal to the boundary condition at $\partial D$. Other choices of boundary conditions at $\partial\Omega$ would yield similar errors near $\partial D$ but different behavior far from $D$ (results not shown).
	
	\begin{table}[H]
		\center
		\begin{tabular}{c|c|c|c|c|c|c}
			\bfseries $\epsilon$ & c=4 &\bfseries k  & c=16 &\bfseries k& c=128 &\bfseries k
			\begin{filecontents*}{mddmt3_l2.csv}
eps,rate=4,k,rate =16,k,rate=128,k
0.2,2.85E-02,0.00 ,2.88E-02,0.00 ,2.88E-02,0.00 
0.1,7.19E-03,1.99 ,7.39E-03,1.96 ,7.40E-03,1.96 
0.05,1.68E-03,2.10 ,1.76E-03,2.07 ,1.76E-03,2.07 
0.025,4.08E-04,2.04 ,4.19E-04,2.07 ,4.20E-04,2.07 
0.0125,9.71E-05,2.07 ,1.03E-04,2.02 ,1.03E-04,2.03 
0.00625,2.48E-05,1.97 ,2.61E-05,1.98 ,2.62E-05,1.98 
\end{filecontents*}
\csvreader{mddmt3_l2.csv}{}
			{\\\hline \csvcoli & \csvcolii & \csvcoliii & \csvcoliv & \csvcolv & \csvcolvi & \csvcolvii}
		\end{tabular}
		\caption{The \(L^{2}\) errors for simulating the time-dependent (diffusion) equation using mDDMt3 with \(h=\epsilon/c\) on a fixed domain.}\label{ddmts1}
	\end{table}
	\begin{table}[H]
		\center
		\begin{tabular}{c|c|c|c|c|c|c}
			\bfseries $\epsilon$& c=4 &\bfseries k  & c=16 &\bfseries k& c=128 &\bfseries k
			\begin{filecontents*}{mddmt3_inf.csv}
eps,rate=4,k,rate =16,k,rate=128,k
0.2,3.87E-02,0.00 ,3.85E-02,0.00 ,4.01E-02,0.00 
0.1,1.00E-02,1.95 ,1.10E-02,1.81 ,1.18E-02,1.76 
0.05,2.15E-03,2.22 ,3.17E-03,1.79 ,3.56E-03,1.73 
0.025,5.00E-04,2.10 ,1.08E-03,1.55 ,1.10E-03,1.69 
0.0125,1.59E-04,1.65 ,3.27E-04,1.72 ,3.41E-04,1.69 
0.00625,4.38E-05,1.86 ,8.86E-05,1.88 ,1.15E-04,1.57 
\end{filecontents*}
\csvreader{mddmt3_inf.csv}{}
			{\\\hline \csvcoli & \csvcolii & \csvcoliii & \csvcoliv & \csvcolv& \csvcolvi & \csvcolvii}
		\end{tabular}
		\caption{The \(L^{\infty}\) errors for simulating the time-dependent (diffusion) equation using mDDMt3 with \(h=\epsilon/c\) on a fixed domain.}\label{ddmts2}
	\end{table}
	
	Next we assume the boundary of the domain is also moving and set \(x_l(t) = -1.111\), \(x_r(t) = 1.111+0.5t\), where \(x_l(t)\) and \(x_r(t)\) represent the left and right hand sides of the domain \(D\), respectively.
	We take \(f=\cos(x)\cos(t)+\cos(x)\sin(t)\), 
	\(g=\cos(x_*)\sin(t)\), where \(x_* = x_l,x_r\), and the exact solution is \(u=\cos(x)\sin(t)\).
	
	In general, neither the signed distance function \(r\) nor the function \(\phi\) are given analytically. 
	Because we need to use the signed distance function in the modified DDMs, we find it convenient to apply the level-set method \cite{Osher-1988,Gibou-2018} to determine both \(r\) and \(\phi\). Therefore, we solve the following Hamilton-Jacobi equation:
	
	\begin{equation}
	\partial_t r+v|\nabla r|=0\label{hj},
	\end{equation}
	using a 5th order upwind WENO scheme \cite{shu11} and a 2nd order Total Variation Diminishing (TVD) Runge-Kutta (RK) method time discretization to obtain an accurate fully discrete solution \cite{gottlieb1998total}.
	In order to keep \(|\nabla r|=1\), we periodically perform reinitialization \cite{min10,rus00} by solving the following equation,
	\begin{equation}
	\partial_\tau r=\text{sgn}_h(r^0)(1-|\nabla r|), \label{rein}
	\end{equation}
	where \(\tau\) is pseudo time and \(\text{sgn}_h(r^0\)) is a smoothed approximation function of the sign of the initial signed distance function \cite{peng1999pde}:
	\begin{equation}	
	\text{sgn}_h(r)=\left\{
	\begin{array}{lr}
	-1 & r<-h\\
	1 & r>h\\
	\frac{r}{h}+\frac{1}{\pi}\sin(\frac{\pi r}{h} ) & \text{else}
	\end{array}
	\right.\label{sign_r}
	\end{equation}
	To determine whether reinitialization is needed or not, we calculate the slope of \(r\) near the boundaries of \(D(t)\). Then we check if the maximum difference between the absolute value of the slope and 1 exceeds a threshold, here taken to be 0.01. If so, then we perform reinitialization.
	
	Once the signed distance function \(r(x,t)\) is obtained, we construct the phase field function \(\phi\) through
	Eq. \eqref{phi}.
	Next, we solve mDDMt3 using the same numerical setup as before. The results are very similar to the stationary domain case. We again find that mDDMt3 is 2nd order accurate in the \(L^2\) norm and between 1.5 and 2nd order accurate in the \(L^\infty\) norm, which is consistent with our analysis.
	{\color{black}See Tabs. \ref{mDDMt3_mov_l2} and \ref{mDDMt3_mov_inf} in App. \ref{res_1d_mov}.}
	
	\section{2D numerical results}\label{2d}
\begin{figure}[H]
		\center
		\begin{subfigure}{.49\textwidth}
			\includegraphics[width=\linewidth]{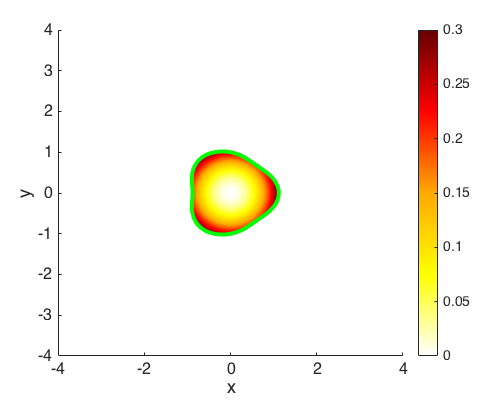}
			\caption{}\label{2dt_0}
		\end{subfigure}
		\begin{subfigure}{.49\textwidth}
			\includegraphics[width=\linewidth]{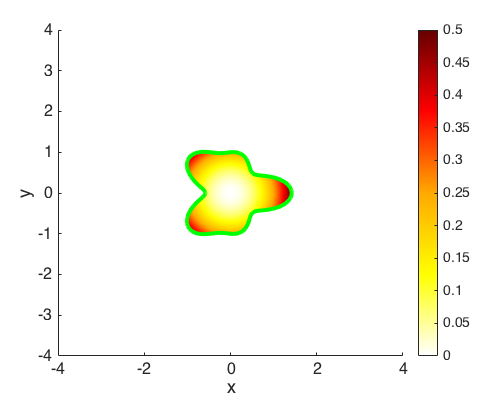}
			\caption{}\label{2dt_1}
		\end{subfigure}
		\begin{subfigure}{.49\textwidth}
			\includegraphics[width=\linewidth]{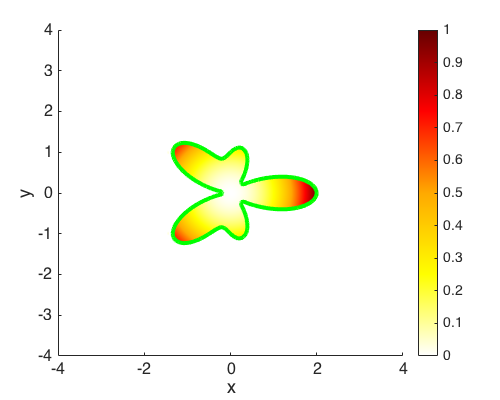}
			\caption{}\label{2dt_2}
		\end{subfigure}
		\begin{subfigure}{.49\textwidth}
			\includegraphics[width=\linewidth]{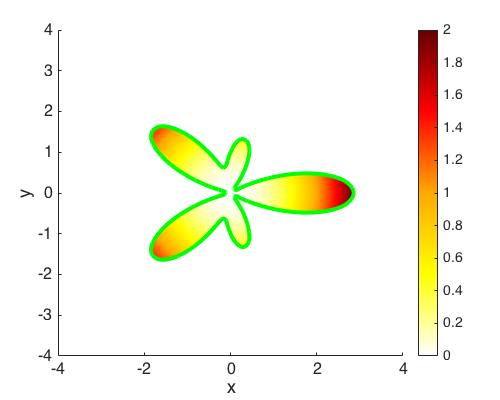}
			\caption{}\label{2dt_3}
		\end{subfigure}
		\begin{subfigure}{.49\textwidth}
			\includegraphics[width=\linewidth]{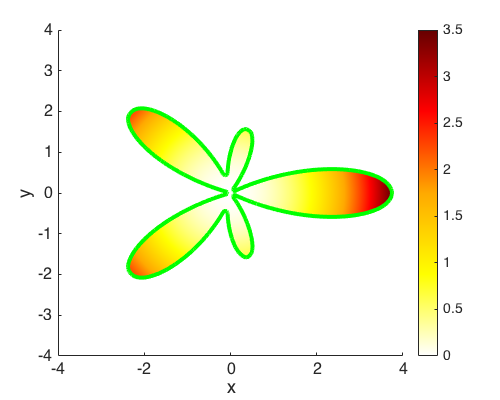}
			\caption{}\label{2dt_4}
		\end{subfigure}
		\begin{subfigure}{.49\textwidth}
			\includegraphics[width=\linewidth]{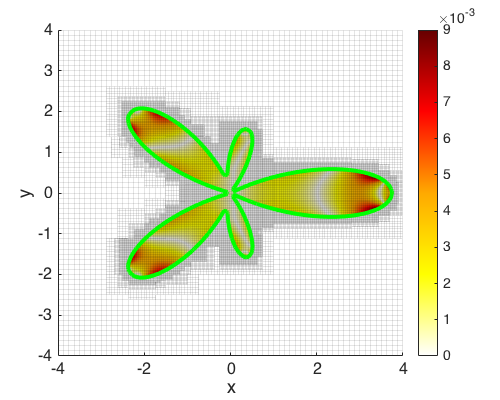}
			\caption{}\label{2dt_5}
		\end{subfigure}
		\caption{Solution of mDDMt3 from Eq. \eqref{tddm1} at different times. The solution is restricted to \(D(t)\) whose boundary $\partial D$ is denoted by the green curve. See text for details. (a): t=0, (b): t=0.5, (c): t=1.0, (d): t=1.5, (e): t=1.9. (f): The point-wise error at t=1.9 together with the block-structured adaptive mesh.} \label{2dt}
	\end{figure}

	Now we consider the 2D diffusion equation Eqs. \eqref{heat1} and \eqref{heat2} on a moving 2D domain \(D(t)\). The initial domain is enclosed by the polar curve (shown in Fig. \ref{2dt_0}):
	\begin{equation}
	r(\theta)=1+0.1\cos(3\theta)+0.02\cos(5\theta). \label{2d_po_dom}
	\end{equation}
	We suppose the velocity of the domain \(\bm{v}(\bm{x},t)\) is given by,
	\begin{equation}
	\bm{v}= (0.2(1+2t)\cos(3\theta)+0.12(1+6t)\cos(5\theta))(\cos(\theta),\sin(\theta)), \label{2d_curve_v}
	\end{equation}
	and solve for the signed distance function \(r(x,t)\) using the level-set method as described in Sec. \ref{1D time dependent}. The phase field function \(\phi\) is determined from the signed distance function \(r\) via Eq. \eqref{phi}. 
	We choose  \(f\) and \(g\) such that the exact solution is \(u(r,\theta,t)=\frac{1}{4}r^2\). Accordingly, we take the initial condition $u_\epsilon(\mathbf{x},0)=\frac{1}{4}r^2$. The extensions of $f$, $g$ off $D$ that are  constant in the normal direction are accomplished by solving \cite{peng1999pde}
	\begin{align}
	\psi_t+\text{sgn}_h(r)\frac{\nabla r}{|\nabla r|}\cdot \nabla \psi=0,
	\label{HJextension}
	\end{align}
	{\color{black} in a region near \(\partial D_t\) with width of \({O}(\epsilon)\)}, where \(\psi = f\text{ or }g\) and sgn\(_h\) is given by Eq. \eqref{sign_r}. {\color{black} Because the Hamilton-Jacobi equations (\ref{hj})-(\ref{rein}) for the signed distance function $r$ and Eq. (\ref{HJextension}) for the extensions of $f$ and $g$ need to be solved only in a narrow band around $\partial D$ and explicit time stepping methods are used, the costs associated with solving these equations are negligible compared to those associated with solving the mDDM equations.}
	
	We solve mDDMt3 (Eq. \eqref{tddm1}) on the larger domain \(\Omega=[-2,2]\times [-2,2]\) using the Crank-Nicholson method and the standard second-order central difference scheme together with adaptive, block-structured mesh refinement as in \cite{wis07}. The implicit equations are solved using a mass-conserving multigrid solver \cite{Feng2018}. {\color{black} Empirically, we found the number of iterations of the multigrid solver is insensitive to the choice of the minimum grid size $h_{fine}$ and only weakly dependent on $\epsilon$. See Tab. \ref{multigrid} in App. \ref{res_2d}.} The \(L^2\) and \(L^{\infty}\) norms of the errors are calculated using 2D extensions of Eqs. (\ref{err}) and (\ref{err_max}). We use an adaptive mesh with a maximum of three levels of refinement. The mesh is refined according to the undivided gradient of $\left(u_\epsilon\phi\right)$, e.g., if $|\nabla\left(u_\epsilon\phi \right)|\ge 10^{-4}/\left(2h\right)$, where $h$ is the local grid size, then the mesh is targeted for block-structured refinement (see \cite{wis07} for details). We start with
	\(\epsilon=0.2\) and the coarse mesh grid size \(h_{coarse}=4/16\), which gives the fine mesh size \(h_{fine}=4/128=\epsilon/6.4\), and the time step \(dt=0.1/128\). Then, the parameters \(\epsilon\), \(h_{coarse}\) and \(dt\) are refined at the same time while the number of levels and \(h_{fine}=\epsilon/6.4\) are fixed.

\begin{figure}[H]
		\center
		\begin{subfigure}{\textwidth}
			\center
			\includegraphics[scale=0.5]{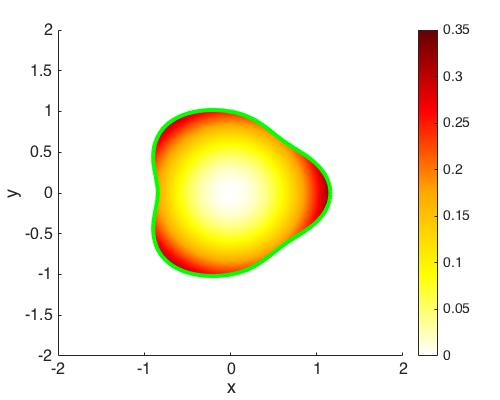}
			\caption{}\label{2d_5e}
		\end{subfigure}
		\begin{subfigure}{.49\textwidth}
			\includegraphics[width=\linewidth]{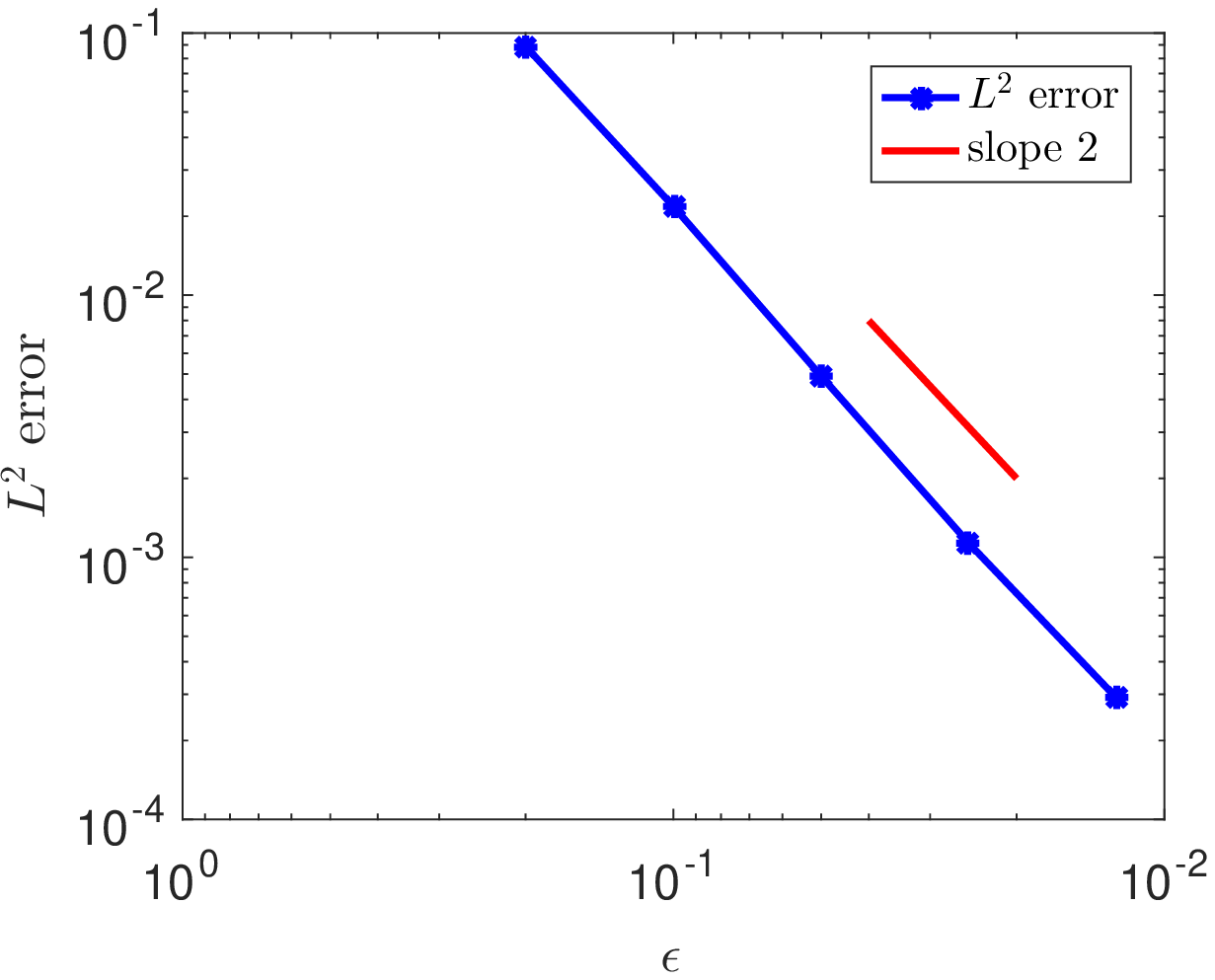}
			\caption{}\label{err_2dt_mov_2}
		\end{subfigure}
		\begin{subfigure}{.49\textwidth}
			\includegraphics[width=\linewidth]{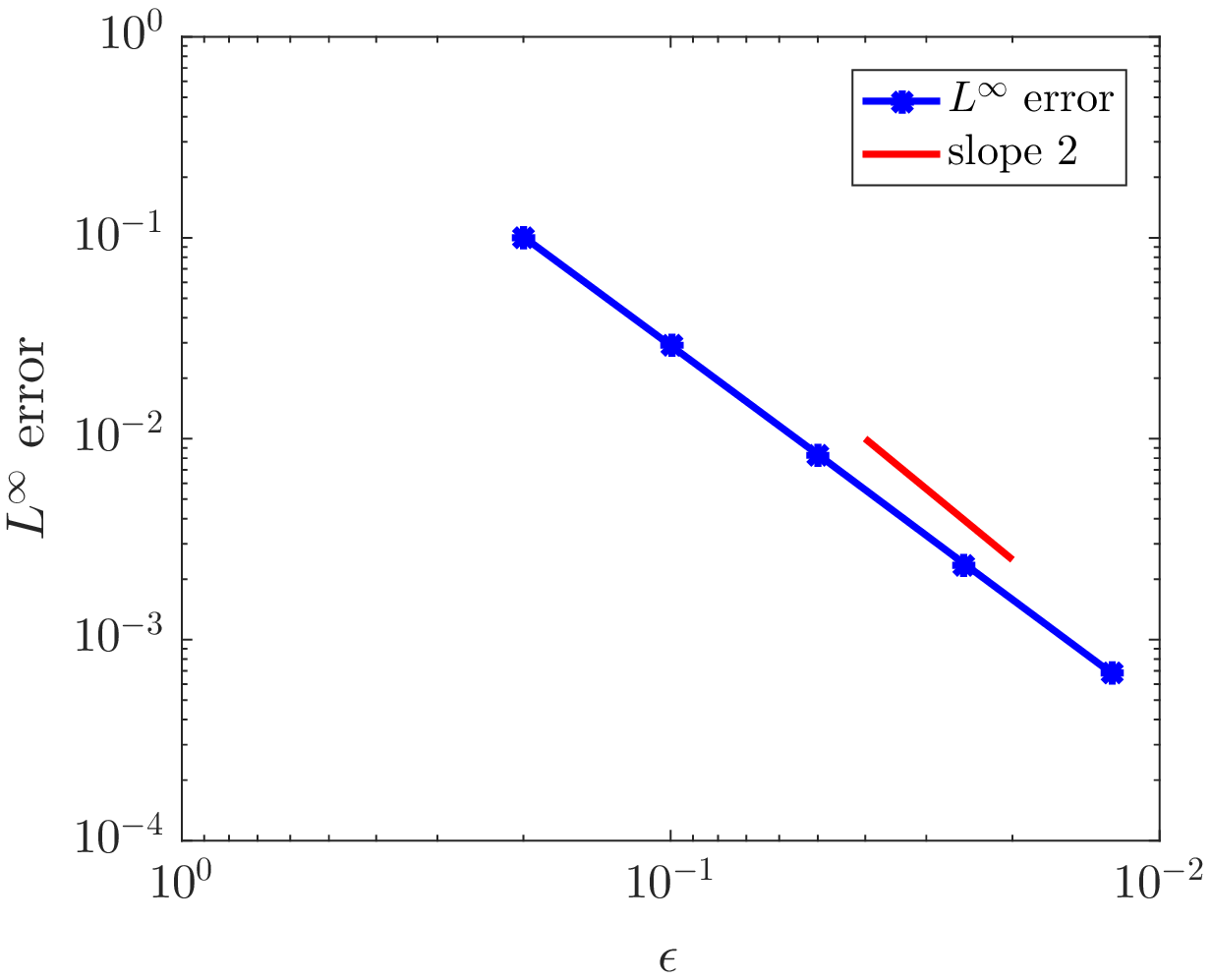}
			\caption{}\label{err_2dt_mov_m}
		\end{subfigure}
		\caption{Error analysis of the 2D diffusion equation using mDDMt3 on the moving domain \(D(t)\), see text for details. The green curve denotes the boundary $\partial D$. (a): The solution at \(t=0.1\) restricted on \(D(t)\) with the boundary contour (green line), (b): The \(L^{2}\) error at t=0.1, (c): The \(L^{\infty}\) error at t=0.1. }\label{err_2dt_mov_all}
	\end{figure}

	Note that for the 2D and 3D numerical numerical implementations, we need to discretize the normal vector \(\bm{n}\) in the mDDMs in the following way in order for the solvers to converge:
	\\[0.5cm]
	\(\bm{n}=\left\{
	\begin{array}{lr}
	-\frac{\nabla \phi}{|\nabla \phi|}, & \text{if}~r\le 0,\\
	0  & \text{else.}
	\end{array}
	\right.
	\)
	\\[0.5cm]
	{\color{black} We believe this is related to the smoothing properties of the relaxation step, but a full analysis is beyond the scope of this paper and will be considered in future work. Effectively this means that the modification in the mDDMs is performed in a region near the boundary but within the physical domain.}
The solution \(u_{\epsilon}\) in the evolving domain is shown at \(t=0.1\) using \(\epsilon=0.025\) in Fig. \ref{2d_5e}. In Figs. \ref{err_2dt_mov_2} and \ref{err_2dt_mov_m}, the \(L^2\) and \(L^{\infty}\) errors are shown at this time ({\color{black}see also Tab. \ref{err_2dt_mov} in App. \ref{res_2d}}). The results suggest mDDMt3 is 2nd order accurate in \(h\), \(\epsilon\) and \(dt\) in \(L^2\) and between 1.5 and 2nd order in \(L^{\infty}\), consistent with our theory and numerical results from the previous section.
	
	In Fig. \ref{2dt}, we present a long time simulation of the dynamics using mDDMt3 (Eq. \eqref{tddm1}) {\color{black} using the domain \(\Omega=[-4,4]\times [-4,4]\)} with  \(\epsilon=0.025\), \(h_{coarse}=8/64\), three levels of mesh refinement and \(dt=0.1/512\). {\color{black} Here the computational domain is larger than that used for the convergence test presented earlier ($[-2,2]\times[-2,2]$) to accommodate the growth of the physical domain $D$ (see Fig. \ref{2dt}). To improve efficiency, we could dynamically increased the size of $\Omega$ to accomodate the growth of $D$. In general, the buffer between $\partial D$ and $\partial \Omega$ should be larger than $\approx 10\epsilon$.}
	We present the solution restricted to the moving domain \(D(t)\) at different times, up to t=1.9. 
	The pointwise error at \(t=1.9\) is shown in Fig. \ref{2dt_5} and the largest error, which is of order \(10^{-3}\), occurs near the right most finger.
	As can be seen from these figures, mDDMt3 is able to accurately simulate solutions on time-dependent, highly complex non-rectangular domains.

{\color{black}
	\section{3D numerical results}\label{3d}	

\subsection{3D Poisson test problem}

We next consider a 3D Poisson equation on a stationary complex domain $D$: \begin{align}
-\boldsymbol{\nabla}\cdot \left( \beta(\boldsymbol{x})\boldsymbol{\nabla}u(\boldsymbol{x})\right) &= f(\boldsymbol{x}),~~~~~~~~\boldsymbol{x}~\in D,\label{3d-psn}\\
u&=p(\boldsymbol{x}). ~~~~~~~~\boldsymbol{x}~\in \partial D.\label{3d-psn-bc}
\end{align}
Following \cite{HELLRUNG20122015}, in which this equation was used as as test problem for the virtual node method,
we take $\beta(x,y,z)=7+x+2y+3z$ and set $f$ and $p$ such that the exact solution is $u(x,y,z)=x e^{y}+e^{z}\sqrt{1+y^2}$. The domain $D$ is bounded by a torus centered at $(0,0,0)$ with major radius $R_{a}=0.6$, minor radius $R_{b}=0.3$, and axis along $(0,-\sin(0.75),\cos(0.75))$. See Fig. \ref{3d-result}(left), which also shows the block-structured adaptive mesh we used.

To approximate Eqs.(\ref{3d-psn}) and (\ref{3d-psn-bc}), we use a version of mDDM3:
\begin{align}
-\boldsymbol{\nabla}\cdot \left( \phi\beta\boldsymbol{\nabla}u_{\epsilon}\right)+\cfrac{|\boldsymbol{\nabla}\phi|}{\epsilon^2} \left(u_{\epsilon}-g-r\boldsymbol{n}\cdot \boldsymbol{\nabla}u_{\epsilon}\right)&=\phi f.\label{3D-DD}
\end{align} 
The signed distance function $r$ is taken to be:
\begin{align}
r &= \sqrt{\left(\left(R_{a}-\sqrt{{r_x}^2+{r_y}^2}\right)^2+{r_z}^2\right)}-R_{b},
\end{align}
where $r_{x}=x$, $r_{y}=y\sin(\theta)+z\cos({\theta})$, $r_{z}=y\cos(\theta)-z\sin({\theta})$, and $\theta=135^{o}$.
The phase field function $\phi$ is determined from the signed distance function $r$ via Eq.~(\ref{phi}). 

We solve Eq.~(\ref{3D-DD}) on the larger domain $\Omega=[-1,1]^3$ using the standard second-order difference scheme together with adaptive, block-structured mesh refinement as in Sec. \ref{2d}. We start with $\epsilon=0.2$ and a uniform grid with mesh size $h=2/16$. We test convergence by decreasing $\epsilon$ by factors of two and concommittantly increasing the number of levels of refinement by one such that the finest mesh size $h_{fine}=\epsilon/3.2$ is fixed.  Fig. \ref{3d-result}(left) shows the torus surface (in grey) with an adaptive mesh of 2 refinement levels, which corresponds to $\epsilon=0.05$. Note that the refined mesh is concentrated near $\partial D$. The mDDM3 solution $u_{\epsilon}$ on $\partial D$ is shown in Fig.~\ref{3d-result}(right). 

The $L^{\infty}$ errors, compared with those using the second order virtual node method \cite{HELLRUNG20122015}, which like other specialized, higher-order embedded boundary methods is more difficult to implement than mDDM3, are shown in Fig.~\ref{3d-error}. The grid size $N$ for mDDM3 is obtained by assuming that we are using an uniform grid $N^3$ so that $h_{fine}=2/N$. The results suggest that the rate of convergence in $N$ of mDDM3 is similar to that of the virtual node method (about 1.8 order), although the absolute errors in mDDM3 are larger by about an order of magnitude. {\color{red} The reasons for this discrepancy are not entirely clear but are likely due to a combination of factors including the choices of the relation between $h_{fine}$ and $\epsilon$ and the scaling of the $|\nabla\phi|$ term in the boundary condition term Eq. (\ref{3d-psn-bc}). In particular, by varying these choices, there are many ways the error could be optimized. For example, for each $h_{fine}$, we could find an optimal choice of $\epsilon$ that minimizes the error. While we tested a few combinations that could reduce the error, the optimal relation between $h_{fine}$ and $\epsilon$ to minimize the error is not clear. However, taking $\epsilon \propto h_{fine}$ as we did robustly yields nearly 2nd order convergence rates. In addition, following \cite{Poulson-2018}, one could introduce an additional parameter $\alpha$ in the boundary condition term in Eq. (\ref{3d-psn-bc}), e.g., replace $|\nabla\phi|$ with $\alpha|\nabla\phi|$, and try to minimize the error through judicious choices of $\alpha$. Other aspects of the solver (e.g., distances over which extensions are performed, adaptive mesh parameters, etc.) could also be optimized. These are subjects for future work. 
}
 
	\begin{figure}[H]
		\center
		\begin{subfigure}{.44\textwidth}
			\includegraphics[width=\linewidth]{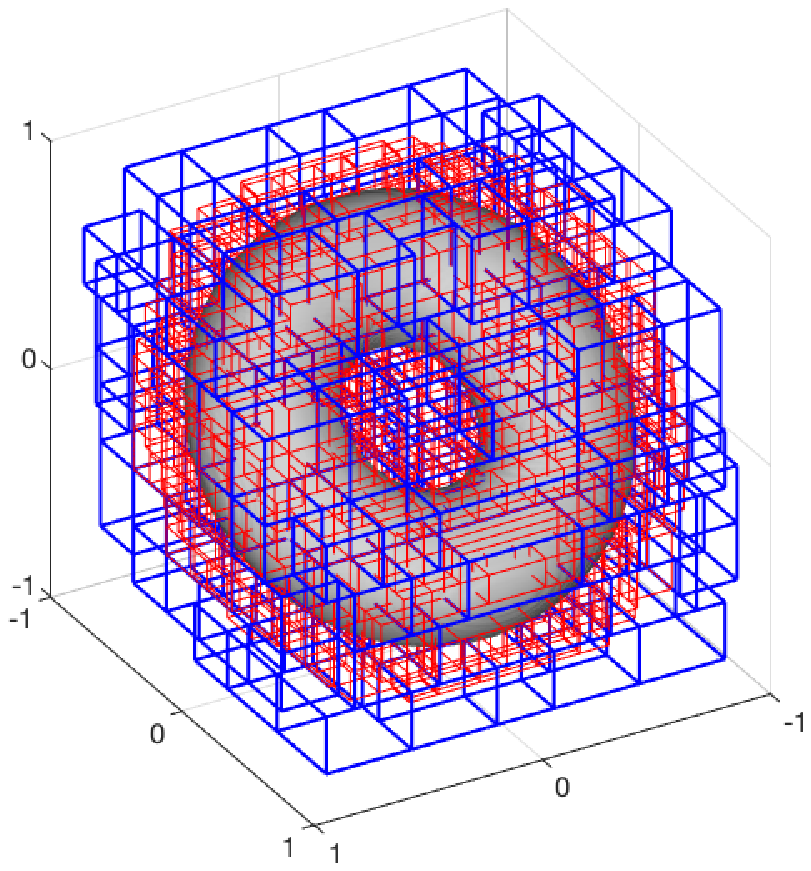}
		\end{subfigure}
		\begin{subfigure}{.44\textwidth}
			\includegraphics[width=\linewidth]{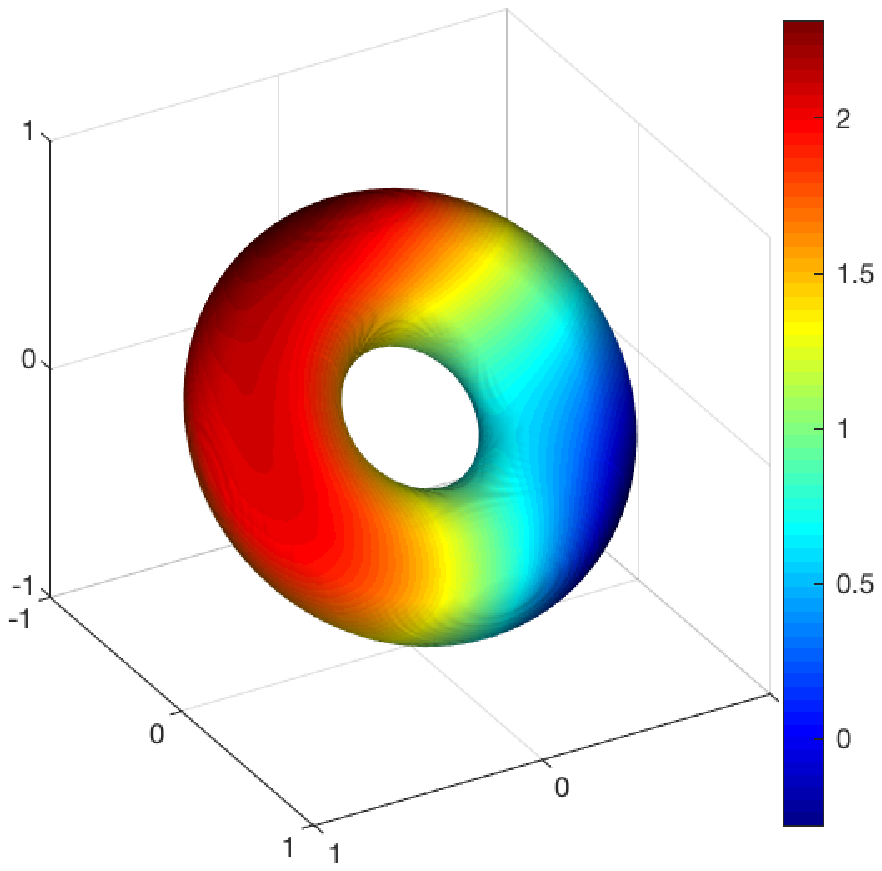}
		\end{subfigure}
		\caption{A 3D example of the Poisson equation in a complex domain $D$ using mDDM3. (left): The boundary of $D$ is the toroidal surface (in grey), which is shown together with a 2-level adaptive mesh that corresponds to the interface thickness $\epsilon=0.05$, (right): the mDDM3 solution $u_{\epsilon}$ at the boundary $\partial D$. See text for details.}
		\label{3d-result}
	\end{figure}

	\begin{figure}[H]
		\center
			\includegraphics[width=.6\linewidth]{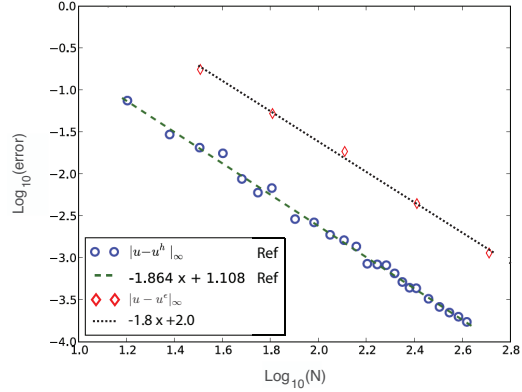}
		\caption{$L^\infty$ errors in solving the 3D Poisson equation on the toroidal domain $D$ from Fig. \ref{3d-result}. The open circles correspond to results from \cite{HELLRUNG20122015} using the virtual node method while the diamonds correspond to results using mDDM3. $N$ is the number of grid points in one dimension for the virtual node method, which uses a uniform mesh, while $N=2/h_{fine}$ for mDDM3 in which an adaptive mesh is used. See text for details.} \label{3d-error}
	\end{figure}

{\color{red}
\section{Simulation of drug perfusion in the human brain}

As a final example, we simulate the perfusion of a drug in a 3D model of a human brain. In particular, we used an anatomical template of the human brain from the ICBM atlas \cite{Atlas,Collins1999,Fonov2009,Fonov2011}.  The atlas provides smooth, probabilistic segmentations for the anatomical structures in the brain.This data is given on a Cartesian $256^3$ voxel grid $\Omega_h$ where the volume of each voxel is 1 $mm^3$. Each voxel in this grid can contain white matter (WM), grey matter (GM), cerebral spinal fluid (CSF) or no brain matter (empty space).

We thresholded the CSF such that in the interior brain voxels contain only CSF or solid tissue (GM/WM) and only those voxels near CSF/solid tissue boundaries contain a combination of CSF and solid tissue. The corresponding volume fractions of the WM and GM are $\phi_W$ and $\phi_G$, which are smooth functions. The CSF region is inferred as the region where $\phi_W+\phi_G< 0.5$ in the brain interior. There is also a thin rim of CSF at the outer boundary of the brain. See Fig. \ref{3d-brain-result} (leftmost column) where the regions from lightest to darkest denote the WM (where $\phi_W\ge 0.5$), the GM (where $\phi_G\ge 0.5$) and CSF. Denote these regions by $D_{W}$, $D_{G}$, and $D_{CSF}$. 

Let $u_W$ and $u_G$ be the drug concentrations in $D_{W}$ and $D_{G}$, respectively. As a simple model, we assume that the drug is delivered through the CSF, which is modeled as a Dirichlet boundary condition for the drug concentration at the boundary between the solid tissue and the CSF regions. Therefore, we model the perfusion of drug by:
\begin{align}
\partial_t u_{W}(\mathbf{x},t)-\boldsymbol{\nabla}\cdot \left(\beta_{W}\boldsymbol{\nabla}u_W \right) +\alpha u_W&=0,~~~~~~~~\mathbf{x}~\in D_W,\label{3d-brain-1}\\
\partial_t u_{G}(\mathbf{x},t)-\boldsymbol{\nabla}\cdot \left(\beta_{G}\boldsymbol{\nabla}u_G \right) +\alpha u_G&=0,~~~~~~~~\mathbf{x}~\in D_G,\label{3d-brain-2}\\
u_W(\mathbf{x},t)=u_G(\mathbf{x},t)&=g(\mathbf{x},t), ~~~~~~~~\mathbf{x}\in \partial D_{CSF},\label{3d-brain-bc}\\
u_W(\mathbf{x},t)=u_G(\mathbf{x},t),& ~~~~~~~~\mathbf{x}\in \partial D_{G}\cap\partial D_{W},\label{3d-brain-bc-2}\\
\beta_W\nabla u_W\cdot\mathbf{n}(\mathbf{x},t)&=\beta_G\nabla u_G\cdot\mathbf{n}(\mathbf{x},t), ~~~~~~~~\mathbf{x}\in \partial D_{G}\cap\partial D_{W},\label{3d-brain-bc-3}
\end{align}
where $\beta_W$, $\beta_G$ are the diffusion coefficients in the WM and GM, respectively, and for simplicity we have assumed that the uptake rate of drug, $\alpha$, is the same in WM and GM.

To approximate Eqs.(\ref{3d-brain-1})-(\ref{3d-brain-bc-3}), we first nondimensionalize using the length scale $L=64~mm$ and time scale $L^2/\bar \beta$, where $\bar \beta$ is a characteristic diffusion coefficient. Then, we formulate a nondimensional diffuse domain model using a time-dependent version of mDDM1: 
\begin{align}
\left(\phi_{S}u_{\epsilon}\right)_{t}-\boldsymbol{\nabla}\cdot \left(\phi_S\beta\boldsymbol{\nabla}u _{\epsilon}\right) &+\alpha\phi_{S} u_{\epsilon}+\cfrac{1-\phi_{S}}{\epsilon^3} \left(u_{\epsilon}-g-r\boldsymbol{n}\cdot \boldsymbol{\nabla}u_{\epsilon}\right)=0,\label{3d-brain-dd}
\end{align} 
where $\beta=\beta_W\phi_W+\beta_G\phi_G$ is a nondimensional diffusion coefficient and $\phi_S$ is the hyperbolic tangent function in Eq. (\ref{phi}) using the signed-distance function $r$ obtained by solving Eq. (\ref{rein}) using the initial condition
\begin{align}
r(\mathbf{x},0) &= \phi_{W}+\phi_{G}-0.5.
\end{align} 
Note that this formulation also utilizes a diffuse domain approximation of the continuity of concentration and flux at the WM/GM interfaces \cite{li09,Lervag2015}. We solve the system using a 3D version of the solver in Sec. \ref{2d} on the uniform $256^3$ grid (e.g., the nondimensional spatial grid size is $h=4/256$). We take $g=1$, $\epsilon=0.03\approx 2h$, $\alpha=0.005$ and time step $\Delta t= 0.001$. The results are shown in Fig. \ref{3d-brain-result}, which depicts the evolution of the concentration field in the solid tissue $\phi_s u_\epsilon$ using different diffusion coefficient ratios $\beta_W/\beta_G$. 

As time evolves, the extent of the drug penetration into the brain tissue depends on $\beta_W/\beta_G$ and by time $t=0.05$, the drug concentration field has equilibrated. In Fig. \ref{3d-brain-result}(a), where $\beta_W=1$ and $\beta_G=0.1$, the drug is mainly confined to regions of the brain around the CSF since most of the CSF is surrounded by GM where the diffusion coefficient is small. For this reason in Fig. \ref{3d-brain-result}(b), where the diffusivity is larger in GM ($\beta_W=0.1$ and $\beta_G=1.0$), the drug penetrates much farther into the brain tissue. The most drug perfusion is observed in Fig. \ref{3d-brain-result}(c) where the diffusivities are matched ($\beta_G=\beta_W=1$). While we considered the diffusivities to be isotropic in the brain tissue, the presence of fibers in the WM can introduce anisotropic diffusion. This can be considered in future work. 

	\begin{figure}[H]
		\center
		\begin{subfigure}{.8\textwidth}
			\includegraphics[width=\linewidth]{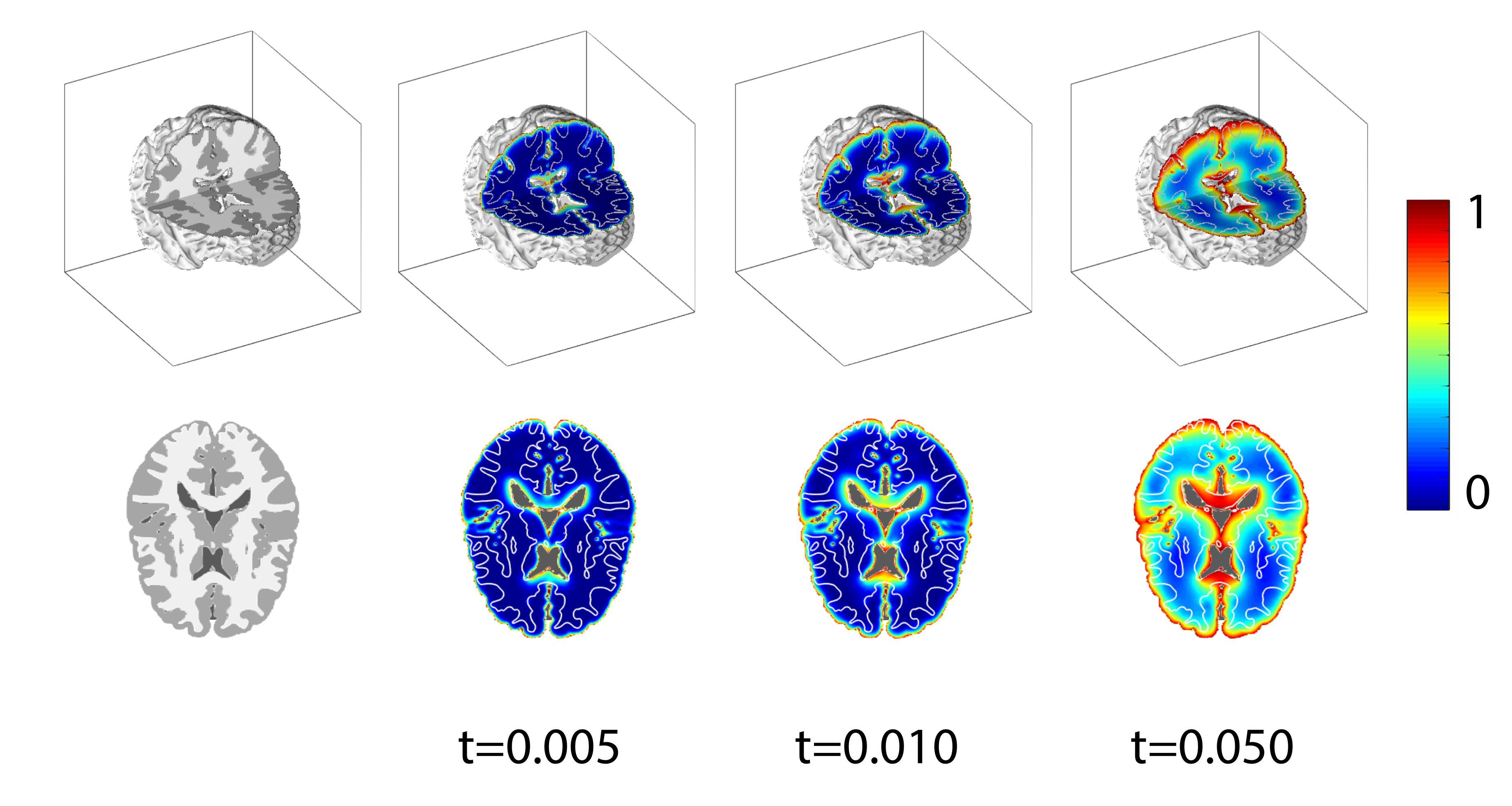}(a)
		\end{subfigure}
		\begin{subfigure}{.8\textwidth}
			\includegraphics[width=\linewidth]{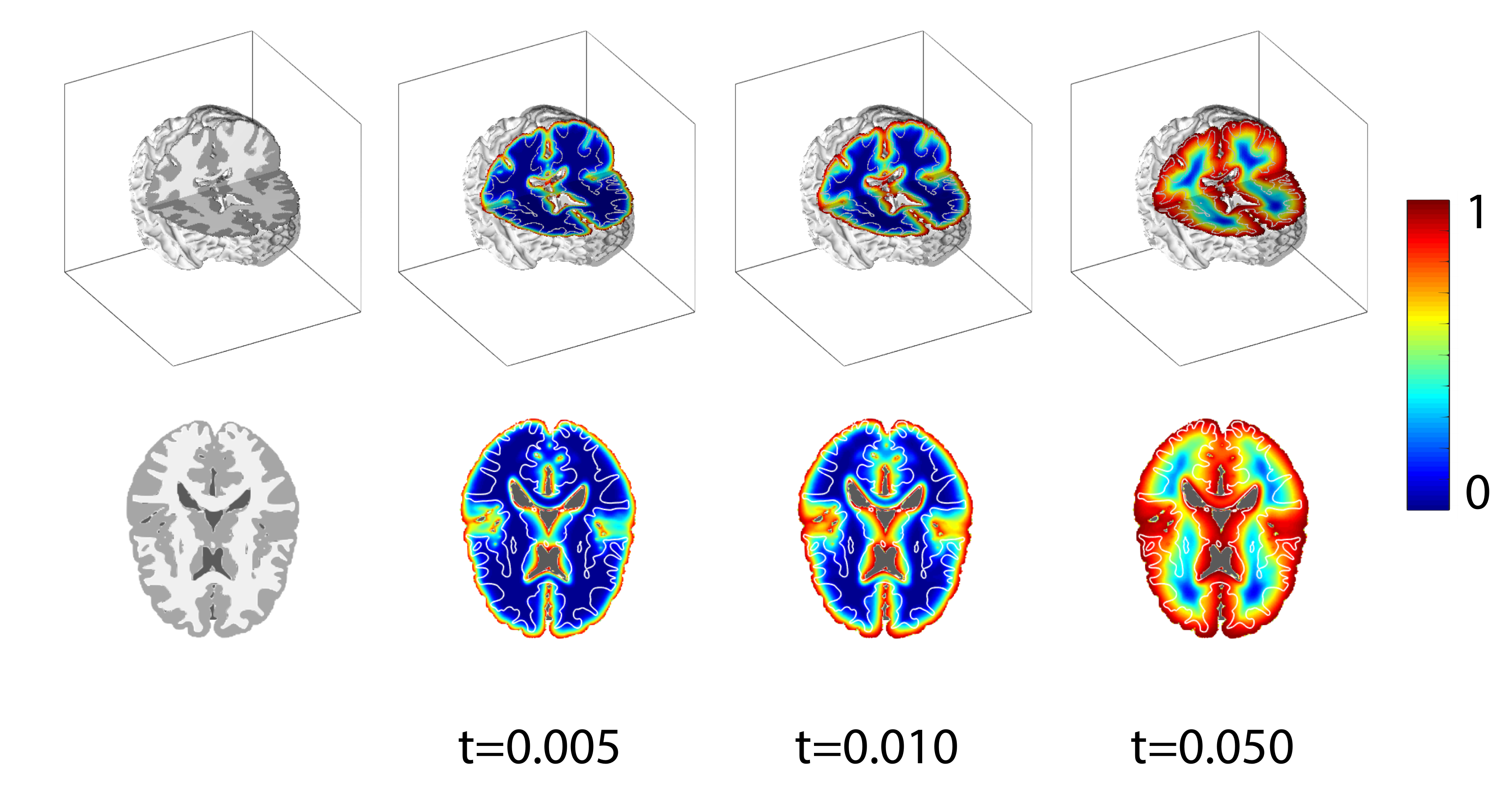}(b)
		\end{subfigure}
		\begin{subfigure}{.8\textwidth}
			\includegraphics[width=\linewidth]{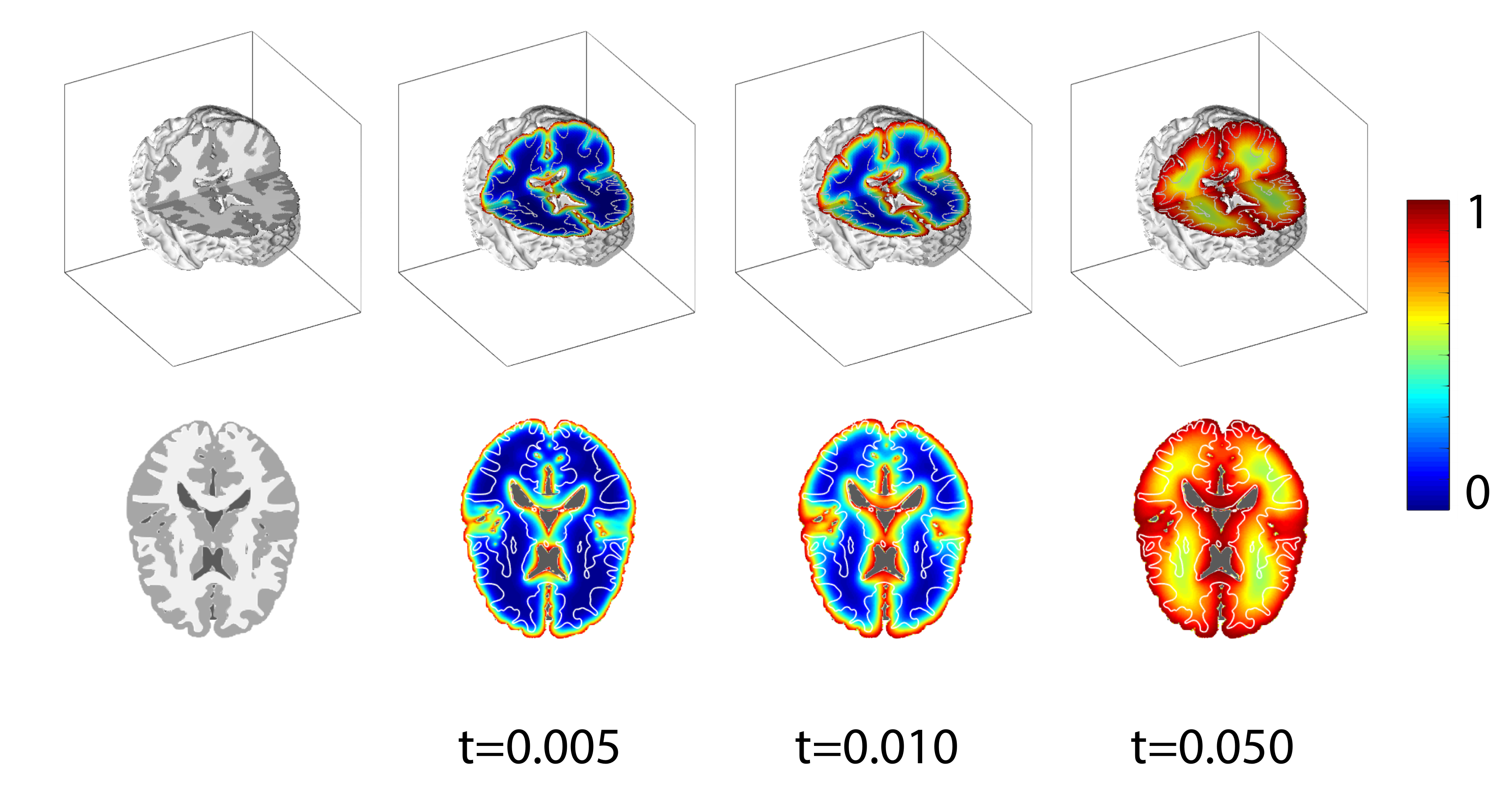}(c)
		\end{subfigure}
		\caption{Drug perfusion in an anatomical template for the human brain from the ICBM atlas \cite{Atlas} (see text). Leftmost column: Light gray (white matter, WM), darker grey (grey matter, GM), darkest grey (cerebral spinal fluid, CSF). A Dirichlet boundary condition is posed at the boundary of the CSF region. The evolution of the drug concentration is shown over time for different ratios of drug diffusivities $\beta_W$ and $\beta_G$ in the WM and GM, respectively. (a): $\beta_{W}/\beta_{G}=10$, (b): $\beta_{W}/\beta_{G}=0.1$, (c): $\beta_{W}/\beta_G=1$.} 
		\label{3d-brain-result}
	\end{figure}
}

}

	\section{Conclusions and future work}\label{conclude}
	
	We have analyzed several diffuse domain methods (DDMs) originally developed in \cite{li09} for both the Poisson equation and the diffusion equation with Dirchlet boundary conditions on stationary and moving domains. {\color{black} Advantages of DDM methods is that they are flexible and easy to implement. The same solver can be used for any domain, whether moving or stationary, and that the methods can be implemented using the tools provided in standard software packages so that one does not need to write specialized code that requires treating the region near the interface differently from that in bulk regions.}
	
	Our analysis reveals why different DDM formulations yield varying degrees of accuracy. Guided by our analysis, we presented new modifications of the DDMs (mDDMs) that provide higher-order accuracy. {\color{black} While the DDM and mDDM systems require the solution of non-constant coefficient equations that introduce extra length scales to the problem, iterative matrix solvers available in most software packages are sufficient to solve the discrete equations. Here, we used a mass-conservative, finite-difference multigrid method for which the number of iterations is not sensitive to the smoothing thickness $\epsilon$ or the local grid size $h$.}
		
	Using a matched asymptotic analysis, two methods, mDDM1 and mDDM3, were shown to be 2nd order accurate in $L^2$, e.g., $O(\epsilon^{2})$ where $\epsilon$ is the diffuse interface smoothing parameter. The analysis shows that the errors in the $L^\infty$ norm scale as $O(\epsilon^{1.5})$. {\color{black}Using numerical simulations in 1D, 2D and 3D for selected test cases}, where
the level-set method was used to implicitly capture the domain movement, to construct the smoothed characteristic function needed by the methods and to perform the needed constant normal extensions, these theoretical predictions were confirmed when the grid size $h\sim \epsilon^{1.5}$ is used to resolve the boundary layer (e.g., $K_3$ in Sec. \ref{aa_mddm1_3}). 
	
	In addition, numerical simulations revealed that in the $L^2$ norm, mDDM1 and mDDM3 and their time-dependent versions (mDDMt1, mDDMt3) are $O(\epsilon^{2})$ even when $h\propto \epsilon$, which is much cheaper computationally and hence more cost-effective to use for 2D and 3D problems. In this case, however, the $L^\infty$ errors were found to be more variable over the range of parameters tested. In $L^\infty$, simulations with \(h=\epsilon/c\) show that when $c$ is small, the schemes are roughly $O(\epsilon^2)$ whereas when $c$ is large the errors are roughly $O(\epsilon^{1.5})$. {\color{black} This occurs because when $c$ is small, the truncation and analytic errors are close in magnitude but of opposite signs so that there is cancellation in the total error and the total error is smaller than would be expected. However, as $h$ decreases, the analytic errors dominate and the total errors scale as ${O}(\epsilon^{1.5})$ as predicted by theory. 
	
	{\color{red} We also compared the performance of mDDM3 in 3D to the virtual node method in \cite{HELLRUNG20122015}, which is more difficult to implement than mDDM3. We found that mDDM3 and the virtual node method converge to the exact solution at a similar rate in $L^\infty$ but that the error is larger using mDDM3 for the same minimum grid size $h_{min}$. While the reasons for this are unclear, we identified ways that the accuracy of mDDM3 could be further improved. However, it is open question as to whether the mDDM3 could be made as accurate as the virtual node method.
	
	 As a final example in 3D, we simulated the dynamics of drug perfusion in an anatomical template of the human brain. Although the geometry is highly complex, no changes in the algorithm are needed. Assuming the drug is administered through the cerebral spinal fluid, which is modeled using a Dirichlet boundary condition, we investigated the effect of different drug diffusivities in the white and grey matter regions of the brain on the spatiotemporal drug distributions.
	 }
		
	Although we have focused on the Poisson and diffusion equations here, the mDDMs can be applied to general elliptic and parabolic partial differential equations following previous work (e.g., \cite{li09,Yu-2012}). For example, consider the following general partial differential equation with Dirichlet boundary conditions in a moving domain \(D(t)\),
	\begin{eqnarray}
	\partial_t u-\nabla \cdot(\bm{A}\nabla u)+\bm{b}\cdot \nabla u+cu=f \text{ in } D(t), 
	\label{eq:general1}
	\\
	u=g \text{ on } \partial D(t),
	\label{eq:general2}
	\end{eqnarray}
	with \(\bm{A}=\bm{A}(u,\nabla u,x,t)\) a positive definite matrix, \(\bm{b}=\bm{b}(u,\nabla u,x,t)\) a vector, \(c=c(u,\nabla u,x,t)\in \mathbb{R}\) and \(f=f(x,t)\).
	Then, the corresponding version of mDDMt3 is given by,
	\begin{equation}
	\partial_t(\phi u)-\nabla\cdot (\phi \bm{A}\nabla u)+\phi \bm{b}\cdot \nabla u+\phi c u+ \frac{|\nabla \phi|}{\epsilon^2}(u-g-r\bm{n}\cdot \nabla u)=\phi f \text{ in } \Omega.
	\end{equation}
	with \(\bm{A}, \bm{b}\) and \(c\) extended coefficients. The other mDDMs can be defined analogously.
	
	While the matched asymptotic analysis was presented here in for the Poisson and diffusion equations 1D, the same results hold in 2D and 3D for arbitrary smooth domains and for general PDEs such as in Eqs. (\ref{eq:general1})-(\ref{eq:general2}). This is because in 2D and 3D, the boundary appears flat at leading order and the effects of curvature appear in lower order terms that ultimately do not affect the expansions to the orders we calculated here, {\color{black} provided that the interface thickness $\epsilon$ is sufficiently small, e.g., smaller than any physically relevant length scale, such as the radius of curvature of the boundary, and the region around the boundary over which functions can be smoothly extended out of the physical domain.} The additional terms in the equations (advection and reaction) also do not affect the expansions we calculated here because they are lower order as well.
		
	In future work, we plan to develop a theoretical proof of convergence for the mDDMs. In addition, we plan to apply the methods developed here to important physical applications including fluid flows in complex geometries, which involves applying the theory developed here to the Navier-Stokes equations, and to the epitaxial growth of graphene heterostructures, which involves solving diffusion equations in domains that evolve according to flux balance conditions determined from the solutions to the diffusion equations themselves.
	\begin{appendices}
	
		\section{Additional numerical results of DDMs for 1D time-independent problems}\label{app_a}
		We have tested the following three DDMs using \(h=\epsilon/4\):
		
		\textbf{DDM1}: \(\nabla(\phi \nabla u)-\frac{1}{\epsilon^3}(1-\phi)(u-g)=\phi f\),
		
		\textbf{DDM2}: \( \phi\Delta u-\frac{1}{\epsilon^2}(1-\phi)(u-g)=\phi f\),
		
		\textbf{DDM3}: \(\nabla(\phi \nabla u)-\frac{1}{\epsilon^2}|\nabla \phi|(u-g)=\phi f\),
		
		\noindent on the following seven cases,
		
		Case1: \(u=\frac{x^2}{2}\);
		
		Case2: \(u=(x^2-1.111)^2\);
		
		Case3: \(u=\frac{1}{x^2+1}\);
		
		Case4: \(u=\cos(x)\);
		
		Case5: \(u=(x^2+1)^2\);
		
		Case6: \(u=\log(x^2+1)\);
		
		Case7: \(u=\sqrt{x^2+1}\).
		
		\noindent The original Poisson equation with Dirichlet boundary condition is again defined on \([-1.111,1.111]\) and we solve the DDMs on a larger domain \([-2,2]\). Tables \ref{afd1} - \ref{afd6} show the \(L^2\) and \(L^\infty\) norms of errors for the three DDMs for the seven cases. In cases 3-7, the convergence orders of DDM1-3 in both \(L^2\) and \(L^\infty\) are similar to those in case 1. However in case 2, where the derivative of the exact solution is 0 at the boundary (\(A=0\)), we observe 2nd order convergence in DDM2 and higher than 1st order convergence in DDM1 and DDM3. Our analysis in Sec. \ref{sec_ana} indicates that \(\bar u_1^{(1)}(0) = \frac{A}{3}(-\ln 6 + \gamma)\), see Eq. \eqref{ddm2_out_sol2}, \(-\frac{A}{6}\ln\epsilon+\frac{A}{3}(-\ln 6+\gamma)\), see Eq. \eqref{ddm1_out_u2}, and \(-\frac{A}{6}\ln\epsilon+\frac{A}{3}(-\ln(\frac{1}{2}\sqrt{\frac{2}{3}})+\gamma)\), see Eq. \eqref{ddm3_out_u2} for DDM2, DDM1 and DDM3, respectively. Note that these all vanish when \(A=0\). Hence, DDMs can achieve higher than 1st order accuracy when \(A=0\). In fact, the errors in both norms for DDM1 and DDM3 are dominated by \({O}(\epsilon^2(\ln(\epsilon))^2)\). This term rises in the next order matching in \(K_2\), that is \(\lim_{z_2\rightarrow -\infty} \hat u_2^{(2)} = \bar u_2^{(2)}(0) + (z_2+\ln \epsilon/6)^2\frac{d^2}{dx^2}\bar u_1^{(0)}(0)+...\). Using an analogous argument to that given in Sec. \ref{ddm13}, we obtain \(\bar u_2^{(2)}(0)\sim {O}((\ln\epsilon)^2)\). In Fig. \ref{d_eps_ddm1}, we plot \(D(\epsilon) = \lim_{z_2\rightarrow -\infty} (u_\epsilon-u)/\epsilon^2\) versus \(\ln(\epsilon)\) using the numerical solution of DDM1 in case 2 and find that it is a quadratic function of \(\ln(\epsilon)\). 
		
		\begin{figure}[H]
			\center
			\begin{subfigure}{.49\textwidth}
				\includegraphics[width=\linewidth]{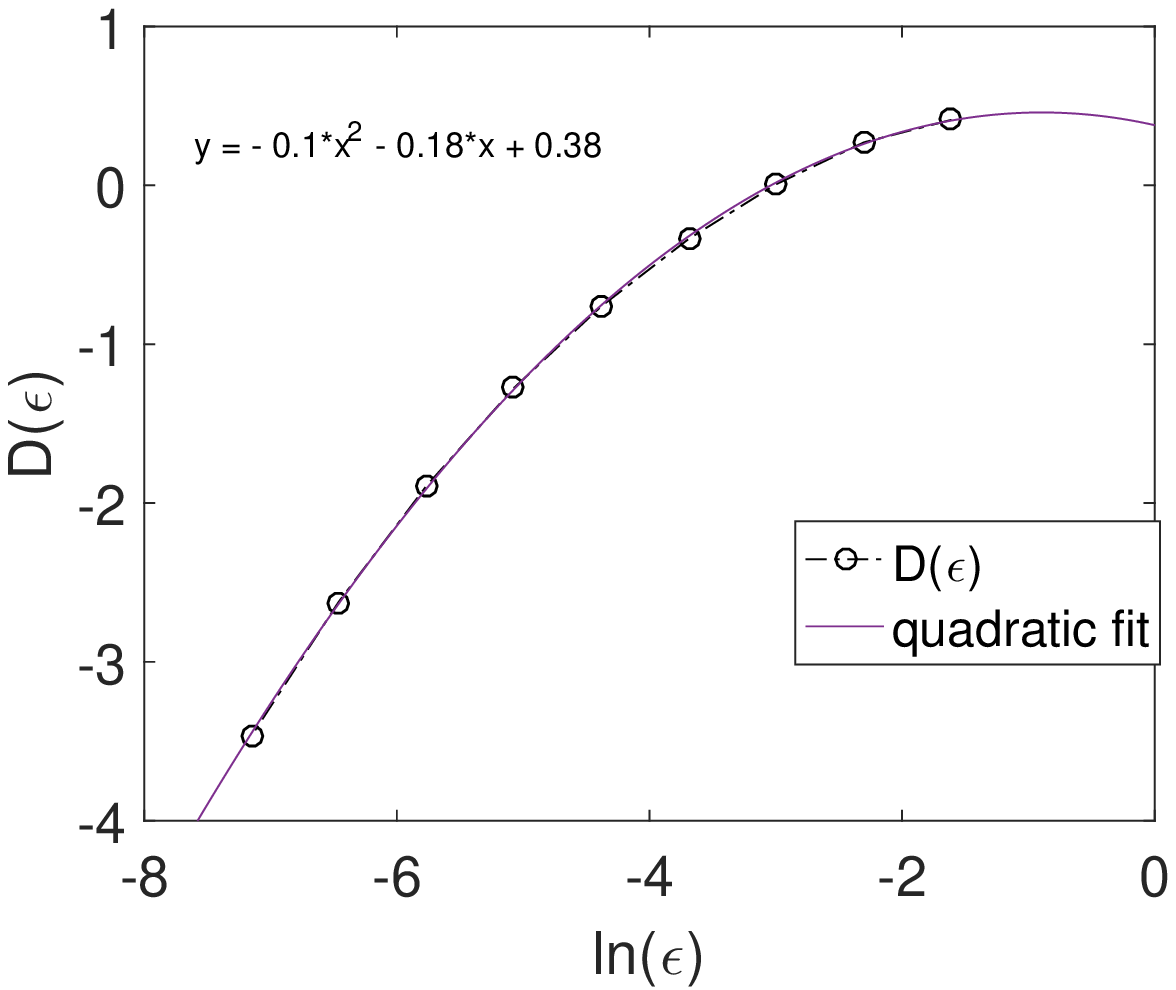}
				\caption{}\label{d_eps_ddm1}
			\end{subfigure}
			\begin{subfigure}{.49\textwidth}
				\includegraphics[width=\linewidth]{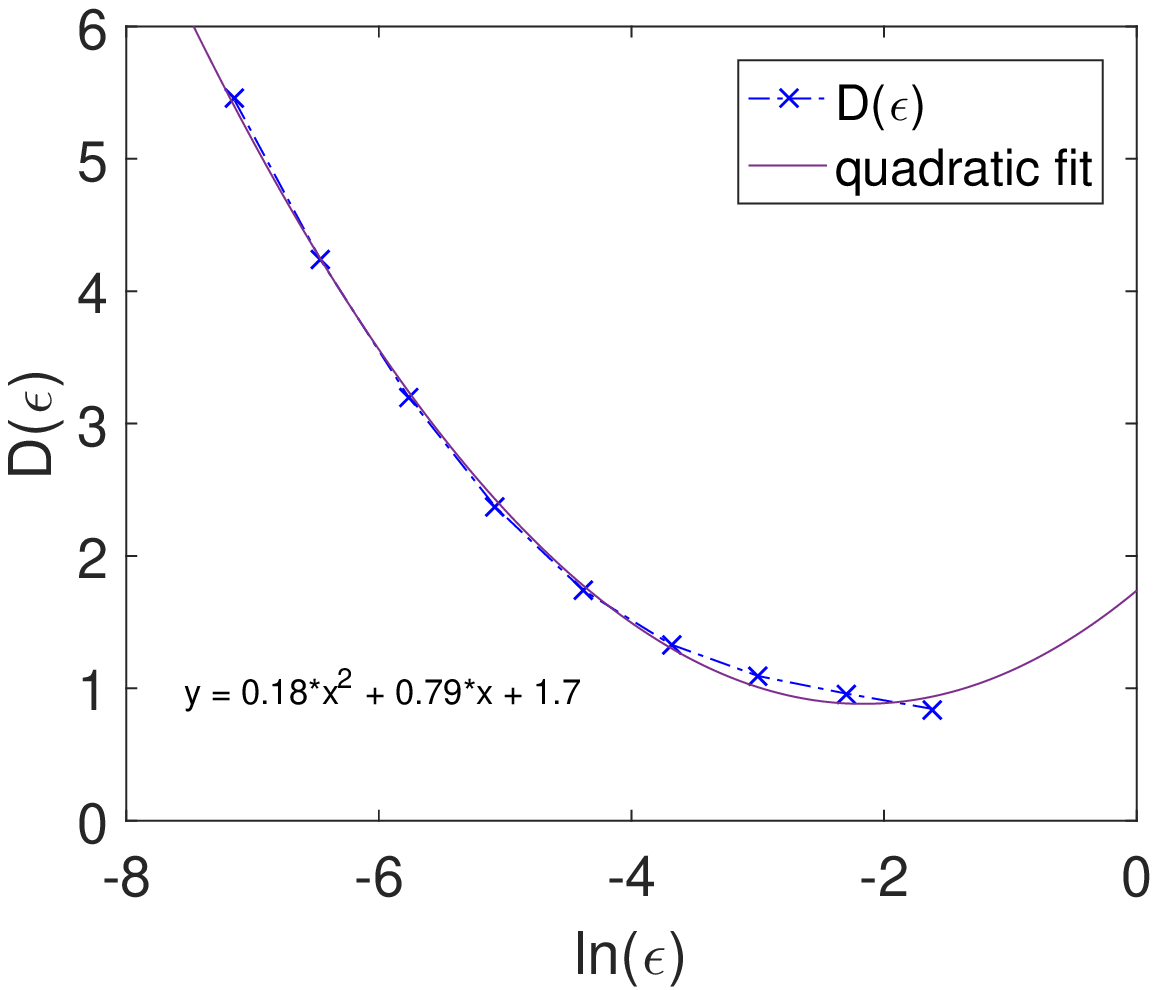}
				\caption{}\label{d_eps_mddm1}
			\end{subfigure}
			\caption{(a): \(D(\epsilon)\) for DDM1 in case 2, (b): \(D(\epsilon)\) for mDDM1 in case 2.}\label{d_eps}
		\end{figure}

		\begin{table}[H]
			\center
			
			\begin{tabular}{c|c|c|c|c|c|c}
				\bfseries \(\epsilon\) & \bfseries case 1 &\bfseries k  & \bfseries case 2 &\bfseries k
				& \bfseries case 3 &\bfseries k
				\begin{filecontents*}{ddm11.csv}
eps,case 1,k,case 2,k,case 3,k,case 4,k,case 5,k,case 6,k,case 7,k
2.00E-01,4.88E-01,0.00 ,1.55E-02,0.00 ,6.79E-02,0.00 ,1.25E-01,0.00 ,4.88E-01,0.00 ,3.03E-01,0.00 ,7.32E-02,0.00 
1.00E-01,1.12E-01,2.12 ,2.75E-03,2.50 ,1.59E-02,2.10 ,2.95E-02,2.08 ,1.13E-01,2.10 ,6.97E-02,2.12 ,1.72E-02,2.09 
5.00E-02,9.88E-03,3.51 ,9.39E-05,4.87 ,1.40E-03,3.51 ,2.63E-03,3.49 ,1.01E-02,3.49 ,6.09E-03,3.52 ,1.52E-03,3.50 
2.50E-02,1.31E-02,-0.41 ,1.93E-04,-1.04 ,1.87E-03,-0.43 ,3.52E-03,-0.42 ,1.33E-02,-0.40 ,8.11E-03,-0.41 ,2.03E-03,-0.42 
1.25E-02,1.41E-02,-0.11 ,1.15E-04,0.75 ,2.02E-03,-0.10 ,3.80E-03,-0.11 ,1.44E-02,-0.11 ,8.72E-03,-0.10 ,2.19E-03,-0.11 
6.25E-03,1.04E-02,0.44 ,4.94E-05,1.22 ,1.49E-03,0.44 ,2.81E-03,0.43 ,1.07E-02,0.43 ,6.44E-03,0.44 ,1.62E-03,0.43 
3.13E-03,6.80E-03,0.62 ,1.85E-05,1.41 ,9.70E-04,0.62 ,1.83E-03,0.62 ,6.95E-03,0.62 ,4.19E-03,0.62 ,1.06E-03,0.62 
1.56E-03,4.16E-03,0.71 ,6.46E-06,1.52 ,5.94E-04,0.71 ,1.12E-03,0.71 ,4.25E-03,0.71 ,2.57E-03,0.71 ,6.46E-04,0.71 
\end{filecontents*}
\csvreader{ddm11.csv}{}
				{\\\hline \csvcoli & \csvcolii & \csvcoliii & \csvcoliv & \csvcolv& \csvcolvi
					& \csvcolvii}
			\end{tabular}
			\begin{tabular}{c|c|c|c|c|c|c|c}
				\bfseries case 4 &\bfseries k
				& \bfseries case 5 &\bfseries k
				& \bfseries case 6 &\bfseries k
				& \bfseries case 7 &\bfseries k
				\csvreader{ddm11.csv}{}
				{\\\hline  \csvcolviii& \csvcolix & \csvcolx& \csvcolxi& \csvcolxii& \csvcolxiii
					& \csvcolxiv& \csvcolxv}
			\end{tabular}
			\caption{The \(L^2\) errors for DDM1}\label{afd1}
		\end{table}
		
		\begin{table}[H]
			\center
			
			\begin{tabular}{c|c|c|c|c|c|c}
				\bfseries \(\epsilon\) & \bfseries case 1 &\bfseries k  & \bfseries case 2 &\bfseries k
				& \bfseries case 3 &\bfseries k
				\begin{filecontents*}{ddm11inf.csv}
eps,case 1,k,case 2,k,case 3,k,case 4,k,case 5,k,case 6,k,case 7,k
2.00E-01,2.12E-01,0.00 ,1.09E-02,0.00 ,5.36E-02,0.00 ,1.04E-01,0.00 ,2.32E-01,0.00 ,1.47E-01,0.00 ,5.82E-02,0.00 
1.00E-01,5.26E-02,2.01 ,1.93E-03,2.49 ,1.28E-02,2.07 ,2.57E-02,2.02 ,5.81E-02,2.00 ,3.60E-02,2.03 ,1.43E-02,2.02 
5.00E-02,7.96E-03,2.72 ,8.82E-05,4.45 ,1.93E-03,2.73 ,3.88E-03,2.72 ,8.83E-03,2.72 ,5.43E-03,2.73 ,2.17E-03,2.72 
2.50E-02,5.89E-03,0.44 ,1.42E-04,-0.69 ,1.46E-03,0.40 ,2.91E-03,0.41 ,6.48E-03,0.45 ,4.06E-03,0.42 ,1.62E-03,0.42 
1.25E-02,6.31E-03,-0.10 ,7.86E-05,0.86 ,1.57E-03,-0.10 ,3.14E-03,-0.11 ,6.96E-03,-0.10 ,4.35E-03,-0.10 ,1.74E-03,-0.10 
6.25E-03,4.66E-03,0.44 ,3.28E-05,1.26 ,1.16E-03,0.44 ,2.32E-03,0.43 ,5.14E-03,0.44 ,3.21E-03,0.44 ,1.29E-03,0.44 
3.13E-03,3.04E-03,0.62 ,1.22E-05,1.43 ,7.53E-04,0.62 ,1.51E-03,0.62 ,3.35E-03,0.62 ,2.09E-03,0.62 ,8.40E-04,0.62 
1.56E-03,1.86E-03,0.71 ,4.20E-06,1.53 ,4.60E-04,0.71 ,9.26E-04,0.71 ,2.05E-03,0.71 ,1.28E-03,0.71 ,5.14E-04,0.71 
\end{filecontents*}
\csvreader{ddm11inf.csv}{}
				{\\\hline \csvcoli & \csvcolii & \csvcoliii & \csvcoliv & \csvcolv& \csvcolvi
					& \csvcolvii}
			\end{tabular}
			\begin{tabular}{c|c|c|c|c|c|c|c}
				\bfseries case 4 &\bfseries k
				& \bfseries case 5 &\bfseries k
				& \bfseries case 6 &\bfseries k
				& \bfseries case 7 &\bfseries k
				\csvreader{ddm11inf.csv}{}
				{\\\hline  \csvcolviii& \csvcolix & \csvcolx& \csvcolxi& \csvcolxii& \csvcolxiii
					& \csvcolxiv& \csvcolxv}
			\end{tabular}
			\caption{The \(L^\infty\) errors for DDM1}\label{afd2}
		\end{table}
		
		\begin{table}[H]
			\center
			
			\begin{tabular}{c|c|c|c|c|c|c}
				\bfseries \(\epsilon\) & \bfseries case 1 &\bfseries k  & \bfseries case 2 &\bfseries k
				& \bfseries case 3 &\bfseries k
				\begin{filecontents*}{ddm22.csv}
eps,case 1,k,case 2,k,case 3,k,case 4,k,case 5,k,case 6,k,case 7,k
2.00E-01,3.42E-01,0.00 ,5.80E-02,0.00 ,4.44E-02,0.00 ,9.27E-02,0.00 ,3.74E-01,0.00 ,1.99E-01,0.00 ,5.26E-02,0.00 
1.00E-01,1.70E-01,1.01 ,1.37E-02,2.09 ,2.28E-02,0.96 ,4.54E-02,1.03 ,1.79E-01,1.06 ,1.01E-01,0.97 ,2.60E-02,1.02 
5.00E-02,8.45E-02,1.01 ,3.31E-03,2.04 ,1.16E-02,0.98 ,2.25E-02,1.02 ,8.77E-02,1.03 ,5.12E-02,0.99 ,1.29E-02,1.01 
2.50E-02,4.17E-02,1.02 ,8.19E-04,2.02 ,5.84E-03,0.99 ,1.12E-02,1.00 ,4.30E-02,1.03 ,2.55E-02,1.01 ,6.44E-03,1.01 
1.25E-02,2.07E-02,1.01 ,2.03E-04,2.01 ,2.93E-03,0.99 ,5.59E-03,1.00 ,2.13E-02,1.01 ,1.27E-02,1.00 ,3.21E-03,1.00 
6.25E-03,1.03E-02,1.00 ,5.07E-05,2.00 ,1.47E-03,1.00 ,2.79E-03,1.00 ,1.06E-02,1.01 ,6.36E-03,1.00 ,1.60E-03,1.00 
3.13E-03,5.16E-03,1.00 ,1.27E-05,2.00 ,7.35E-04,1.00 ,1.39E-03,1.00 ,5.29E-03,1.00 ,3.18E-03,1.00 ,8.02E-04,1.00 
1.56E-03,2.58E-03,1.00 ,3.16E-06,2.00 ,3.68E-04,1.00 ,6.97E-04,1.00 ,2.64E-03,1.00 ,1.59E-03,1.00 ,4.01E-04,1.00 
\end{filecontents*}
\csvreader{ddm22.csv}{}
				{\\\hline \csvcoli & \csvcolii & \csvcoliii & \csvcoliv & \csvcolv& \csvcolvi
					& \csvcolvii}
			\end{tabular}
			\begin{tabular}{c|c|c|c|c|c|c|c}
				\bfseries case 4 &\bfseries k
				& \bfseries case 5 &\bfseries k
				& \bfseries case 6 &\bfseries k
				& \bfseries case 7 &\bfseries k
				\csvreader{ddm22.csv}{}
				{\\\hline  \csvcolviii& \csvcolix & \csvcolx& \csvcolxi& \csvcolxii& \csvcolxiii
					& \csvcolxiv& \csvcolxv}
			\end{tabular}
			\caption{The \(L^2\) errors for DDM2}\label{afd3}
		\end{table}
		
		\begin{table}[H]
			\center
			
			\begin{tabular}{c|c|c|c|c|c|c}
				\bfseries \(\epsilon\) & \bfseries case 1 &\bfseries k  & \bfseries case 2 &\bfseries k
				& \bfseries case 3 &\bfseries k
				\begin{filecontents*}{ddm22inf.csv}
eps,case 1,k,case 2,k,case 3,k,case 4,k,case 5,k,case 6,k,case 7,k
2.00E-01,1.65E-01,0.00 ,3.82E-02,0.00 ,3.53E-02,0.00 ,7.84E-02,0.00 ,1.94E-01,0.00 ,1.06E-01,0.00 ,4.34E-02,0.00 
1.00E-01,7.95E-02,1.05 ,9.03E-03,2.08 ,1.81E-02,0.96 ,3.83E-02,1.03 ,9.08E-02,1.10 ,5.26E-02,1.01 ,2.13E-02,1.03 
5.00E-02,3.88E-02,1.04 ,2.19E-03,2.04 ,9.11E-03,0.99 ,1.88E-02,1.03 ,4.36E-02,1.06 ,2.61E-02,1.01 ,1.05E-02,1.02 
2.50E-02,1.91E-02,1.02 ,5.42E-04,2.02 ,4.61E-03,0.98 ,9.39E-03,1.00 ,2.13E-02,1.04 ,1.30E-02,1.01 ,5.22E-03,1.01 
1.25E-02,9.51E-03,1.00 ,1.35E-04,2.01 ,2.33E-03,0.98 ,4.71E-03,0.99 ,1.06E-02,1.01 ,6.50E-03,1.00 ,2.62E-03,1.00 
6.25E-03,4.82E-03,0.98 ,3.38E-05,2.00 ,1.19E-03,0.97 ,2.40E-03,0.97 ,5.34E-03,0.98 ,3.31E-03,0.98 ,1.33E-03,0.98 
3.13E-03,2.40E-03,1.00 ,8.42E-06,2.00 ,5.93E-04,1.00 ,1.20E-03,1.00 ,2.66E-03,1.01 ,1.65E-03,1.00 ,6.64E-04,1.00 
1.56E-03,1.20E-03,1.01 ,2.10E-06,2.00 ,2.96E-04,1.00 ,5.96E-04,1.01 ,1.32E-03,1.01 ,8.22E-04,1.01 ,3.31E-04,1.01 
\end{filecontents*}
\csvreader{ddm22inf.csv}{}
				{\\\hline \csvcoli & \csvcolii & \csvcoliii & \csvcoliv & \csvcolv& \csvcolvi
					& \csvcolvii}
			\end{tabular}
			\begin{tabular}{c|c|c|c|c|c|c|c}
				\bfseries case 4 &\bfseries k
				& \bfseries case 5 &\bfseries k
				& \bfseries case 6 &\bfseries k
				& \bfseries case 7 &\bfseries k
				\csvreader{ddm22inf.csv}{}
				{\\\hline  \csvcolviii& \csvcolix & \csvcolx& \csvcolxi& \csvcolxii& \csvcolxiii
					& \csvcolxiv& \csvcolxv}
			\end{tabular}
			\caption{The \(L^\infty\) errors for DDM2}\label{afd4}
		\end{table}
		
		\begin{table}[H]
			\center
			
			\begin{tabular}{c|c|c|c|c|c|c}
				\bfseries \(\epsilon\) & \bfseries case 1 &\bfseries k  & \bfseries case 2 &\bfseries k
				& \bfseries case 3 &\bfseries k
				\begin{filecontents*}{ddm33.csv}
h,case 1,k,case 2,k,case 3,k,case 4,k,case 5,k,case 6,k,case 7,k
2.00E-01,3.82E-02,0.00 ,1.37E-02,0.00 ,6.82E-03,0.00 ,1.10E-02,0.00 ,3.57E-02,0.00 ,2.61E-02,0.00 ,6.36E-03,0.00 
1.00E-01,9.45E-02,-1.31 ,7.22E-03,0.92 ,1.44E-02,-1.08 ,2.59E-02,-1.23 ,9.39E-02,-1.39 ,6.02E-02,-1.21 ,1.50E-02,-1.23 
5.00E-02,8.03E-02,0.24 ,3.12E-03,1.21 ,1.17E-02,0.30 ,2.16E-02,0.26 ,8.08E-02,0.22 ,5.02E-02,0.26 ,1.25E-02,0.26 
2.50E-02,5.39E-02,0.58 ,1.19E-03,1.39 ,7.79E-03,0.59 ,1.46E-02,0.57 ,5.46E-02,0.56 ,3.35E-02,0.58 ,8.39E-03,0.57 
1.25E-02,3.33E-02,0.69 ,4.18E-04,1.51 ,4.80E-03,0.70 ,9.01E-03,0.69 ,3.39E-02,0.69 ,2.06E-02,0.70 ,5.19E-03,0.69 
6.25E-03,1.97E-02,0.76 ,1.40E-04,1.58 ,2.83E-03,0.76 ,5.34E-03,0.76 ,2.01E-02,0.75 ,1.22E-02,0.76 ,3.07E-03,0.76 
3.13E-03,1.14E-02,0.79 ,4.49E-05,1.64 ,1.63E-03,0.80 ,3.07E-03,0.80 ,1.16E-02,0.79 ,7.02E-03,0.80 ,1.77E-03,0.80 
1.56E-03,6.42E-03,0.82 ,1.43E-05,1.66 ,9.18E-04,0.83 ,1.74E-03,0.82 ,6.57E-03,0.82 ,3.96E-03,0.82 ,9.99E-04,0.82 
\end{filecontents*}
\csvreader{ddm33.csv}{}
				{\\\hline \csvcoli & \csvcolii & \csvcoliii & \csvcoliv & \csvcolv& \csvcolvi
					& \csvcolvii}
			\end{tabular}
			\begin{tabular}{c|c|c|c|c|c|c|c}
				\bfseries case 4 &\bfseries k
				& \bfseries case 5 &\bfseries k
				& \bfseries case 6 &\bfseries k
				& \bfseries case 7 &\bfseries k
				\csvreader{ddm33.csv}{}
				{\\\hline  \csvcolviii& \csvcolix & \csvcolx& \csvcolxi& \csvcolxii& \csvcolxiii
					& \csvcolxiv& \csvcolxv}
			\end{tabular}
			\caption{The \(L^2\) errors for DDM3}\label{afd5}
		\end{table}
		
		\begin{table}[H]
			\center
			
			\begin{tabular}{c|c|c|c|c|c|c}
				\bfseries \(\epsilon\) & \bfseries case 1 &\bfseries k  & \bfseries case 2 &\bfseries k
				& \bfseries case 3 &\bfseries k
				\begin{filecontents*}{ddm33inf.csv}
h,case 1,k,case 2,k,case 3,k,case 4,k,case 5,k,case 6,k,case 7,k
2.00E-01,4.73E-02,0.00 ,9.52E-03,0.00 ,1.06E-02,0.00 ,2.27E-02,0.00 ,5.44E-02,0.00 ,3.11E-02,0.00 ,1.26E-02,0.00 
1.00E-01,4.42E-02,0.10 ,4.92E-03,0.95 ,1.15E-02,-0.11 ,2.19E-02,0.06 ,4.75E-02,0.19 ,3.12E-02,-0.01 ,1.22E-02,0.04 
5.00E-02,3.66E-02,0.27 ,2.08E-03,1.24 ,9.17E-03,0.32 ,1.80E-02,0.28 ,3.99E-02,0.25 ,2.54E-02,0.30 ,1.00E-02,0.28 
2.50E-02,2.44E-02,0.59 ,7.82E-04,1.41 ,6.08E-03,0.59 ,1.21E-02,0.58 ,2.67E-02,0.58 ,1.68E-02,0.59 ,6.71E-03,0.58 
1.25E-02,1.50E-02,0.70 ,2.72E-04,1.52 ,3.73E-03,0.70 ,7.45E-03,0.70 ,1.65E-02,0.70 ,1.03E-02,0.70 ,4.14E-03,0.70 
6.25E-03,8.84E-03,0.76 ,9.04E-05,1.59 ,2.20E-03,0.76 ,4.40E-03,0.76 ,9.74E-03,0.76 ,6.09E-03,0.76 ,2.44E-03,0.76 
3.13E-03,5.09E-03,0.80 ,2.89E-05,1.64 ,1.26E-03,0.80 ,2.53E-03,0.80 ,5.61E-03,0.79 ,3.50E-03,0.80 ,1.41E-03,0.80 
1.56E-03,2.87E-03,0.82 ,9.14E-06,1.66 ,7.12E-04,0.83 ,1.43E-03,0.82 ,3.17E-03,0.82 ,1.98E-03,0.82 ,7.94E-04,0.82 
\end{filecontents*}
\csvreader{ddm33inf.csv}{}
				{\\\hline \csvcoli & \csvcolii & \csvcoliii & \csvcoliv & \csvcolv& \csvcolvi
					& \csvcolvii}
			\end{tabular}
			\begin{tabular}{c|c|c|c|c|c|c|c}
				\bfseries case 4 &\bfseries k
				& \bfseries case 5 &\bfseries k
				& \bfseries case 6 &\bfseries k
				& \bfseries case 7 &\bfseries k
				\csvreader{ddm33inf.csv}{}
				{\\\hline  \csvcolviii& \csvcolix & \csvcolx& \csvcolxi& \csvcolxii& \csvcolxiii
					& \csvcolxiv& \csvcolxv}
			\end{tabular}
			\caption{The \(L^\infty\) errors for DDM3}\label{afd6}
		\end{table}
		\section{Additional numerical results of mDDMs for 1D time-independent problems}\label{app_b}
		We have tested the following three mDDMs using \(h=\epsilon^{1.5}/4\):
		
		\textbf{mDDM1}: \(\nabla(\phi \nabla u)-\frac{1}{\epsilon^3}(1-\phi)(u-g-r\bm{n}\cdot \nabla u)=\phi f\),
		
		\textbf{mDDM2}: \( \phi\Delta u-\frac{1}{\epsilon^2}(1-\phi)(u-g-r\bm{n}\cdot\nabla u)=\phi f\),
		
		\textbf{mDDM3}: \(\nabla(\phi \nabla u)-\frac{1}{\epsilon^2}|\nabla \phi|(u-g- r\bm{n}\cdot \nabla u)=\phi f\),
		
		\noindent on the seven cases given in Appendix \ref{app_a}. The problem setup and the numerical discretization are analogous to those in Appendix \ref{app_a}. In case 3-7, mDDM1-3 perform analogously as those in case 1. However in case 2, in which \(A=0\), we observe similar behavior as those corresponding DDMs presented in Appendix \ref{app_a}. Our analysis in Sec. \ref{mddm2_ana} indicates \(\bar u_1^{(1)}(0) \approx-A/2.92\) (Eq. \eqref{mddm2_bar_u1}) and \(\hat u_1^{(1)}(0)\approx Ae^{1/36}/2.92\) (Eq. \eqref{mddm2_hat_u1}) for mDDM2, which both vanish when \(A=0\). Hence, mDDM2 can achieve 2nd order accuracy when \(A=0\). As for mDDM1 and mDDM3, the errors are again dominated by \({O}(\epsilon^2(\ln(\epsilon))^2)\), whose coefficient is not affected by our modification in the next order matching when \(A=0\).
		In Fig. \ref{d_eps_mddm1}, we plot \(D(\epsilon) = \lim_{z_2\rightarrow -\infty} (u_\epsilon-u)/\epsilon^2\) versus \(\ln(\epsilon)\) using the numerical solution of mDDM1 in case 2 and find that it is a quadratic function of \(\ln(\epsilon)\). 
		
		\begin{table}[H]
			\center
			
			\begin{tabular}{c|c|c|c|c|c|c}
				\bfseries \(\epsilon\) & \bfseries case 1 &\bfseries k  & \bfseries case 2 &\bfseries k
				& \bfseries case 3 &\bfseries k
				\begin{filecontents*}{mddm11.csv}
eps,case 1,k,case 2,k,case 3,k,case 4,k,case 5,k,case 6,k,case 7,k
2.00E-01,4.88E-01,0.00 ,3.11E-02,0.00 ,7.00E-02,0.00 ,1.25E-01,0.00 ,4.81E-01,0.00 ,3.07E-01,0.00 ,7.37E-02,0.00 
1.00E-01,1.50E-01,1.70 ,9.11E-03,1.77 ,2.19E-02,1.68 ,3.96E-02,1.66 ,1.49E-01,1.69 ,9.44E-02,1.70 ,2.31E-02,1.68 
5.00E-02,4.46E-02,1.75 ,2.64E-03,1.79 ,6.57E-03,1.73 ,1.20E-02,1.73 ,4.45E-02,1.74 ,2.81E-02,1.75 ,6.94E-03,1.73 
2.50E-02,1.26E-02,1.83 ,8.14E-04,1.70 ,1.87E-03,1.81 ,3.40E-03,1.82 ,1.25E-02,1.83 ,7.94E-03,1.82 ,1.97E-03,1.82 
1.25E-02,3.34E-03,1.91 ,2.69E-04,1.60 ,5.05E-04,1.89 ,9.09E-04,1.90 ,3.31E-03,1.92 ,2.12E-03,1.90 ,5.25E-04,1.90 
6.25E-03,8.61E-04,1.95 ,9.25E-05,1.54 ,1.33E-04,1.92 ,2.36E-04,1.95 ,8.44E-04,1.97 ,5.53E-04,1.94 ,1.36E-04,1.95 
\end{filecontents*}
\csvreader{mddm11.csv}{}
				{\\\hline \csvcoli & \csvcolii & \csvcoliii & \csvcoliv & \csvcolv& \csvcolvi
					& \csvcolvii}
			\end{tabular}
			\begin{tabular}{c|c|c|c|c|c|c|c}
				\bfseries case 4 &\bfseries k
				& \bfseries case 5 &\bfseries k
				& \bfseries case 6 &\bfseries k
				& \bfseries case 7 &\bfseries k
				\csvreader{mddm11.csv}{}
				{\\\hline  \csvcolviii& \csvcolix & \csvcolx& \csvcolxi& \csvcolxii& \csvcolxiii
					& \csvcolxiv& \csvcolxv}
			\end{tabular}
			\caption{The \(L^2\) errors for mDDM1}\label{af1}
		\end{table}
		
		\begin{table}[H]
			\center
			
			\begin{tabular}{c|c|c|c|c|c|c}
				\bfseries \(\epsilon\) & \bfseries case 1 &\bfseries k  & \bfseries case 2 &\bfseries k
				& \bfseries case 3 &\bfseries k
				\begin{filecontents*}{mddm11inf.csv}
eps,case 1,k,case 2,k,case 3,k,case 4,k,case 5,k,case 6,k,case 7,k
2.00E-01,2.10E-01,0.00 ,2.06E-02,0.00 ,5.51E-02,0.00 ,1.05E-01,0.00 ,2.26E-01,0.00 ,1.49E-01,0.00 ,5.85E-02,0.00 
1.00E-01,6.59E-02,1.67 ,5.93E-03,1.80 ,1.71E-02,1.69 ,3.29E-02,1.67 ,7.11E-02,1.67 ,4.64E-02,1.68 ,1.83E-02,1.67 
5.00E-02,1.98E-02,1.74 ,1.70E-03,1.80 ,5.12E-03,1.74 ,9.90E-03,1.73 ,2.14E-02,1.73 ,1.39E-02,1.74 ,5.51E-03,1.73 
2.50E-02,5.58E-03,1.83 ,5.22E-04,1.70 ,1.45E-03,1.82 ,2.81E-03,1.82 ,6.01E-03,1.83 ,3.94E-03,1.82 ,1.56E-03,1.82 
1.25E-02,1.55E-03,1.85 ,1.72E-04,1.60 ,3.92E-04,1.89 ,7.71E-04,1.86 ,1.70E-03,1.82 ,1.07E-03,1.88 ,4.28E-04,1.87 
6.25E-03,4.91E-04,1.65 ,5.91E-05,1.54 ,1.22E-04,1.68 ,2.45E-04,1.65 ,5.42E-04,1.65 ,3.39E-04,1.66 ,1.36E-04,1.65 
\end{filecontents*}
\csvreader{mddm11inf.csv}{}
				{\\\hline \csvcoli & \csvcolii & \csvcoliii & \csvcoliv & \csvcolv& \csvcolvi
					& \csvcolvii}
			\end{tabular}
			\begin{tabular}{c|c|c|c|c|c|c|c}
				\bfseries case 4 &\bfseries k
				& \bfseries case 5 &\bfseries k
				& \bfseries case 6 &\bfseries k
				& \bfseries case 7 &\bfseries k
				\csvreader{mddm11inf.csv}{}
				{\\\hline  \csvcolviii& \csvcolix & \csvcolx& \csvcolxi& \csvcolxii& \csvcolxiii
					& \csvcolxiv& \csvcolxv}
			\end{tabular}
			\caption{The \(L^\infty\) errors for mDDM1}\label{af2}
		\end{table}
		
		\begin{table}[H]
			\center
			
			\begin{tabular}{c|c|c|c|c|c|c}
				\bfseries \(\epsilon\) & \bfseries case 1 &\bfseries k  & \bfseries case 2 &\bfseries k
				& \bfseries case 3 &\bfseries k
				\begin{filecontents*}{mddm22.csv}
eps,case 1,k,case 2,k,case 3,k,case 4,k,case 5,k,case 6,k,case 7,k
2.00E-01,2.85E-01,0.00 ,3.73E-02,0.00 ,3.75E-02,0.00 ,7.65E-02,0.00 ,3.07E-01,0.00 ,1.68E-01,0.00 ,4.36E-02,0.00 
1.00E-01,1.41E-01,1.01 ,8.77E-03,2.09 ,1.92E-02,0.97 ,3.77E-02,1.02 ,1.48E-01,1.05 ,8.48E-02,0.98 ,2.16E-02,1.01 
5.00E-02,6.98E-02,1.01 ,2.13E-03,2.04 ,9.70E-03,0.98 ,1.87E-02,1.01 ,7.23E-02,1.03 ,4.25E-02,1.00 ,1.08E-02,1.01 
2.50E-02,3.46E-02,1.01 ,5.24E-04,2.02 ,4.88E-03,0.99 ,9.33E-03,1.00 ,3.57E-02,1.02 ,2.12E-02,1.00 ,5.37E-03,1.00 
1.25E-02,1.73E-02,1.00 ,1.30E-04,2.01 ,2.45E-03,1.00 ,4.66E-03,1.00 ,1.77E-02,1.01 ,1.06E-02,1.00 ,2.68E-03,1.00 
6.25E-03,8.62E-03,1.00 ,3.23E-05,2.01 ,1.23E-03,1.00 ,2.33E-03,1.00 ,8.84E-03,1.00 ,5.31E-03,1.00 ,1.34E-03,1.00 
\end{filecontents*}
\csvreader{mddm22.csv}{}
				{\\\hline \csvcoli & \csvcolii & \csvcoliii & \csvcoliv & \csvcolv& \csvcolvi
					& \csvcolvii}
			\end{tabular}
			\begin{tabular}{c|c|c|c|c|c|c|c}
				\bfseries case 4 &\bfseries k
				& \bfseries case 5 &\bfseries k
				& \bfseries case 6 &\bfseries k
				& \bfseries case 7 &\bfseries k
				\csvreader{mddm22.csv}{}
				{\\\hline  \csvcolviii& \csvcolix & \csvcolx& \csvcolxi& \csvcolxii& \csvcolxiii
					& \csvcolxiv& \csvcolxv}
			\end{tabular}
			\caption{The \(L^2\) errors for mDDM2}\label{af3}
		\end{table}
		
		\begin{table}[H]
			\center
			
			\begin{tabular}{c|c|c|c|c|c|c}
				\bfseries \(\epsilon\) & \bfseries case 1 &\bfseries k  & \bfseries case 2 &\bfseries k
				& \bfseries case 3 &\bfseries k
				\begin{filecontents*}{mddm22inf.csv}
eps,case 1,k,case 2,k,case 3,k,case 4,k,case 5,k,case 6,k,case 7,k
2.00E-01,1.33E-01,0.00 ,2.50E-02,0.00 ,2.96E-02,0.00 ,6.43E-02,0.00 ,1.55E-01,0.00 ,8.67E-02,0.00 ,3.56E-02,0.00 
1.00E-01,6.48E-02,1.04 ,5.86E-03,2.09 ,1.52E-02,0.96 ,3.17E-02,1.02 ,7.35E-02,1.07 ,4.34E-02,1.00 ,1.76E-02,1.02 
5.00E-02,3.19E-02,1.02 ,1.42E-03,2.05 ,7.67E-03,0.98 ,1.57E-02,1.01 ,3.57E-02,1.04 ,2.16E-02,1.00 ,8.73E-03,1.01 
2.50E-02,1.59E-02,1.01 ,3.52E-04,2.01 ,3.88E-03,0.98 ,7.89E-03,1.00 ,1.77E-02,1.02 ,1.09E-02,1.00 ,4.38E-03,1.00 
1.25E-02,7.93E-03,1.00 ,8.73E-05,2.01 ,1.95E-03,1.00 ,3.94E-03,1.00 ,8.78E-03,1.01 ,5.43E-03,1.00 ,2.19E-03,1.00 
6.25E-03,3.96E-03,1.00 ,2.18E-05,2.00 ,9.77E-04,1.00 ,1.97E-03,1.00 ,4.39E-03,1.00 ,2.72E-03,1.00 ,1.09E-03,1.00 
\end{filecontents*}
\csvreader{mddm22inf.csv}{}
				{\\\hline \csvcoli & \csvcolii & \csvcoliii & \csvcoliv & \csvcolv& \csvcolvi
					& \csvcolvii}
			\end{tabular}
			\begin{tabular}{c|c|c|c|c|c|c|c}
				\bfseries case 4 &\bfseries k
				& \bfseries case 5 &\bfseries k
				& \bfseries case 6 &\bfseries k
				& \bfseries case 7 &\bfseries k
				\csvreader{mddm22inf.csv}{}
				{\\\hline  \csvcolviii& \csvcolix & \csvcolx& \csvcolxi& \csvcolxii& \csvcolxiii
					& \csvcolxiv& \csvcolxv}
			\end{tabular}
			\caption{The \(L^\infty\) errors for mDDM2}\label{af4}
		\end{table}
		
		\begin{table}[H]
			\center
			
			\begin{tabular}{c|c|c|c|c|c|c}
				\bfseries \(\epsilon\) & \bfseries case 1 &\bfseries k  & \bfseries case 2 &\bfseries k
				& \bfseries case 3 &\bfseries k
				\begin{filecontents*}{mddm33.csv}
h,case 1,k,case 2,k,case 3,k,case 4,k,case 5,k,case 6,k,case 7,k
2.00E-01,1.55E-01,0.00 ,3.89E-02,0.00 ,2.80E-02,0.00 ,4.42E-02,0.00 ,1.44E-01,0.00 ,1.07E-01,0.00 ,2.57E-02,0.00 
1.00E-01,3.77E-02,2.04 ,1.38E-02,1.50 ,7.15E-03,1.97 ,1.08E-02,2.03 ,3.31E-02,2.12 ,2.69E-02,1.99 ,6.31E-03,2.02 
5.00E-02,8.54E-03,2.14 ,4.87E-03,1.50 ,1.80E-03,1.99 ,2.51E-03,2.10 ,6.78E-03,2.29 ,6.50E-03,2.05 ,1.48E-03,2.09 
2.50E-02,1.89E-03,2.17 ,1.72E-03,1.50 ,4.71E-04,1.94 ,5.86E-04,2.10 ,1.25E-03,2.44 ,1.59E-03,2.03 ,3.48E-04,2.09 
1.25E-02,4.16E-04,2.19 ,5.88E-04,1.55 ,1.27E-04,1.89 ,1.38E-04,2.09 ,1.91E-04,2.71 ,4.01E-04,1.99 ,8.28E-05,2.07 
6.25E-03,8.74E-05,2.25 ,1.94E-04,1.60 ,3.46E-05,1.88 ,3.20E-05,2.11 ,1.39E-05,3.78 ,1.01E-04,1.98 ,1.96E-05,2.08 
\end{filecontents*}
\csvreader{mddm33.csv}{}
				{\\\hline \csvcoli & \csvcolii & \csvcoliii & \csvcoliv & \csvcolv& \csvcolvi
					& \csvcolvii}
			\end{tabular}
			\begin{tabular}{c|c|c|c|c|c|c|c}
				\bfseries case 4 &\bfseries k
				& \bfseries case 5 &\bfseries k
				& \bfseries case 6 &\bfseries k
				& \bfseries case 7 &\bfseries k
				\csvreader{mddm33.csv}{}
				{\\\hline  \csvcolviii& \csvcolix & \csvcolx& \csvcolxi& \csvcolxii& \csvcolxiii
					& \csvcolxiv& \csvcolxv}
			\end{tabular}
			\caption{The \(L^2\) errors for mDDM3}\label{af5}
		\end{table}
		
		\begin{table}[H]
			\center
			
			\begin{tabular}{c|c|c|c|c|c|c}
				\bfseries \(\epsilon\) & \bfseries case 1 &\bfseries k  & \bfseries case 2 &\bfseries k
				& \bfseries case 3 &\bfseries k
				\begin{filecontents*}{mddm33inf.csv}
h,case 1,k,case 2,k,case 3,k,case 4,k,case 5,k,case 6,k,case 7,k
2.00E-01,8.23E-02,0.00 ,2.58E-02,0.00 ,2.18E-02,0.00 ,4.10E-02,0.00 ,8.94E-02,0.00 ,5.78E-02,0.00 ,2.29E-02,0.00 
1.00E-01,2.23E-02,1.88 ,8.98E-03,1.52 ,5.82E-03,1.91 ,1.12E-02,1.88 ,2.40E-02,1.89 ,1.58E-02,1.88 ,6.22E-03,1.88 
5.00E-02,5.96E-03,1.91 ,3.14E-03,1.51 ,1.56E-03,1.90 ,2.99E-03,1.90 ,6.39E-03,1.91 ,4.22E-03,1.90 ,1.67E-03,1.90 
2.50E-02,2.03E-03,1.55 ,1.10E-03,1.51 ,5.17E-04,1.59 ,1.02E-03,1.56 ,2.21E-03,1.53 ,1.41E-03,1.58 ,5.64E-04,1.56 
1.25E-02,6.15E-04,1.72 ,3.76E-04,1.55 ,1.56E-04,1.73 ,3.07E-04,1.72 ,6.72E-04,1.72 ,4.27E-04,1.73 ,1.71E-04,1.72 
6.25E-03,2.29E-04,1.42 ,1.24E-04,1.60 ,5.73E-05,1.44 ,1.14E-04,1.43 ,2.52E-04,1.41 ,1.58E-04,1.43 ,6.35E-05,1.43 
\end{filecontents*}
\csvreader{mddm33inf.csv}{}
				{\\\hline \csvcoli & \csvcolii & \csvcoliii & \csvcoliv & \csvcolv& \csvcolvi
					& \csvcolvii}
			\end{tabular}
			\begin{tabular}{c|c|c|c|c|c|c|c}
				\bfseries case 4 &\bfseries k
				& \bfseries case 5 &\bfseries k
				& \bfseries case 6 &\bfseries k
				& \bfseries case 7 &\bfseries k
				\csvreader{mddm33inf.csv}{}
				{\\\hline  \csvcolviii& \csvcolix & \csvcolx& \csvcolxi& \csvcolxii& \csvcolxiii
					& \csvcolxiv& \csvcolxv}
			\end{tabular}
			\caption{The \(L^\infty\) errors for mDDM3}\label{af6}
		\end{table}
		\section{Validation of the asymptotic analysis}\label{app_c}
		Here we present validations of our asymptotic analysis using the numerical results for the seven cases.
		In Tab. \ref{asy_th2}, we present \(\bar u_1^{(1)}\)(0) from DDM2 obtained from both our asymptotic analysis theory (Eq. \eqref{ddm2_out_sol2}) and the numerical results. In Tab. \ref{ubar1_mddm13}, we compare the slope of \(C(\epsilon)\) computed numerically (through an analogous linear fit as in Sec. \ref{ddm13}) with that derived from our asymptotic analysis theory (\(-A/6\)) for DDM1 and DDM3. In Tab. \ref{uhat1_mddm13}, we show \(\hat u_3^{(1.5)}(0)\) from our asymptotic theory (\(-A/\sqrt{2\pi}\) for mDDM1 and \(-A/\sqrt{6\pi}\) for mDDM3) together with those calculated from the numerical results for mDDM1 and mDDM3. in Tab. \ref{ubar1__mddm2}, \(\bar u_1^{(1)}(0)\) and \(\hat u^{(1)}(0)\) from mDDM2 obtained from both our asymptotic analysis theory (Eqs. \eqref{mddm2_bar_u1} and \eqref{mddm2_hat_u1}) and the numerical results are presented. Clearly our theory is consistent with the numerical results.
		
		\begin{table}[H]
			\center
			\begin{tabular}{c|c|c}
				\bfseries case & \bfseries Theory &\bfseries Numerics  
				\begin{filecontents*}{ubar1_mddm2.csv}
rate,theory,numeric
case 1,-0.450 ,-0.450 
case 2,0.000 ,0.000 
case 3,0.180 ,0.183 
case 4,0.363 ,0.368 
case 5,-4.020 ,-4.072 
case 6,-0.403 ,-0.408 
case 7,-0.301 ,-0.305 
\end{filecontents*}
\csvreader{ubar1_mddm2.csv}{}
				{\\\hline \csvcoli & \csvcolii & \csvcoliii }
			\end{tabular}
			\caption{Comparisons between \(\bar u_1^{(1)}\) from theory and that from numerical results.}\label{asy_th2}
		\end{table}
		\begin{table}[H]
			\center
			\begin{tabular}{c|c|c|c}
				\bfseries case & \bfseries Theory &\bfseries Numerics (DDM1)& \bfseries Numerics (DDM3)  
				\begin{filecontents*}{ubar1_mddm13.csv}
rate,predict,numeric,numeric
case 1,-0.185 ,-0.186 ,-0.188 
case 2,0.000 ,0.000 ,0.000 
case 3,0.074 ,0.075 ,0.075 
case 4,0.149 ,0.150 ,0.152 
case 5,-1.655 ,-1.667 ,-1.686 
case 6,-0.166 ,-0.167 ,-0.168 
case 7,-0.124 ,-0.125 ,-0.126 
\end{filecontents*}
\csvreader{ubar1_mddm13.csv}{}
				{\\\hline \csvcoli & \csvcolii & \csvcoliii &\csvcoliv }
			\end{tabular}
			\caption{Comparisons between the slope of \(C(\epsilon)\) from theory and that from numerical results.}\label{ubar1_mddm13}
		\end{table}
		\begin{table}[H]
			\center
			\begin{tabular}{c||c|c||c|c}
				&\multicolumn{2}{c||}{mDDM1}&\multicolumn{2}{c}{mDDM3}\\ \hline
				case&Theory&Numerics&Theory&Numerics
				\begin{filecontents*}{uhat1_mddm13.csv}
rate,predict,numeric,predict ,numeric
case 1,-0.443 ,-0.442 ,-0.256 ,-0.257 
case 2,0.000 ,0.001 ,0.000 ,-0.001 
case 3,0.178 ,0.177 ,0.103 ,0.103 
case 4,0.356 ,0.357 ,0.206 ,0.207 
case 5,-3.961 ,-3.952 ,-2.287 ,-2.299 
case 6,-0.397 ,-0.396 ,-0.229 ,-0.230 
case 7,-0.297 ,-0.296 ,-0.171 ,-0.172 
\end{filecontents*}
\csvreader[]{uhat1_mddm13.csv}{}
				{\\\hline \csvcoli & \csvcolii & \csvcoliii & \csvcoliv & \csvcolv}
			\end{tabular}
			\caption{Comparisons between \(\hat{u}_3^{(1.5)}(0)\) from asymptotic theory and from numerical simulations of mDDM1 and mDDM3.}\label{uhat1_mddm13}
		\end{table}
		\begin{table}[H]
			\center
			\begin{tabular}{c||c|c||c|c}
				&\multicolumn{2}{c||}{\(\bar u^{(1)}(0)\)}&\multicolumn{2}{c}{\(\hat u^{(1)}(0)\)}\\ \hline
				case&Predictions&Numerics&Predictions&Numerics
				\begin{filecontents*}{ubar1__mddm2.csv}
rate,predict,numeric,predict ,numeric
case 1,-0.381 ,-0.380 ,-0.391 ,-0.391 
case 2,0.000 ,0.000 ,0.000 ,0.000 
case 3,0.152 ,0.152 ,0.157 ,0.157 
case 4,0.307 ,0.307 ,0.316 ,0.315 
case 5,-3.401 ,-3.399 ,-3.496 ,-3.495 
case 6,-0.341 ,-0.340 ,-0.350 ,-0.350 
case 7,-0.255 ,-0.254 ,-0.262 ,-0.262 
\end{filecontents*}
\csvreader[]{ubar1__mddm2.csv}{}
				{\\\hline \csvcoli & \csvcolii & \csvcoliii & \csvcoliv & \csvcolv}
			\end{tabular}
			\caption{Comparisons between the asymptotic theory and the numerical results for \(\bar u^{(1)}(0)\) and \(\hat u^{(1)}(0)\) from mDDM2.}\label{ubar1__mddm2}
		\end{table}

		\section{Derivation of the solution to Eq. \eqref{mddm2_ode}}\label{app_d}
		Recall the homogeneous ordinary differential equation from Eq. \eqref{mddm2_ode} in the main text,
		\begin{equation}
		y''-e^{6x}(y-xy')=0.\notag
		\end{equation}
		Clearly, \(y_1=x\) is one of the linearly independent solutions.
		We derive the other solution through a reduction of order. We assume \(y_2(x)=v(x)y_1(x)=xv(x)\) and plug into the equation to get:
		\begin{equation}
		xv''+2v'+x^2e^{6x}v'=0, \label{mddm2_ode1}
		\end{equation}
		which gives \(v(x)=\int\frac{e^{e^{6x}(1-6x)/36}}{x^2}dx\) for \(x\neq 0\). Hence,
		\begin{equation}
		y_2=x\int\frac{e^{e^{6x}(1-6x)/36}}{x^2}dx = -e^{e^{6x}(1-6x)/36}-x\int h(x)dx \text{ for } x\neq 0,
		\end{equation}
		where \(h(x)=e^{e^{6x}(1-6x)/36+6x}\). It is not hard to verify that \(y_2\) is a solution to the equation for all \(x\) including 0. Although the anti-derivative of \(h(x)\) is not an elementary function, \(h(x)\in L^1(-\infty,+\infty)\) and \(\int_{-\infty}^{+\infty} h(x)dx\approx 2.92\). Let \(H'(x)=h(x)\), then \(y_2\) can be written as
		\begin{equation}
		y_2=-e^{e^{6x}(1-6x)/36}-x\int_0^xh(t)dt-H(0)x.
		\end{equation}
		Note that \(H(0)x\) is linearly dependent with respect to \(y_1\), thus \(y_2\) can be simplified as 
		\begin{equation}
		y_2=-e^{e^{6x}(1-6x)/36}-x\int_0^xh(t)dt.
		\end{equation} 
		Hence, the general solution to Eq. \eqref{mddm2_ode} is
		\begin{align}
		y=C_1x+C_2(-e^{e^{6x}(1-6x)/36}-x\int_0^xh(t)dt).
		\end{align}

		{\color{black}

		\section{Analysis of the truncation and analytic errors}\label{res_trunc_analytic}

In Sec. \ref{1D time dependent}, we showed that the errors using mDDMt3 to solve time-dependent problems could actually be smaller using $h=\epsilon/c$ when $c$ is small (e.g., $c=4$) than when $c$ is large (e.g., $c=128$). Recall Fig. \ref{comp_c_mddmt3_00625_sta} and Tables \ref{ddmts1} and \ref{ddmts2}. Here, we present an analysis of the numerical results to explain this surprising behavior. Because the errors are dominated by the space discretization, we consider, for simplicity, the analogous time-independent problem in case 1 using mDDM3 (see App. \ref{app_a} for a definition of case 1; this is also the problem considered in Sec. \ref{num_1d_ddm}). 

In Fig. \ref{soln_c_mddm3} below, we show that for the time-independent problem in case 1, the solution using mDDM3 with $c=4$ is also closer to the exact solution than those with larger values of $c$. Indeed, the behavior of the solution and the errors are very similar to the time-dependent case presented in the main text in Sec. \ref{1D time dependent}. The actual relative errors $||u_{h,\epsilon}-u||_\infty/||u||_\infty$ are given in Table \ref{inf_err}. To understand this behavior write the total error as
$$
||u_{h,\epsilon}-u||_\infty=||\left(u_{h,\epsilon}-u_{\epsilon}\right)+\left(u_{\epsilon}-u\right)||_\infty,
$$
where $u_\epsilon$ is the solution of the continuous modified equation. The first term on the right hand side is the truncation error and the second term on the right is the analytic error. 
The truncation error can be approximated using the consecutive errors $E^\infty_{h,h/2,\epsilon}=u_{h,\epsilon}-u_{h/2,\epsilon}$ evaluated at the point $x_{h,h/2}^*=argmax|u_{h,\epsilon}-u_{h/2,\epsilon}|$. The consecutive errors with $h=\epsilon/c$ are shown in Table \ref{h_err} below. Note that there is no norm in the evaluation of $E^\infty_{h,h/2,\epsilon}$ and that all the consecutive errors are positive. The truncation error can then be approximated by $u_{h,\epsilon}-u_\epsilon\approx\sum_{j=0}^\infty E^\infty_{h/2^j,h/2^{j+1},\epsilon}$, which is thus also positive. Next, as shown in Sec. \ref{aa_mddm1_3}, the analytic error near the boundary of the physical domain $D$ can be approximated to leading order as $u_\epsilon-u\sim - \frac{A\epsilon^{1.5}}{\sqrt{6\pi}}<0$. where $A$ is the normal derivative of the analytic solution at the boundary ($A=1.111$ for the problem in case 1). The leading-order analytic relative errors $\frac{u_\epsilon-u}{||u||_\infty}$ are shown in the last row of Table \ref{h_err} and are of a similar magnitude as the truncation errors (at least for $\epsilon$ smaller than 0.1 and $c$ is small) but have the opposite sign. For moderate sizes of $h$ and $\epsilon$, the truncation and analytic errors partially cancel one another. As $h$ decreases, however, the analytic error tends to dominate and the convergence rate decreases from close to 2 for larger $h$ (smaller $c$) to 1.5 for small $h$ (large $c$). A careful analysis shows that when $\epsilon$ is small (e.g., 0.05, 0.025, 0.0125), the combination of truncation and leading-order asymptotic analytic errors provides a good estimate of the total error. For larger $\epsilon$,  the estimate is not as good because the leading-order asymptotic error does not provide as good an estimate for the overall analytic error.

	\begin{figure}[H]
	\includegraphics[width=.75\linewidth]{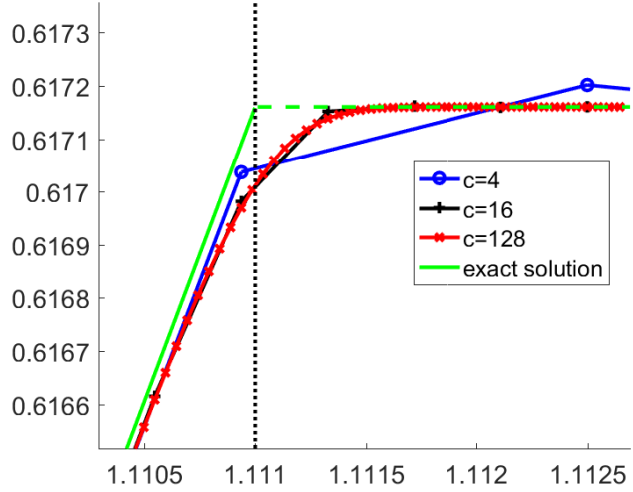}
	\caption{Numerical solutions of mDDM3 near the boundary x=1.111 for the time-independent problem in case 1 with different c such that $c=\epsilon/h$. The green dashed line shows the extension of the exact solution outside the domain (e.g., constant in the normal direction). The vertical dotted line represents the boundary of the physical domain (x=1.111). Observe that the mDDM3 solution with c=4 is closer to the exact solution than those with larger c}\label{soln_c_mddm3}
\end{figure}

\begin{table}[H]
	\center
	\begin{tabular}{c|c|c|c|c|c|c}
		\bfseries $\epsilon$& c=4 &\bfseries k  & c=16 &\bfseries k& c=128 &\bfseries k
		\begin{filecontents*}{mddm3_inf1.csv}
eps,rate=4,k,rate =16,k,rate=128,k
0.2,8.19E-02,0.00 ,8.15E-02,0.00 ,8.29E-02,0.00 
0.1,2.05E-02,2.00 ,2.22E-02,1.88 ,2.38E-02,1.80 
0.05,4.26E-03,2.26 ,6.34E-03,1.81 ,7.12E-03,1.74 
0.025,9.22E-04,2.21 ,2.17E-03,1.55 ,2.20E-03,1.69 
0.0125,2.93E-04,1.66 ,6.55E-04,1.73 ,6.82E-04,1.69 
0.00625,8.76E-05,1.74 ,1.77E-04,1.88 ,2.30E-04,1.57 
\end{filecontents*}
\csvreader{mddm3_inf1.csv}{}
		{\\\hline \csvcoli & \csvcolii & \csvcoliii & \csvcoliv & \csvcolv& \csvcolvi & \csvcolvii}
	\end{tabular}
	\caption{The total $L^\infty$ relative error $\frac{||u_{h,\epsilon}-u||_\infty}{||u||_\infty}$ with $h=\epsilon/c$, as labeled, for the time-independent problem in case 1 using mDDM3.}
	\label{inf_err}
\end{table}

		\begin{table}[H]
 	\footnotesize 
	\center
	\begin{tabular}{c||c|c||c|c||c|c||c|c||c|c}
	    &\multicolumn{2}{c}{\(\epsilon=0.2\)}&\multicolumn{2}{c}{\(\epsilon=0.1\)}&\multicolumn{2}{c}{\(\epsilon=0.05\)}&\multicolumn{2}{c}{\(\epsilon=0.025\)}&\multicolumn{2}{c}{\(\epsilon=0.0125\)}\\ \hline
	 c&\(E^{\infty,r}_{h,h/2,\epsilon}\)&\(\bm{k}\)&\(E^{\infty,r}_{h,h/2,\epsilon}\)&\(\bm{k}\)&\(E^{\infty,r}_{h,h/2,\epsilon}\)&\(\bm{k}\)&\(E^{\infty,r}_{h,h/2,\epsilon}\)&\(\bm{k}\)&\(E^{\infty,r}_{h,h/2,\epsilon}\)&\(\bm{k}\)
		\begin{filecontents*}{pos_or_neg.csv}
2,,,,,,,,,,
4,2.90E-03,,2.24E-03,,1.41E-03,,3.92E-04,,3.29E-04,
8,5.37E-04,2.43,4.14E-04,2.43,2.36E-04,2.58,1.54E-04,1.34,1.50E-04,1.13
16,1.15E-04,2.22,7.61E-05,2.44,6.09E-05,1.96,5.62E-05,1.46,4.13E-05,1.86
32,2.73E-05,2.07,1.87E-05,2.03,1.54E-05,1.98,1.60E-05,1.81,1.27E-05,1.7
64,6.73E-06,2.02,4.64E-06,2.01,4.23E-06,1.87,4.09E-06,1.97,3.29E-06,1.94
128,1.68E-06,2,1.16E-06,2.01,1.06E-06,2,1.03E-06,1.99,8.30E-07,1.99
256,4.24E-07,1.99,2.85E-07,2.02,2.63E-07,2.01,2.57E-07,2,2.11E-07,1.98
512,1.06E-07,1.99,6.71E-08,2.09,6.36E-08,2.05,6.46E-08,1.99,5.30E-08,1.99
\end{filecontents*}
\csvreader[]{pos_or_neg.csv}{}
		{\\\hline \csvcoli & \csvcolii & \csvcoliii & \csvcoliv & \csvcolv & \csvcolvi & \csvcolvii & \csvcolviii & \csvcolix & \csvcolx & \csvcolxi} \\ \hline  \\ \hline Analytic error  & -4.58E-02	& &	-1.62E-02& &		-5.7E-03& &		-2.02E-03	& &	-7.16E-04
		
	\end{tabular}
	\caption{The consecutive relative truncation errors $E^{\infty,r}_{h,h/2,\epsilon}/||u||_\infty$ with $c=\epsilon/h$ as labeled for the time-independent problem in case 1 using mDDM3. The results show that for fixed $\epsilon$, the consecutive errors are positive and that the scheme converges with 2nd order accuracy in $h$. The last row shows the leading-order asymptotic analytic relative error $\frac{u_\epsilon-u}{||u||_\infty}\sim-\frac{2\epsilon^{1.5}}{\sqrt{6\pi}}$.}\label{h_err}
\end{table}

		\section{Numerical results for mDDMt3 for 1D time-dependent problems on moving domains}\label{res_1d_mov}

		In Tables \ref{mDDMt3_mov_l2} and \ref{mDDMt3_mov_inf} below, we present the $L^2$ and $L^\infty$ errors for the problem considered in Sec. \ref{1D time dependent} and order of convergence $\mathbf{k}$ with $h=\epsilon/c$.

		\begin{table}[H]
			\center
			\begin{tabular}{c|c|c|c|c|c|c}
				\bfseries \(\epsilon\) & c = 4&\bfseries k  & c = 16 &\bfseries k & c = 128 &\bfseries k
				\begin{filecontents*}{mddmt3_mov_l2.csv}
eps,rate=4,k,rate =16,k,rate=128,k
0.2,3.07E-02,0.00 ,3.07E-02,0.00 ,3.07E-02,0.00 
0.1,7.92E-03,1.95 ,8.05E-03,1.93 ,8.06E-03,1.93 
0.05,1.90E-03,2.06 ,1.96E-03,2.04 ,1.96E-03,2.04 
0.025,4.74E-04,2.00 ,4.79E-04,2.03 ,4.79E-04,2.03 
0.0125,1.18E-04,2.01 ,1.21E-04,1.99 ,1.21E-04,1.99 
0.00625,2.97E-05,1.99 ,3.12E-05,1.96 ,3.13E-05,1.95 
\end{filecontents*}
\csvreader{mddmt3_mov_l2.csv}{}
				{\\\hline \csvcoli & \csvcolii & \csvcoliii & \csvcoliv & \csvcolv& \csvcolvi & \csvcolvii}
			\end{tabular}
			\caption{The \(L^{2}\) errors for simulating the time-dependent (diffusion) equation using mDDMt3 with \(h=\epsilon/c\) on a moving domain.}\label{mDDMt3_mov_l2}
		\end{table}
		\begin{table}[H]
			\center
			\begin{tabular}{c|c|c|c|c|c|c}
				\bfseries \(\epsilon\) & c = 4&\bfseries k  & c = 16 &\bfseries k & c = 128 &\bfseries k
				\begin{filecontents*}{mddmt3_mov_inf.csv}
eps,rate=4,k,rate =16,k,rate=128,k
0.2,4.38E-02,0.00 ,4.35E-02,0.00 ,4.51E-02,0.00 
0.1,1.16E-02,1.92 ,1.25E-02,1.80 ,1.35E-02,1.74 
0.05,2.57E-03,2.17 ,3.62E-03,1.79 ,4.05E-03,1.74 
0.025,6.36E-04,2.01 ,1.23E-03,1.56 ,1.25E-03,1.70 
0.0125,2.01E-04,1.66 ,3.69E-04,1.74 ,3.85E-04,1.70 
0.00625,5.01E-05,2.00 ,1.00E-04,1.88 ,1.29E-04,1.58 
\end{filecontents*}
\csvreader{mddmt3_mov_inf.csv}{}
				{\\\hline \csvcoli & \csvcolii & \csvcoliii & \csvcoliv & \csvcolv& \csvcolvi & \csvcolvii}
			\end{tabular}
			\caption{The \(L^{\infty}\) errors for simulating the time-dependent (diffusion) equation using mDDMt3 with \(h=\epsilon/c\) on a moving domain.}\label{mDDMt3_mov_inf}
		\end{table}
		\section{Numerical results for mDDMt3 in a moving domain in 2D}\label{res_2d}
		
	Table \ref{multigrid} below shows the number of multigrid iterations needed to solve a time-independent version of the problem considered in Sec. \ref{2d} using mDDM3, where the time derivative in the diffusion equation is dropped and the corresponding Poisson problem is solved in the $t=0$ domain. The number of multigrid iterations needed to achieve a tolerance of $10^{-10}$ is independent of minimum grid size $h_{fine}$ and weakly dependent on $\epsilon$. Fewer iterations are required for the time-dependent problem using mDDMt3 (e.g., number of iterations decreases as the time step decreases) but the dependence on $h$ and $\epsilon$ is very similar. The results presented in Sec. \ref{2d} correspond to using $h$ and $\epsilon$ along the diagonal entries in the table, e.g., $(\epsilon,h)=(0.2, 1/32),~(0.1, 1/64),~(0.05, 1/128),~(0.025,1/256),~(0.0125,1/512)$.
		
		\begin{table}[H]
	\center
	\begin{tabular}{c|c|c|c|c|c}
		\bfseries $\epsilon$& $h_{fine}$=1/32 &\bfseries $h_{fine}=1/64$  & $h_{fine}=1/128$ &\bfseries $h_{fine}=1/256$ & $h_{fine}=1/512$ 
		\begin{filecontents*}{num_of_iter.csv}
eps,1/32,1/64,1/128,1/256,1/512
0.2,9,9,9,9,9
0.1,10,10,10,10,10
0.05,12,13,13,13,13
0.025,13,14,15,15,15
0.0125,14,15,16,16,17
\end{filecontents*}
\csvreader{num_of_iter.csv}{}
		{\\\hline \csvcoli & \csvcolii & \csvcoliii & \csvcoliv & \csvcolv& \csvcolvi}
		\\ \hline
	\end{tabular}
	\caption{The number of multigrid iterations needed to achieve a tolerance of $10^{-10}$ for a time-independent problem analogous to that presented in Sec. \ref{2d}. See text. The number of iterations is independent of $h$ and weakly dependent on $\epsilon$.}\label{multigrid}
\end{table}

	Table \ref{err_2dt_mov} below shows the values of the $L^2$ and $L^\infty$ errors and the order of convergence $\mathbf{k}$ with $h_{fine}=\epsilon/6.4$ for mDDMt3 applied to the 2D problem in a moving domain described in Sec. \ref{2d}.
		
		\begin{table}[H]
			\center
			\begin{tabular}{c|c|c|c|c}
				{\bfseries \(\epsilon\)} & {\bfseries \(L^{2}\) error} &\bfseries k  & {\bfseries \(L^{\infty}\) error }&\bfseries k
				\begin{filecontents*}{21.csv}
eps,l^2,k,l^inf,k
0.2,8.89E-02,0.00 ,1.01E-01,0.00 
0.1,2.18E-02,2.03 ,2.91E-02,1.80 
0.05,4.91E-03,2.15 ,8.35E-03,1.80 
0.025,1.13E-03,2.12 ,2.37E-03,1.82 
0.0125,2.90E-04,1.96 ,6.78E-04,1.81 
\end{filecontents*}
\csvreader{21.csv}{}
				{\\\hline \csvcoli & \csvcolii & \csvcoliii & \csvcoliv & \csvcolv}
			\end{tabular}
			\caption{The \(L^2\) and \(L^{\infty}\) errors for simulating the 2D time-dependent diffusion equation at t=0.1 using mDDMt3 on the moving domain \(D(t)\).}\label{err_2dt_mov}
		\end{table}
}

	\end{appendices}

	\bibliography{mybibfile}

\begin{thebibliography}{10}
\expandafter\ifx\csname url\endcsname\relax
  \def\url#1{\texttt{#1}}\fi
\expandafter\ifx\csname urlprefix\endcsname\relax\def\urlprefix{URL }\fi
\expandafter\ifx\csname href\endcsname\relax
  \def\href#1#2{#2} \def\path#1{#1}\fi

\bibitem{000395210500018}
L.~Chen, H.~Wei, M.~Wen, {An interface-fitted mesh generator and virtual
  element methods for elliptic interface problems}, {JOURNAL OF COMPUTATIONAL
  PHYSICS} {334} ({2017}) {327--348}.
\newblock \href {http://dx.doi.org/{10.1016/j.jcp.2017.01.004}}
  {\path{doi:{10.1016/j.jcp.2017.01.004}}}.

\bibitem{Galenko-2018}
P.~Galenko, D.~Alexandrov, E.~Titova, The boundary integral theory for slow and
  rapid curved solid/liquid interfaces propagating into binary systems, Phil.
  Trans. Roy. Soc. A 376~(2113) (2018) 20170218.

\bibitem{Feischl-2015}
M.~Feischl, T.~Fuehrer, N.~Heuer, M.~Karkulik, D.~Praetorius, Adaptive boundary
  element methods a posteriori error estimators, adaptivity, convergence, and
  implementation, Arch. Comp. Meth. Eng. 22~(3) (2015) 309--389.

\bibitem{000402481300001}
T.~Askham, A.~J. Cerfon, {An adaptive fast multipole accelerated Poisson solver
  for complex geometries}, {JOURNAL OF COMPUTATIONAL PHYSICS} {344} ({2017})
  {1--22}.
\newblock \href {http://dx.doi.org/{10.1016/j.jcp.2017.04.063}}
  {\path{doi:{10.1016/j.jcp.2017.04.063}}}.

\bibitem{A1994MX36800005}
R.~GLOWINSKI, T.~PAN, J.~PERIAUX, {A FICTITIOUS DOMAIN METHOD FOR DIRICHLET
  PROBLEM AND APPLICATIONS}, {COMPUTER METHODS IN APPLIED MECHANICS AND
  ENGINEERING} {111}~({3-4}) ({1994}) {283--303}.
\newblock \href {http://dx.doi.org/{10.1016/0045-7825(94)90135-X}}
  {\path{doi:{10.1016/0045-7825(94)90135-X}}}.

\bibitem{000343821300007}
A.~Massing, M.~G. Larson, A.~Logg, M.~E. Rognes, {A Stabilized Nitsche
  Fictitious Domain Method for the Stokes Problem}, {JOURNAL OF SCIENTIFIC
  COMPUTING} {61}~({3}) ({2014}) {604--628}.
\newblock \href {http://dx.doi.org/{10.1007/s10915-014-9838-9}}
  {\path{doi:{10.1007/s10915-014-9838-9}}}.

\bibitem{000326596600021}
S.~Gallier, E.~Lemaire, L.~Lobry, F.~Peters, {A fictitious domain approach for
  the simulation of dense suspensions}, {JOURNAL OF COMPUTATIONAL PHYSICS}
  {256} ({2014}) {367--387}.
\newblock \href {http://dx.doi.org/{10.1016/j.jcp.2013.09.015}}
  {\path{doi:{10.1016/j.jcp.2013.09.015}}}.

\bibitem{000336620900010}
N.~Noormohammadi, B.~Boroomand, {A fictitious domain method using equilibrated
  basis functions for harmonic and bi-harmonic problems in physics}, {JOURNAL
  OF COMPUTATIONAL PHYSICS} {272} ({2014}) {189--217}.
\newblock \href {http://dx.doi.org/{10.1016/j.jcp.2014.04.011}}
  {\path{doi:{10.1016/j.jcp.2014.04.011}}}.

\bibitem{000226822300011}
R.~Mittal, G.~Iaccarino, {Immersed boundary methods}, {ANNUAL REVIEW OF FLUID
  MECHANICS} {37} ({2005}) {239--261}.
\newblock \href {http://dx.doi.org/{10.1146/annurev.fluid.37.061903.175743}}
  {\path{doi:{10.1146/annurev.fluid.37.061903.175743}}}.

\bibitem{000330577100008}
J.~Favier, A.~Revell, A.~Pinelli, {A Lattice Boltzmann-Immersed Boundary method
  to simulate the fluid interaction with moving and slender flexible objects},
  {JOURNAL OF COMPUTATIONAL PHYSICS} {261} ({2014}) {145--161}.
\newblock \href {http://dx.doi.org/{10.1016/j.jcp.2013.12.052}}
  {\path{doi:{10.1016/j.jcp.2013.12.052}}}.

\bibitem{000341311000010}
A.~Calderer, S.~Kang, F.~Sotiropoulos, {Level set immersed boundary method for
  coupled simulation of air/water interaction with complex floating
  structures}, {JOURNAL OF COMPUTATIONAL PHYSICS} {277} ({2014}) {201--227}.
\newblock \href {http://dx.doi.org/{10.1016/j.jcp.2014.08.010}}
  {\path{doi:{10.1016/j.jcp.2014.08.010}}}.

\bibitem{000397072800006}
D.~B. Stein, R.~D. Guy, B.~Thomases, {Immersed Boundary Smooth Extension
  (IBSE): A high-order method for solving incompressible flows in arbitrary
  smooth domains}, {JOURNAL OF COMPUTATIONAL PHYSICS} {335} ({2017})
  {155--178}.
\newblock \href {http://dx.doi.org/{10.1016/j.jcp.2017.01.010}}
  {\path{doi:{10.1016/j.jcp.2017.01.010}}}.

\bibitem{000169338900012}
G.~Tryggvason, B.~Bunner, A.~Esmaeeli, D.~Juric, N.~Al-Rawahi, W.~Tauber,
  J.~Han, S.~Nas, Y.~Jan, {A front-tracking method for the computations of
  multiphase flow}, {JOURNAL OF COMPUTATIONAL PHYSICS} {169}~({2}) ({2001})
  {708--759}.
\newblock \href {http://dx.doi.org/{10.1006/jcph.2001.6726}}
  {\path{doi:{10.1006/jcph.2001.6726}}}.

\bibitem{000398874700008}
M.~Irfan, M.~Muradoglu, {A front tracking method for direct numerical
  simulation of evaporation process in a multiphase system}, {JOURNAL OF
  COMPUTATIONAL PHYSICS} {337} ({2017}) {132--153}.
\newblock \href {http://dx.doi.org/{10.1016/j.jcp.2017.02.036}}
  {\path{doi:{10.1016/j.jcp.2017.02.036}}}.

\bibitem{000412788000084}
H.~Shahin, S.~Mortazavi, {Three-dimensional simulation of microdroplet
  formation in a co-flowing immiscible fluid system using front tracking
  method}, {JOURNAL OF MOLECULAR LIQUIDS} {243} ({2017}) {737--749}.
\newblock \href {http://dx.doi.org/{10.1016/j.molliq.2017.08.082}}
  {\path{doi:{10.1016/j.molliq.2017.08.082}}}.

\bibitem{A1997XQ32700018}
C.~Hirt, A.~Amsden, J.~Cook, {An arbitrary Lagrangian-Eulerian computing method
  for all flow speeds (Reprinted from the Journal of Computational Physics, vol
  14, pg 227-253, 1974)}, {JOURNAL OF COMPUTATIONAL PHYSICS} {135}~({2})
  ({1997}) {203--216}.
\newblock \href {http://dx.doi.org/{10.1006/jcph.1997.5702}}
  {\path{doi:{10.1006/jcph.1997.5702}}}.

\bibitem{000381585100032}
A.~J. Barlow, P.-H. Maire, W.~J. Rider, R.~N. Rieben, M.~J. Shashkov,
  {Arbitrary Lagrangian-Eulerian methods for modeling high-speed compressible
  multimaterial flows}, {JOURNAL OF COMPUTATIONAL PHYSICS} {322} ({2016})
  {603--665}.
\newblock \href {http://dx.doi.org/{10.1016/j.jcp.2016.07.001}}
  {\path{doi:{10.1016/j.jcp.2016.07.001}}}.

\bibitem{000360087400046}
W.~Bo, M.~Shashkov, {Adaptive reconnection-based arbitrary Lagrangian Eulerian
  method}, {JOURNAL OF COMPUTATIONAL PHYSICS} {299} ({2015}) {902--939}.
\newblock \href {http://dx.doi.org/{10.1016/j.jcp.2015.07.032}}
  {\path{doi:{10.1016/j.jcp.2015.07.032}}}.

\bibitem{000414478700020}
C.~Klingenberg, G.~Schnucke, Y.~Xia, {An Arbitrary Lagrangian-Eulerian Local
  Discontinuous Galerkin Method for Hamilton-Jacobi Equations}, {JOURNAL OF
  SCIENTIFIC COMPUTING} {73}~({2-3, SI}) ({2017}) {906--942}.
\newblock \href {http://dx.doi.org/{10.1007/s10915-017-0471-2}}
  {\path{doi:{10.1007/s10915-017-0471-2}}}.

\bibitem{A1994PB02700004}
R.~LEVEQUE, Z.~LI, {THE IMMERSED INTERFACE METHOD FOR ELLIPTIC-EQUATIONS WITH
  DISCONTINUOUS COEFFICIENTS AND SINGULAR SOURCES}, {SIAM JOURNAL ON NUMERICAL
  ANALYSIS} {31}~({4}) ({1994}) {1019--1044}.
\newblock \href {http://dx.doi.org/{10.1137/0731054}}
  {\path{doi:{10.1137/0731054}}}.

\bibitem{000432512900032}
R.~Hu, Z.~Li, {Error analysis of the immersed interface method for Stokes
  equations with an interface}, {APPLIED MATHEMATICS LETTERS} {83} ({2018})
  {207--211}.
\newblock \href {http://dx.doi.org/{10.1016/j.aml.2018.03.034}}
  {\path{doi:{10.1016/j.aml.2018.03.034}}}.

\bibitem{000431399600012}
S.~Amat, Z.~Li, J.~Ruiz, {On an New Algorithm for Function Approximation with
  Full Accuracy in the Presence of Discontinuities Based on the Immersed
  Interface Method}, {JOURNAL OF SCIENTIFIC COMPUTING} {75}~({3}) ({2018})
  {1500--1534}.
\newblock \href {http://dx.doi.org/{10.1007/s10915-017-0596-3}}
  {\path{doi:{10.1007/s10915-017-0596-3}}}.

\bibitem{000432769500031}
Z.~Li, M.-C. Lai, X.~Peng, Z.~Zhang, {A least squares augmented immersed
  interface method for solving Navier-Stokes and Darcy coupling equations},
  {COMPUTERS \& FLUIDS} {167} ({2018}) {384--399}.
\newblock \href {http://dx.doi.org/{10.1016/j.compfluid.2018.03.032}}
  {\path{doi:{10.1016/j.compfluid.2018.03.032}}}.

\bibitem{li2013adaptive}
Z.~Li, P.~Song, Adaptive mesh refinement techniques for the immersed interface
  method applied to flow problems, Computers \& structures 122 (2013) 249--258.

\bibitem{Kolahdouz2020}
E.~Kolahdouz, A.~Pal Singh~Bhalla, C.~B.A., B.~Griffith, An immersed interface
  method for discrete surfaces, J. Comput. Phys. 400 (2020) 108854.

\bibitem{000081366600002}
R.~Fedkiw, T.~Aslam, B.~Merriman, S.~Osher, {A non-oscillatory Eulerian
  approach to interfaces in multimaterial flows (the ghost fluid method)},
  {JOURNAL OF COMPUTATIONAL PHYSICS} {152}~({2}) ({1999}) {457--492}.
\newblock \href {http://dx.doi.org/{10.1006/jcph.1999.6236}}
  {\path{doi:{10.1006/jcph.1999.6236}}}.

\bibitem{000362379300017}
B.~Lalanne, L.~R. Villegas, S.~Tanguy, F.~Risso, {On the computation of viscous
  terms for incompressible two-phase flows with Level Set/Ghost Fluid Method},
  {JOURNAL OF COMPUTATIONAL PHYSICS} {301} ({2015}) {289--307}.
\newblock \href {http://dx.doi.org/{10.1016/j.jcp.2015.08.036}}
  {\path{doi:{10.1016/j.jcp.2015.08.036}}}.

\bibitem{000375799200041}
L.~R. Villegas, R.~Alis, M.~Lepilliez, S.~Tanguy, {A Ghost Fluid/Level Set
  Method for boiling flows and liquid evaporation: Application to the
  Leidenfrost effect}, {JOURNAL OF COMPUTATIONAL PHYSICS} {316} ({2016})
  {789--813}.
\newblock \href {http://dx.doi.org/{10.1016/j.jcp.2016.04.031}}
  {\path{doi:{10.1016/j.jcp.2016.04.031}}}.

\bibitem{000418229800020}
Z.~Ge, J.-C. Loiseau, O.~Tammisola, L.~Brandt, {An efficient mass-preserving
  interface-correction level set/ghost fluid method for droplet suspensions
  under depletion forces}, {JOURNAL OF COMPUTATIONAL PHYSICS} {353} ({2018})
  {435--459}.
\newblock \href {http://dx.doi.org/{10.1016/j.jcp.2017.10.046}}
  {\path{doi:{10.1016/j.jcp.2017.10.046}}}.

\bibitem{000180706400035}
D.~Ingram, D.~Causon, C.~Mingham, {Developments in Cartesian cut cell methods},
  {MATHEMATICS AND COMPUTERS IN SIMULATION} {61}~({3-6}) ({2003}) {561--572},
  {2nd IMACS Conference on Mathematical Modelling and Computational Methods in
  Mechanics, Physics, Biomechanics and Geodynamics, PLZEN, CZECH REPUBLIC, JUN
  19-25, 2001}.
\newblock \href {http://dx.doi.org/{10.1016/S0378-4754(02)00107-6}}
  {\path{doi:{10.1016/S0378-4754(02)00107-6}}}.

\bibitem{000432481000009}
N.~Gokhale, N.~Nikiforakis, R.~Klein, {A dimensionally split Cartesian cut cell
  method for hyperbolic conservation laws}, {JOURNAL OF COMPUTATIONAL PHYSICS}
  {364} ({2018}) {186--208}.
\newblock \href {http://dx.doi.org/{10.1016/j.jcp.2018.03.005}}
  {\path{doi:{10.1016/j.jcp.2018.03.005}}}.

\bibitem{000428483000007}
F.~Nikfarjam, Y.~Cheny, O.~Botella, {The LS-STAG immersed boundary/cut-cell
  method for non-Newtonian flows in 3D extruded geometries}, {COMPUTER PHYSICS
  COMMUNICATIONS} {226} ({2018}) {67--80}.
\newblock \href {http://dx.doi.org/{10.1016/j.cpc.2018.01.006}}
  {\path{doi:{10.1016/j.cpc.2018.01.006}}}.

\bibitem{000427393800011}
B.~Muralidharan, S.~Menon, {Simulation of moving boundaries interacting with
  compressible reacting flows using a second-order adaptive Cartesian cut-cell
  method}, {JOURNAL OF COMPUTATIONAL PHYSICS} {357} ({2018}) {230--262}.
\newblock \href {http://dx.doi.org/{10.1016/j.jcp.2017.12.030}}
  {\path{doi:{10.1016/j.jcp.2017.12.030}}}.

\bibitem{000306367500004}
R.~I. Saye, J.~A. Sethian, {Analysis and applications of the Voronoi Implicit
  Interface Method}, {JOURNAL OF COMPUTATIONAL PHYSICS} {231}~({18}) ({2012})
  {6051--6085}.
\newblock \href {http://dx.doi.org/{10.1016/j.jcp.2012.04.004}}
  {\path{doi:{10.1016/j.jcp.2012.04.004}}}.

\bibitem{000358796700042}
A.~Guittet, M.~Lepilliez, S.~Tanguy, F.~Gibou, {Solving elliptic problems with
  discontinuities on irregular domains - the Voronoi Interface Method},
  {JOURNAL OF COMPUTATIONAL PHYSICS} {298} ({2015}) {747--765}.
\newblock \href {http://dx.doi.org/{10.1016/j.jcp.2015.06.026}}
  {\path{doi:{10.1016/j.jcp.2015.06.026}}}.

\bibitem{000283202200001}
T.-P. Fries, T.~Belytschko, {The extended/generalized finite element method: An
  overview of the method and its applications}, {INTERNATIONAL JOURNAL FOR
  NUMERICAL METHODS IN ENGINEERING} {84}~({3}) ({2010}) {253--304}.
\newblock \href {http://dx.doi.org/{10.1002/nme.2914}}
  {\path{doi:{10.1002/nme.2914}}}.

\bibitem{000266373500001}
T.~Belytschko, R.~Gracie, G.~Ventura, {A review of extended/generalized finite
  element methods for material modeling}, {MODELLING AND SIMULATION IN
  MATERIALS SCIENCE AND ENGINEERING} {17}~({4}).
\newblock \href {http://dx.doi.org/{10.1088/0965-0393/17/4/043001}}
  {\path{doi:{10.1088/0965-0393/17/4/043001}}}.

\bibitem{000402217000004}
J.~Carlos~Martinez, L.~V. Vanegas~Useche, M.~A. Wahab, {Numerical prediction of
  fretting fatigue crack trajectory in a railway axle using XFEM},
  {INTERNATIONAL JOURNAL OF FATIGUE} {100}~({1}) ({2017}) {32--49}.
\newblock \href {http://dx.doi.org/{10.1016/j.ijfatigue.2017.03.009}}
  {\path{doi:{10.1016/j.ijfatigue.2017.03.009}}}.

\bibitem{bedrossian2010second}
J.~Bedrossian, J.~H. Von~Brecht, S.~Zhu, E.~Sifakis, J.~M. Teran, A second
  order virtual node method for elliptic problems with interfaces and irregular
  domains, Journal of Computational Physics 229~(18) (2010) 6405--6426.

\bibitem{HELLRUNG20122015}
J.~L. Hellrung, L.~Wang, E.~Sifakis, J.~M. Teran,
  \href{http://www.sciencedirect.com/science/article/pii/S0021999111006784}{A
  second order virtual node method for elliptic problems with interfaces and
  irregular domains in three dimensions}, Journal of Computational Physics
  231~(4) (2012) 2015 -- 2048.
\newblock \href {http://dx.doi.org/https://doi.org/10.1016/j.jcp.2011.11.023}
  {\path{doi:https://doi.org/10.1016/j.jcp.2011.11.023}}.
\newline\urlprefix\url{http://www.sciencedirect.com/science/article/pii/S0021999111006784}

\bibitem{shirokoff2015sharp}
D.~Shirokoff, J.-C. Nave, A sharp-interface active penalty method for the
  incompressible navier--stokes equations, Journal of Scientific Computing
  62~(1) (2015) 53--77.

\bibitem{Kockelkoren-2003}
J.~Kockelkoren, H.~Levine, W.-J. Rappel, Computational approach for modeling
  intra- and computational approach for modeling intra- and extracellular
  dynamics, Phys. Rev. E 68 (2003) 037702.

\bibitem{Bueno-Orovio-2005}
A.~Bueno-Orovio, V.~Perez-Garcia, Spectral smoothed boundary methods: The role
  of external boundary conditions, Numer. Meth. Partial Diff. Eqns. 22 (2005)
  435--448.

\bibitem{Bueno-Orovio-2006}
A.~Bueno-Orovio, V.~Perez-Garcia, F.~Fenton, Spectral methods for partial
  differential equations in irregular domains: The spectral smoothed boundary
  method, SIAM J. Sci. Comput. 28 (2006) 886--900.

\bibitem{Ratz-2006}
A.~Raetz, A.~Voigt, Pde's on surfaces --- a diffuse interface approach, Comm.
  Math. Sci. 4~(3) (2006) 575--590.

\bibitem{li09}
X.~Li, J.~Lowengrub, A.~Raetz, A.~Voigt, {Solving PDEs in complex geometries: a
  diffuse domain approach}, {COMMUNICATIONS IN MATHEMATICAL SCIENCES} {7}~({1})
  ({2009}) {81--107}.

\bibitem{tei09}
K.~E. Teigen, X.~Li, J.~Lowengrub, F.~Wang, A.~Voigt, A diffuse-interface
  approach for modeling transport, diffusion and adsorption/desorption of
  material quantities on a deformable interface, {COMMUNICATIONS IN
  MATHEMATICAL SCIENCES} {7}~({4}) ({2009}) {1009--1037}.

\bibitem{Yu-2012}
H.-C. Yu, H.-Y. Chen, K.~Thornton, Extended smoothed boundary method for
  solving partial differential equations with general boundary conditions on
  complex boundaries, {MODELLING AND SIMULATION IN MATERIALS SCIENCE AND
  ENGINEERING} 20 (2012) 075008.

\bibitem{Poulson-2018}
S.~Poulsen, P.~Voorhees, Smoothed boundary method for diffusion-related partial
  differential equations in complex geometries, Int. J. Comput. Meth. 15~(3)
  (2018) 1850014.

\bibitem{Fenton-2005}
F.~Fenton, E.~Cherry, A.~Karma, W.-J. Rappel, Modeling wave propagation in
  realistic heart geometries using the phase-field method, Chaos 15 (2005)
  013502.

\bibitem{Aland-2011a}
S.~Aland, C.~Landsberg, R.~Mueller, F.~Stenger, M.~Bobeth, A.~Langheinrich,
  A.~Voigt, Adaptive diffuse domain approach for calculating mechanically
  induced deformation of trabecular bone, Comp. Meth. Biomech. Biomed. Eng.
  17~(1) (2011) 31--38.

\bibitem{Camley-2013}
B.~Camley, Y.~Zhao, B.~Li, H.~Levine, W.-J. Rappel, Periodic migration in a
  physical model of cells on micropatterns, Phys. Rev. Lett. 111 (2013) 158102.

\bibitem{Chen-2014}
Y.~Chen, J.~Lowengrub, Tumor growth in complex, evolving microenvironmental
  geometries: A diffuse domain approach, J. Theor. Biol. 361 (2014) 14--30.

\bibitem{Ratz-2014}
A.~Raetz, M.~Roeger, Symmetry breaking in a bulk-surface reaction-diffusion
  model for signaling networks, Nonlinearity 27 (2014) 1805--1827.

\bibitem{Ratz-2015}
A.~Raetz, A new diffuse-interface model for step flow in epitaxial growth, IMA
  J. Appl. Math. 80~(3) (2015) 697--711.

\bibitem{Lowengrub-2016}
J.~Lowengrub, J.~Allard, S.~Aland, Numerical simulation of endocytosis: Viscous
  flow driven by membranes with non-uniformly distributed curvature-inducing
  molecules, J. Comput. Phys. 309 (2016) 112--128.

\bibitem{Camley-2017}
B.~Camley, Y.~Zhao, B.~Li, H.~Levine, W.-J. Rappel, Crawling and turning in a
  minimal reaction-diffusion cell motility model: Coupling cell shape and
  biochemistry, Phys. Rev. E 95 (2017) 012401.

\bibitem{Lipkova2019}
J.~Lipkova, Computational modelling in neuro-onconolgy, Ph.D. Thesis, T.U.
  Munich.

\bibitem{teigen2011diffuse}
K.~E. Teigen, P.~Song, J.~Lowengrub, A.~Voigt, {A diffuse-interface method for
  two-phase flows with soluble surfactants}, {JOURNAL OF COMPUTATIONAL PHYSICS}
  {230}~({2}) ({2011}) {375--393}.
\newblock \href {http://dx.doi.org/{10.1016/j.jcp.2010.09.020}}
  {\path{doi:{10.1016/j.jcp.2010.09.020}}}.

\bibitem{ala10}
S.~Aland, J.~Lowengrub, A.~Voigt, {Two-phase flow in complex geometries: A
  diffuse domain approach}, {CMES-COMPUTER MODELING IN ENGINEERING \& SCIENCES}
  {57}~({1}) ({2010}) {77--107}.

\bibitem{Aland-2011}
S.~Aland, J.~Lowengrub, A.~Voigt, A continuum model of colloid-stabilized
  interfaces, Phys. Fluids 23~(6) (2011) 062103.

\bibitem{Aland-2012}
S.~Aland, J.~Lowengrub, A.~Voigt, Particles at fluid-fluid interfaces: A new
  navier-stokes-cahn-hilliard surface-phase-field-crystal model, Phys. Rev. E
  86~(4) (2012) 046321.

\bibitem{Aland-2014}
S.~Aland, S.~Egerer, J.~Lowengrub, A.~Voigt, Diffuse interface models of
  locally inextensible vesicles in a viscous fluid, J. Comput. Phys. 277 (2014)
  32--47.

\bibitem{Yu-2016}
H.-C. Yu, M.-J. Choe, G.~Amatucci, Y.-M. Chiang, K.~Thornton, Smoothed boundary
  method for simulating bulk and grain boundary transport in complex
  polycrystalline microstructures, Comput. Mater. Sci. 121 (2016) 14--22.

\bibitem{Ratz-2016}
A.~Raetz, Diffuse-interface approximations of osmosis free boundary problems
  diffuse-interface approximations of osmosis free boundary problems, SIAM J.
  Appl. Math. 76~(3) (2016) 910--929.

\bibitem{Hong-2016}
L.~Hong, L.~Liang, S.~Bhattacharyya, W.~Xing, L.-Q. Chen, Anisotropic li
  intercalation in a lixfepo4 nano-particle: a spectral smoothed boundary
  phase-field model, Phys, Chem. Chem. Phys. 18 (2016) 9537--9543.

\bibitem{Chadwick-2018}
A.~Chadwick, J.~Steward, R.~Du, K.~Thornton, Numerical modeling of localized
  corrosion using phase-field and smoothed boundary methods, J. Electrochemical
  Soc. 165~(10) (2018) C633--C646.

\bibitem{Vey2007}
S.~Vey, A.~Voigt, Amdis: Adaptive multidimensional simulations, Comput. Vis.
  Sci. 10 (2007) 57--67.

\bibitem{Rossinelli2015}
D.~Rossinelli, B.~Hejazialhosseini, W.~van Rees, M.~Gazzola, M.~Bergdorf,
  P.~Koumoutsakos, Mrag-i2d: Multi-resolution adapted grids for remeshed vortex
  methods on multicore architectures, J. Comput. Phys. 288 (2015) 1--18.

\bibitem{Feng2018}
W.~Feng, Z.~Guo, J.~Lowengrub, S.~Wise, Mass-conservative adaptive fas
  multigrid solver for cell-centered finite difference methods on
  block-structured, locally-cartesian grids, J. Comput. Phys. 352 (2018)
  463--497.

\bibitem{fra12}
S.~Franz, H.-G. Roos, R.~G{\"a}rtner, A.~Voigt, A note on the convergence
  analysis of a diffuse-domain approach, {COMPUTATIONAL METHODS IN APPLIED
  MATHEMATICS} 12~(2) (2012) 153--167.

\bibitem{Abels-2015}
H.~Abels, K.~Lam, B.~Stinner, Analysis of the diffuse domain approach for a
  bulk-surface coupled pde system, SIAM J. Math. Anal. 47~(5) (2015)
  3687--3725.

\bibitem{Schlottbom-2016}
M.~Schlottbom, Error analysis of a diffuse interface method for elliptic
  problems with dirichlet boundary conditions, Appl. Num. Math. 109 (2016)
  109--122.

\bibitem{Burger-2017}
M.~Burger, O.~Elvetun, M.~Schlottbom, Analysis of the diffuse domain method for
  second order elliptic boundary value problems, Found. Comput. Math. 17 (2017)
  627--674.

\bibitem{000356209300006}
K.~Y. Lerv{\aa}g, J.~Lowengrub, {Analysis of the diffuse-domain method for
  solving PDEs in complex geometries}, {COMMUNICATIONS IN MATHEMATICAL
  SCIENCES} {13}~({6}) ({2015}) {1473--1500}.
\newblock \href {http://dx.doi.org/{10.4310/CMS.2015.v13.n6.a6}}
  {\path{doi:{10.4310/CMS.2015.v13.n6.a6}}}.

\bibitem{Osher-1988}
S.~Osher, J.~A. Sethian, Fronts propagating with curvature-dependent
  speed-algorithms based on hamilton-jacobi formulations, J. Comput. Phys.
  79~(1) (1988) 12--49.

\bibitem{Gibou-2018}
F.~Gibou, R.~Fedkiw, S.~Osher, A review of level-set methods and some recent
  applications, J. Comput. Phys. 353 (2018) 82--109.

\bibitem{hoffman2001numerical}
J.~D. Hoffman, S.~Frankel, Numerical methods for engineers and scientists, CRC
  press, 2001.

\bibitem{bender2013advanced}
C.~M. Bender, S.~A. Orszag, Advanced mathematical methods for scientists and
  engineers I: Asymptotic methods and perturbation theory, Springer Science \&
  Business Media, 2013.

\bibitem{shu11}
G.~Jiang, C.~Shu, {Efficient implementation of weighted ENO schemes}, {JOURNAL
  OF COMPUTATIONAL PHYSICS} {126}~({1}) ({1996}) {202--228}.
\newblock \href {http://dx.doi.org/{10.1006/jcph.1996.0130}}
  {\path{doi:{10.1006/jcph.1996.0130}}}.

\bibitem{gottlieb1998total}
S.~Gottlieb, C.~Shu, {Total variation diminishing Runge-Kutta schemes},
  {MATHEMATICS OF COMPUTATION} {67}~({221}) ({1998}) {73--85}.
\newblock \href {http://dx.doi.org/{10.1090/S0025-5718-98-00913-2}}
  {\path{doi:{10.1090/S0025-5718-98-00913-2}}}.

\bibitem{min10}
C.~Min, {On reinitializing level set functions}, {JOURNAL OF COMPUTATIONAL
  PHYSICS} {229}~({8}) ({2010}) {2764--2772}.
\newblock \href {http://dx.doi.org/{10.1016/j.jcp.2009.12.032}}
  {\path{doi:{10.1016/j.jcp.2009.12.032}}}.

\bibitem{rus00}
G.~Russo, P.~Smereka, {A remark on computing distance functions}, {JOURNAL OF
  COMPUTATIONAL PHYSICS} {163}~({1}) ({2000}) {51--67}.
\newblock \href {http://dx.doi.org/{10.1006/jcph.2000.6553}}
  {\path{doi:{10.1006/jcph.2000.6553}}}.

\bibitem{peng1999pde}
D.~Peng, B.~Merriman, S.~Osher, H.~Zhao, M.~Kang, {A PDE-based fast local level
  set method}, {JOURNAL OF COMPUTATIONAL PHYSICS} {155}~({2}) ({1999})
  {410--438}.
\newblock \href {http://dx.doi.org/{10.1006/jcph.1999.6345}}
  {\path{doi:{10.1006/jcph.1999.6345}}}.

\bibitem{wis07}
S.~Wise, J.~Kim, J.~Lowengrub, {Solving the regularized, strongly anisotropic
  Cahn-Hilliard equation by an adaptive nonlinear multigrid method}, {JOURNAL
  OF COMPUTATIONAL PHYSICS} {226}~({1}) ({2007}) {414--446}.
\newblock \href {http://dx.doi.org/{10.1016/j.jcp.2007.04.020}}
  {\path{doi:{10.1016/j.jcp.2007.04.020}}}.

\bibitem{Atlas}
Icbm 152 nonlinear atlases version 2009, bic.mni.mcgill.ca.

\bibitem{Collins1999}
D.~Collins, A.~Zijdenbos, W.~Baare, A.~Evans, Animal+insect: Improved cortical
  structure segmentation, IPMI Lect. Notes Comput. Sci. 1613/1999 (1999)
  210--223.

\bibitem{Fonov2009}
V.~Fonov, A.~Evans, K.~Botteron, R.~McKinstry, C.~Almli, D.~Collins, Unbiased
  nonlinear average age-appropriate brain templates from birth to adulthood,
  NeuroImage 47 (2009) S102, organization for Human Brain Mapping 2009 Annual
  Meeting.

\bibitem{Fonov2011}
V.~Fonov, A.~Evans, K.~Botteron, C.~Almli, R.~McKinstry, D.~Collins, BDCG,
  Unbiased nonlinear average age-appropriate atlases for pediatric studies,
  Comput. Vis. Sci. 54 (2011) 1053--8119.

\bibitem{Lervag2015}
K.~Lervag, J.~Lowengrub, Analysis of the diffuse-domain method for solving pdes
  in complex geometries, Commun. Math. Sci. 6 (2015) 1473--1500.

\end{thebibliography}
	
\end{document}